\documentclass[12pt,leqno]{amsart}
\usepackage{amsmath,amscd,amsthm,amsxtra,amssymb}
\usepackage{epsfig,graphics} 
\usepackage[usenames,dvipsnames]{color}

\usepackage{amssymb,latexsym}
\usepackage{mathrsfs}
\usepackage{comment}
\usepackage[poly,all]{xy}
\usepackage{marginnote}
\usepackage{xspace}

\usepackage[colorlinks=true, pdfstartview=FitV, linkcolor=blue,citecolor=blue,urlcolor=blue]{hyperref}

\usepackage{yfonts}
\usepackage{enumerate}

\usepackage[usenames,dvipsnames,svgnames,table]{xcolor}

\usepackage[normalem]{ulem}  

\usepackage[colorinlistoftodos]{todonotes}

\newdir{ >}{{}*!/-10pt/@{>}}

\allowdisplaybreaks[3]

\newcommand{\arxiv}[1]{\href{http://arxiv.org/abs/#1}{\texttt{arXiv:#1}}}

\newlength{\mylength}
\renewcommand{\le}{\leqslant}
\renewcommand{\ge}{\geqslant}

\renewcommand{\geq}{\ge}
\renewcommand{\leq}{\le}

\theoremstyle{plain}
\newtheorem{thm}{Theorem}[section]
\newtheorem{mthm}[thm]{Main\,Theorem}
\newtheorem{prop}[thm]{Proposition}
\newtheorem{coro}[thm]{Corollary}
\newtheorem{lem}[thm]{Lemma}

\newtheorem{conj}{Conjecture}

\theoremstyle{definition}
\newtheorem{ex}[thm]{Example}
\newtheorem{remark}[thm]{Remark}
\newtheorem{definition}[thm]{Definition}
\newtheorem*{convention}{Convention}
\newtheorem{question}[conj]{Qusetion}
\newtheorem{problem}[conj]{Problem}

\newcommand{\nc}{\newcommand}

\newenvironment{answer}
{\noindent{\bf Answer}\hs{1ex}}
{\hfill \qedsymbol}
\nc{\Prop}{\begin{prop}}
\nc{\enprop}{\end{prop}}
\nc{\Lemma}{\begin{lem}}
\nc{\enlemma}{\end{lem}}
\nc{\Ex}{\begin{ex}}
\nc{\enex}{\end{ex}}
\nc{\Th}{\begin{thm}}
\nc{\enth}{\end{thm}}
\nc{\Def}{\begin{definition}}
\nc{\edf}{\end{definition}}
\nc{\Conj}{\begin{conj}}
\nc{\enconj}{\end{conj}}
\nc{\Quest}{\begin{question}}
\nc{\enquest}{\end{question}}
\nc{\Rem}{\begin{remark}}
\nc{\enrem}{\end{remark}}
\nc{\Ans}{\begin{answer}}
\nc{\enans}{\end{answer}}
\nc{\Prob}{\begin{problem}}
\nc{\enprob}{\end{problem}}
\nc{\MTh}{\begin{mthm}}
\nc{\enmth}{\end{mthm}}

\newenvironment{red}
{\relax\color{red}}
{\hspace*{.5ex}\relax}

\nc{\berm}{\ber{}\marginnote{\fbox{\scshape\lowercase{M}}}}

\newenvironment{jaune}{\relax\marginnote{\ber \scalebox{.6}{\sc{To be deleted}} \er}
  \color{Orchid}}{\hspace*{.5ex}\relax}

\nc{\bj}{\begin{jaune}}
\nc{\ej}{\end{jaune}}

\newenvironment{magenta}
{\relax\color{magenta}}
{\hspace*{.5ex}\relax}
\newcommand{\bemg}{\begin{magenta}}
\newcommand{\emg}{\end{magenta}}

\nc{\bey}{\begin{jaune}}
  \nc{\ey}{\end{jaune}}

\nc{\on}{\operatorname}

\newcommand{\Q}{\mathbb {Q}}
\newcommand{\Z}{\ms{2mu}{\mathbb Z}}
\newcommand{\R}{{\mathbb R}}

\newcommand{\A}{{\mathbf A}}

\newcommand{\D}{\mathscr{D}\ms{1mu}}

\newcommand{\one}{{\bf{1}}}
\newcommand{\seteq}{\mathbin{:=}}

\newcommand{\hd}{{\mathrm{hd}}}      					 
\newcommand{\To}[1][{\hs{0.8ex}}]{\xrightarrow{\ms{7mu}{#1}\ms{7mu}}}

\newcommand{\Sp}{\mathrm{span}_{\R_{\ge0}}}  	

\newcommand{\g}{\ms{1mu}\mathfrak{g}\ms{1mu}}
\newcommand{\n}{\mathfrak{n}}

\newcommand{\Hom}{\operatorname{Hom}}
\newcommand{\HOM}{\on{\mathrm{H{\scriptstyle OM}}}}
\newcommand{\END}{\on{\mathrm{E\scriptstyle ND}}\ms{.1mu}}

\nc{\Rnor}{\mathrm{R}^{\mathrm{norm}}}

\newcommand{\isoto}[1][]{\mathop{\xrightarrow%
[{\raisebox{.3ex}[0ex][.3ex]{$\scriptstyle{#1}$}}]%
{{\raisebox{-.6ex}[0ex][-.6ex]{$\mspace{5mu}\sim\mspace{5mu}$}}}}}

\newcommand{\Mod}{\on{Mod}}
\newcommand{\gmod}{\text{-}\mathrm{gmod}}
\newcommand{\gMod}{\text{-}\mathrm{gMod}}
\newcommand{\proj}{\text{-}\mathrm{proj}}

\nc{\F}{\mathrm{F}}

\newcommand{\conv}[1][]{
\underset{\raisebox{.5ex}{$\scriptstyle{#1}$}}{\mathbin{\scalebox{1.1}{$\mspace{1.5mu}\circ\mspace{1.5mu}$}}}}
\newcommand{\hconv}{\mathbin{\scalebox{.9}{$\nabla$}}}

\newcommand{\sconv}{\mathbin{\scalebox{.9}{$\Delta$}}}

\renewcommand{\Im}{\on{Im}}
\newcommand{\de}{\on{\textfrak{d}}}

\newcommand{\cmA}{\cartan}  
\newcommand{\wlP}{\mathsf{P}}   
\newcommand{\rlQ}{\mathsf{Q}}   
\newcommand{\weyl}{\mathsf{W}}  
\newcommand{\prD}{\Delta_+}            
\newcommand{\nrD}{\Delta_-}            
\newcommand{\sg}{\mathfrak{S}}   
\newcommand{\Po}{\wlP}
\newcommand{\rtlp}{\rtl_+}
\nc{\prt}{\prD}
\nc{\qQ}{Q}
\newcommand{\bQ}{\overline{\qQ}}

\newcommand\Aq[1][{\g^+}]{A_q(#1)}

\newcommand{\wt}{\mathrm{wt}\ms{1.5mu}} 		
\newcommand{\bR}{\mathbf{k}} 		
\nc{\corp}{\bR}
\newcommand{\catC}{ \mathscr{C}}  	
\newcommand{\tcatC}{ \widetilde{\mathscr{C}}}  	
\newcommand{\catT}{ \mathcal{T}}  	
\newcommand{\lT}{ \widetilde{\mathcal{T}}}  	

\newcommand{\dM}{ \mathsf{M }}              
\newcommand{\dC}{ \mathsf{C }}              
\newcommand{\gW}{\mathsf{W}}
\newcommand{\sgW}{\mathsf{W}^*}
\newcommand{\tf}{{\widetilde{f}}}
\newcommand{\te}{{\widetilde{e}\ms{1.5mu}}}  		
\newcommand{\tF}{\widetilde{\mathrm{F}}} 		
\newcommand{\tE}{\widetilde{\mathrm{E}}} 		
\newcommand{\tEs}{\widetilde{\mathrm{E}}^*} 		
\newcommand{\tFs}{\widetilde{\mathrm{F}}^*} 		
\newcommand{\tEm}{\widetilde{\mathrm{E}}^{\hskip 0.1em \rm max}}  		
\newcommand{\tEsm}{\widetilde{\mathrm{E}}^{*{ \hskip 0.1em  \rm max}}}  		
\newcommand{\ep}{\varepsilon}  		
\newcommand{\ph}{\varphi}  		
\newcommand{\trivialM}{\mathbf{1}} 	

\newcommand{\Ht}{\mathrm{ht}} 		
\newcommand{\coRl}{\mathrm{R}^{{\ell}}} 

\newcommand{\La}{\Lambda} 			
\newcommand{\tLa}{\widetilde{\Lambda}} 			
\newcommand{\Res}{\mathrm{Res}\ms{1mu}} 			

\nc{\Ma}{{\ms{1.5mu}\mathsf{M}}}
\nc{\Na}{\mathsf{N}}
\nc{\Xa}{\mathsf{X}}
\nc{\Ya}{\mathsf{Y}}
\nc{\Laa}{\mathsf{L}}
\newcommand{\Mm}{{\mathsf{M}}}
\newcommand{\z}[1][{\Ma}]{{z_{#1}}}

\newcommand{\zM}{{z_\Ma}}
\newcommand{\zN}{{z_\Na}}

\newcommand{\triv}{{\mathbf{1}}}   				
\newcommand{\id}{\ms{2mu}{\mathsf{id}}\ms{1mu}}   				
 
\newcommand{\dphi}{{\phi}}   				

\nc{\be}{\begin{enumerate}}
\newcommand{\bnum}{\be[{\rm(i)}]}
\newcommand{\bna}{\be[{\rm(a)}]}

\newcommand{\rtl}{\rlQ}

\newcommand{\etens}{\boxtimes}

\newcommand{\rmat}[1]{\ms{1mu}{\mathbf{r}}_%
{\mspace{-2mu}\raisebox{-.6ex}{${\scriptstyle{#1}}$}}}

\newcommand{\shc}{{\ms{2mu}\mathcal{C}}}

\newcommand{\tC}{\widetilde{\catC}}

\newcommand{\Ob}{\on{Ob}}

\newcommand{\Irr}{\mathrm{Irr}}

\nc{\ms}{\mspace}
\nc{\cl}{\colon}
\nc{\ro}{{\rm (}}
\nc{\rf}{{\rm )}\xspace}
\nc{\noi}{\noindent}
\nc{\bl}{\bigl(}
\nc{\br}{\bigr)}

\newcommand{\cA}{{\mathcal A}}

\newcommand{\bbZ}{\mathbb Z}

\nc{\ify}{\infty}
\nc{\uqm}{U^-_q(\g)}
\nc{\uq}{U_q(\g)}
\nc{\til}{\tilde}
\nc{\wtil}{\widetilde}
\nc{\vep}{\varepsilon}
\nc{\vp}{\varphi}
\nc{\eit}{\tilde{e}_i}
\nc{\fit}{\tilde{f}_i}
\nc{\ot}{\otimes}

\def\del{\delta}
\def\mapright#1{\smash{\mathop{\longrightarrow}\limits^{#1}}}
\def\TY(#1,#2,#3){#1^{(#2)}_{#3}}
\def\qq{\qquad}
\def\q{\quad}

\newenvironment{myequationn}
{\relax\setlength{\arraycolsep}{1pt}\begin{eqnarray*}}
{\end{eqnarray*}}

\newenvironment{myequation}
{\relax\setlength{\arraycolsep}{1pt}\begin{eqnarray}}
{\end{eqnarray}}

\nc{\eq}{\begin{myequation}}
\nc{\eneq}{\end{myequation}}
\nc{\eqn}{\begin{myequationn}}
\nc{\eneqn}{\end{myequationn}}

\newenvironment{myarray}[1]{\relax\setlength{\arraycolsep}{1pt}
\begin{array}{#1}}{\end{array}\relax}

\newcommand{\ba}{\begin{myarray}}
\newcommand{\ea}{\end{myarray}}

\nc{\hs}{\hspace*}
\nc{\vs}{\vspace*}
\nc{\set}[2]{\left\{{#1}\mid{#2}\right\}}
\nc{\snoi}{\smallskip\noi}
\nc{\mnoi}{\medskip\noi}
\nc{\al}{\alpha}
\nc{\rmz}{\setminus\{0\}}
\nc{\tens}
[1][]{\mathbin{\otimes{\ms{2mu}}_{\raise1.5ex\hbox to-.1em{}#1}}}
\nc{\vphi}{\varphi}
\nc{\ee}{\end{enumerate}}
\nc{\la}{\lambda}
\nc{\bc}{\begin{cases}}
\nc{\cct}{\mathbin{\scalebox{1.2}{$\star$}}}
\nc{\ec}{\end{cases}}
\nc{\qtq}[1][and]{\quad\text{#1}\quad}
\nc{\qt}[1]{\quad\text{#1}}
\nc{\dual}{{\displaystyle{\ms{1mu}\star}}}
\nc{\wle}{\preceq}
\nc{\epito}{\twoheadrightarrow}
\nc{\epiTo}[1][]{\xymatrix@C=4ex{{}\ar@{->>}[r]^-{#1}&{}}}
\nc{\Proof}{\begin{proof}}
\nc{\lan}{\langle}
\nc{\ran}{\rangle}
\nc{\ang}[1]{\lan{#1}\ran}
\nc{\QED}{\end{proof}}
\nc{\soplus}{\mathbin{\raisebox{.1ex}{\scalebox{.65}{\raisebox{.4ex}{$\displaystyle\bigoplus$}}}}}
\nc{\eps}{\varepsilon}
\nc{\supp}{\on{supp}}
\nc{\sct}{strongly commute\xspace}
\nc{\scts}{strongly commutes\xspace}
\nc{\bce}{\eta}			
\nc{\height}[1]{\on{ht}(\ms{.5mu}{#1}\ms{.5mu})}
\nc{\braid}{{\ms{1mu}\mathrm{br}}}
\nc{\gp}{\mathfrak{p}}
\nc{\wtl}{\wlP}
\nc{\ra}{real and admits an affinization}
\nc{\ras}{real and admit affinizations}
\nc{\Cor}{\begin{coro}}
\nc{\encor}{\end{coro}}
\nc{\shf}{\mathcal{F}}
\nc{\Cw}[1][{w}]{\catC_{{#1}}}
\nc{\tCw}[1][{w}]{\widetilde{\catC}_{{#1}}} 
\nc{\tCwv}[1][{w,v}]{\widetilde{\catC}_{{#1}}} 
\nc{\akew}[1][1ex]{\rule[-1ex]{#1}{0ex}}
\nc{\ake}[1][2ex]{\rule[-1ex]{0ex}{#1}}
\nc{\akete}[1][-1ex]{\rule[{#1}]{0ex}{1ex}}
\nc{\tRm}{(R\gmod)\widetilde{\mbox{$\ake[2.5ex]\akew[.9ex]$}}}
\nc{\monoTo}[1][]{\xymatrix{\ar@{>->}[r]^-{{#1}}&}}
\nc{\monoto}[1][]{\rightarrowtail}
\nc{\tX}{\widetilde{X}}
\nc{\corps}{\corp}
\nc{\tL}{\widetilde{L}}
\nc{\prtl}{\rtl_+}
\nc{\nrtl}{\rtl_-}
\nc{\nrt}{\rt_-}
\nc{\tK}{\widetilde{K}}
\nc{\tep}{\widetilde\ep}
\nc{\teps}{\widetilde\ep}
\nc{\tmu}{\widetilde \mu} 
\nc{\teta}{\widetilde\eta}
\nc{\ga}{\mathfrak{a}}
\nc{\scbul}{{\,\raise1pt\hbox{$\scriptscriptstyle\bullet$}\,}}
\nc{\bwr}{\mbox{\large$\wr$}}
\nc{\tR}{{\widetilde{\mathrm{R}}}}
\nc{\lS}{\mathsf{S}}
\nc{\lZ}{\mathcal{Z}}
\nc{\prolim}[1][]{\mathop{\varprojlim}\limits_{{#1}}}
\nc{\Sym}{\mathfrak{S}}
\nc{\sym}{\sg}
\newcounter{myc}
\newcounter{mycc}
\nc{\tCs}[1][w]{ \widetilde{\mathscr{C}}^{\ms{2mu}*}_{#1}}
\nc{\txi}{\tilde{\xi}}
\nc{\rl}{\rlQ}
\nc{\sfC}{\mathsf{C}}
\nc{\cor}{{\ms{1mu}\mathbf{k}\ms{1mu}}}
\nc{\Pro}{\on{Pro}}
\newcommand{\proolim}[1][]{\mathop{\text{``}\ms{2mu}\varprojlim\ms{-5mu}\text{''}}\limits_{#1}}
\nc{\hM}{\widehat{\mathsf{M}}}
\nc{\aff}{\mathrm{aff}}
\nc{\rDa}{{\mathscr{D}_\aff}}
\nc{\st}[1]{\left\{{#1}\right\}}
\nc{\W}{\mathsf{W}}
\nc{\rt}{\Delta}
\nc{\pwtl}{\wtl_+}
\nc{\rev}{{\mathrm{rev}}}
\nc{\E}{\mathrm{E}}

\nc{\Qt}[1][w]{\mathrm{Q}_{#1}}
\nc{\Ctr}{\mathsf{C}}
\nc{\Ctrs}{{\mathsf{C}^*}}
\nc{\Dynkin}{\Delta}
\nc{\cartan}{\mathsf{C}}
\nc{\sfc}{\mathsf{c}}
\nc{\sfa}{\mathsf{a}}
\nc{\SW}{\mathrm{K}}
\nc{\hSW}{\widehat{\mathrm{K}}}
\nc{\refl}{\mathscr{S}}
\nc{\Rre}{\mathrm{R}^{\mathrm{ren}}}
\nc{\Rpre}{\mathrm{R}^{\mathrm{ren}\;'}}
\nc{\bRre}{\ol{\mathrm{R}}^{\ms{2mu}\mathrm{ren}}}
\nc{\sfd}{\ms{1mu}\mathsf{d}\ms{1mu}}
\nc{\shm}{\mathcal{M}}
\nc{\sht}{\mathcal{T}}
\nc{\rank}{\mathrm{rank}}
\nc{\Da}{{\D}\ms{-2.8mu}\raisebox{-.35ex}{$\scriptstyle\mathrm{aff}$}}
\nc{\lDa}{\Da^{-1}}
\nc{\bchi}{{\scalebox{.9}{\mbox{$\mathscr{E}$}}}}
\nc{\bchis}{\bchi{}^{\ms{2mu}*}}
\nc{\Daf}{\mathscr{D}}
\nc{\Laf}{\mathscr{L}}
\nc{\tLaf}{\widetilde{\Laf}}
\nc{\wtaf}{\mathscr{W}\ms{-3mu}{\mathit{t}}}
\nc{\res}[1][]{\mathop\star\limits_{\raisebox{.4ex}{$\scriptstyle #1$}}\ms{2mu}}
\nc{\hchi}{\widehat{\chi}}
\nc{\convaff}{\mathop{\scalebox{1.1}{$\mspace{1.5mu}\circ\mspace{1.5mu}$}}\limits}
\nc{\cvb}{CVB\xspace}
\nc{\svelt}{essentailly samll\xspace}
\nc{\Proc}{\on{Pro}_{\mathrm{coh}}}
\nc{\sha}{\mathcal{A}}
\nc{\Ker}{\on{Ker}}
\nc{\Coker}{\on{Coker}}
\nc{\Aff}[1][z]{\on{Aff}_{\ms{1mu}#1}}
\nc{\scb}{\scalebox}
\nc{\afr}{affreal\xspace}
\nc{\epifrom}{\ms{-5mu}\xymatrix@C=3ex{{}&{}\ar@{->>}[l]}\ms{-5mu}}
\nc{\Mid}{\bigm|}
\nc{\ol}{\overline}
\nc{\bpsi}{\ol{\psi}}
\nc{\Rat}[1][z]{\on{Raff}_{\ms{1mu}#1}}
\nc{\Indlim}{\varinjlim\limits}

\nc{\inddlim}{\mathop{\mbox{``{$\ms{1mu}\varinjlim$}''}}\limits}
\nc{\Rmat}{\mathrm{R}\ms{1mu}}
\nc{\Runi}{\mathrm{R}^{\mathrm{univ}}}
\nc{\Modg}{\mathrm{Mod}_{\mathrm{gr}}}
\nc{\KO}{quasi-rigid\xspace}
\nc{\hF}{\widehat{\F}}
\nc{\Modc}{\Mod_{\mathrm{coh}}}
\nc{\e}{\mathrm{e}}
\nc{\Idx}{\mathsf{\Lambda}}
\nc{\hA}{\widehat{A}}
\nc{\prood}{\mathop{\text{``}\prod\text{''}}\limits}

\nc{\hrefl}{\widehat{\mathscr{S}}}
\nc{\ev}{\mathrm{ev}}
\nc{\coev}{\mathrm{coev}}
\nc{\ihom}{\mathcal{H}om}
\nc{\tY}{\widetilde{Y}}
\nc{\tensz}{\tens[z]\ms{-3.5mu}}

\nc{\tRre}{\widetilde{\mathrm{R}}^{\mathrm{ren}}}
\nc{\dg}{\mathbf{\lambda}}
\nc{\htens}{\hconv}
\nc{\stens}{\sconv}
\nc{\Modgc}{\mathrm{Modg}_{\mathrm{coh}}}
\nc{\Aut}{\mathrm{Aut}}
\nc{\Rd}[1][\dg]{R_{#1}\gmod}
\nc{\nn}{\nonumber}
\nc{\Dual}{\mathrm{D}\ms{1mu}}
\nc{\DA}[1][A]{\ms{1mu}\mathrm{D}_{{#1}}}
\nc{\DmA}[1][A]{\ms{1mu}\Dual^{-1}_{{#1}}}
\nc{\tensa}{\tens[A]\ms{-3mu}}
\nc{\tensc}{\tens[{\ms{3mu}\cor}]\ms{-3mu}}

\nc{\ble}{\preccurlyeq}
\nc{\bge}{\succcurlyeq}
\nc{\Cwv}[1][{w,v}]{\catC_{#1}}

\nc{\afn}{affine object\xspace}
\nc{\afns}{affine objects\xspace}
\nc{\subafn}{affine subobject\xspace}
\nc{\Afns}{Affine objects\xspace}
\nc{\mr}{\mathrm{r}}
\nc{\ml}{\mathrm{l}}
\nc{\LQ}{\mathscr{L}}
\nc{\RQ}{\mathscr{R}}
\nc{\tM}{\widetilde{M}}
\nc{\tN}{\widetilde{N}}
\nc{\convz}[1][z]{\underset{#1}{\circ}}
\nc{\uw}[1][w]{\underline{#1}}
\nc{\Bw}[1][{\uw}]{\mathbf{B}_{#1}}
\nc{\CP}[1][{\uw}]{\mathrm{K}_{#1}} 
\nc{\CPP}[1][{\uw'}]{\mathrm{K}_{#1}}
\nc{\CK}[1][{w',w}]{\mathrm{K}'_{#1}} 
\nc{\Es}[1][i]{\mathrm{E}^*_{{#1}}\ms{1mu}}
\nc{\Est}[1][{w',w}]{\mathbf{E}^*_{#1}}
\nc{\CBw}[1][w]{\mathfrak{B}_{#1}}
\nc{\CB}{\mathfrak{B}}
\nc{\Em}[1][i]{E^{\max}_{#1}}
\nc{\Esm}[1][i]{\tE^{*\;\max}_{#1}}
\nc{\car}{\mathrm{ch}}
\nc{\tdC}{{\widetilde{\dC}}}
\nc{\rert}{\rt^{\mathrm{real}}}
\nc{\Qti}[1][i]{\Qt(\ang{#1})}
\nc{\Esi}[1][i]{\mathbf{E}^*_{#1}}
\nc{\iz}{\ang{i}_z}
\nc{\by}{\mathrel{/}}
\nc{\epss}{\eps^*}
\nc{\Qr}[1][w,s_i]{\mathrm{Q}^{\mathrm r}_{#1}}
\nc{\Cs}[1][w]{\catC_{{#1}}^{\ms{2mu}*}}
\nc{\ssim}{\raisebox{-.8ex}[.5ex][.2ex]{$\sim$}}
\nc{\rootl}{\mathsf{Q}}
\nc{\nconv}{\mathop{\mbox{\large $\odot$}}}
\nc{\bnam}{\be[{\rm(a)}]}
\nc{\wb}[1]{\mbox{$\rule[-1.1ex]{0ex}{2ex}#1$}}
\nc{\bcirc}{\mathbin{\raisebox{.3ex}{\scalebox{.65}{$\displaystyle\bigcirc$}}}}
\nc{\op}{{\mathrm{op}}}
\nc{\ttC}{\widetilde{C}}
\nc{\tLaa}{\widetilde{\Laa}}
\nc{\il}{{i_\ell}}
\nc{\bout}[1]{\ber\sout{ #1 }\er}

\numberwithin{equation}{section}

\title[Crystal structure on Localized category]%
{Crystal Structure of Localized Quantum Unipotent Coordinate Category}

\author[M. Kashiwara]{Masaki Kashiwara}
\thanks{The research of M.\ Kashiwara
was supported by Grant-in-Aid for Scientific Research (B) 23K20206,
Japan Society for the Promotion of Science.}
\address[M. Kashiwara]{%
Kyoto University Institute for Advanced Study, Research Institute
for Mathematical Sciences, Kyoto University, Kyoto 606-8502, Japan}

\email{masaki@kurims.kyoto-u.ac.jp}
\author[T. Nakashima]{Toshiki Nakashima}
\thanks{%
The research of T.Nakashima is supported by 
Grant-in-Aid for Scientific Research (C) 20K03564, 
Japan Society for the Promotion of Science.}
\address[T.Nakashima]{Division of Mathematics, 
Sophia University, Kioicho 7-1, Chiyoda-ku, Tokyo 102-8554,
Japan}
\email{toshiki@sophia.ac.jp}


\makeatletter
\@namedef{subjclassname@2020}{\textup{2020} Mathematics Subject Classification}
\makeatother


\date{September 2, 2025}

\begin{document}

\maketitle

\begin{abstract}
A localized quantum unipotent
coordinate category $\tCw$
 associated with a Weyl group element $w$
 is a rigid monoidal category which is obtained by applying  the localization process 
to the category $\Cw$, a subcategory of the category $R\gmod$ of 
 finite-dimensional graded modules of a quiver Hecke algebra $R$.
We shall show that the family of the isomorphism classes (up to grading shifts) of  
simple objects in $\tCw$ possesses a
crystal structure and it is isomorphic 
to the cellular crystal associated with $w$. As an application of this result, 
we shall show the connectedness of the crystal graph of an arbitrary cellular crystal. 
\end{abstract}

\tableofcontents

\section{Introduction}

One of the most potent tools in the crystal base theory (\cite{K1}) is 
the so-called ``crystal operators" $\eit, \fit$. These operators
 act on $\uq$-modules or the nilpotent subalgebra $\uqm\subset \uq$, 
where $\uq$ is a quantum algebra over the field $\Q(q)$ associated with 
a symmetrizable Kac-Moody Lie algebra $\g$. 
Let $(L,B)$ be the crystal basis of a $\uq$-module $M$ or $\uqm$, 
where $L$ is a free module over the local subring $\cA\subset\Q(q)$ at $q=0$ and $B$ is a $\Q$-basis of $L/qL$. 
The operators $\eit$ and $\fit$ act on $L$, 
the quotient space $L/qL$, and even the basis $B$. 
This induces an oriented graph structure on the basis $B$, 
known as the ``crystal graph". 
Once the crystal base theory is constructed, 
a fundamental and significant problem is to realize its crystal graph using various methods, 
such as tableaux realizations (\cite{KN}), path realizations of perfect crystals
(\cite{KMN1,KMN2}), geometric realizations (\cite{K-Saito}), 
categorical realizations (\cite{LV11},\cite{N2}), etc.

 \medskip
 On the other hand, in \cite{KKOP21, KKOP22, KKOP23}
 the localization of a monoidal category is studied.
 It is analogous to the standard localization of a commutative ring by a  multiplicative set.
To be more precise, 
for a commutative ring $\cor$ and 
a $\cor$-linear monoidal category
$(\catT,\tens)$, 
a left braider in $\catT$ is defined as a $(C, \coRl_C)$, 
consisting of an object $C$, and morphisms
\[
  \coRl_C(X)\cl C\tens X \to X\tens C
  \]
which are functorial in $X \in \catT$ and satisfy 
certain compatibility conditions with the tensor product $\tens$. 
A family $\{(C_i, \coRl_{C_i}\}_{i}$ of left braiders 
in $\catT$ is called a real commuting family if
$\coRl_{C_i}(C_i) \in \cor^\times \id_{C_i \tens C_i}$,
and $\coRl_{C_j}(C_i) \circ \coRl_{C_i}(C_j) \in \cor^\times \id_{C_i \tens C_j}$
for any $i,j$. 
Then, one can construct a $\cor$-linear monoidal category $(\lT, \tens)$ and a monoidal functor $\Phi \cl \catT \to \lT$ such that the objects $\Phi(C_i)$ are invertible, and the morphisms $\Phi(\coRl_{C_i}(X)) \cl \Phi(C_i) \tens \Phi(X) \to \Phi(X) \tens \Phi(C_i)$ are isomorphisms for all $i$ and all $X \in \catT$.


Let $R$ be a quiver Hecke algebra associated with a Kac-Moody Lie algebra $\g$, and let $R\gmod$ denote the category of finite-dimensional graded $R$-modules. For an element $w$ of the Weyl group $\weyl$, there exists a significant and intriguing subcategory $\catC_w$ of $R\gmod$ (see Section \ref{Sec: Cwv} for the precise definition). It has been established that the Grothendieck ring $K(\catC_w)$ of $\catC_w$ is isomorphic to the quantum unipotent coordinate ring $A_q(\n(w))$ associated with $w$ (see \cite{KKKO18} and \cite{KKOP18}). The algebra $A_q(\n(w))$ can be interpreted as a $q$-deformation of the coordinate ring of the unipotent subgroup $N(w)$, and $A_q(\n(w))$ possesses a quantum cluster algebra structure. It was demonstrated in \cite{KKKO18} that the category $\catC_w$ provides a monoidal categorification of $A_q(\n(w))$ as a quantum cluster algebra when $\g$ is symmetric.



In \cite{KKOP21}, it is shown that there exist graded left braiders in $R\gmod$:
$$
( \dM(w\La_i,\La_i),  \coRl_{\dM(w\La_i,\La_i)}, \dphi^l_{\dM(w\La_i,\La_i)} ) \qquad \text{ for any $i\in I$ }
$$ 
consisting of the \emph{determinantial module} $\dM( w\La_i, \La_i )$ (see \S\;\ref{subsec:determinantial})
and the morphism induced from \emph{R-matrices}
$$
\coRl_{\dM(w\La_i,\La_i)}(X) \cl \dM(w\La_i,\La_i) \conv  X  \longrightarrow
X  \conv   \dM(w\La_i,\La_i)\qt{ for any $X \in R\gmod$}
$$
of homogeneous degree
$- (w \La_i + \La_i, \wt(X)) $.

Under the categorification, the determinantial module $\dM(w\La_i, \La_i)$ corresponds to the unipotent quantum minor $D(w\La_i, \La_i)$, which is a frozen variable of the quantum cluster algebra $A_q(\n(w))$. Thus, the localization $\tcatC_w$ of $\catC_w$ with respect to the aforementioned left braiders categorifies the localization $A_q(\n(w))[D(w\La_i, \La_i)^{-1}; i \in I]$ of $A_q(\n(w))$ at the frozen variables $D(w\La_i, \La_i)$ ($i \in I$). Notably, through this categorification, the Grothendieck ring $K(\tcatC_w)$ can be interpreted as a quantization of the coordinate ring of the unipotent cell $N^w$ associated with $w$ (see \cite[Theorem 4.13]{KO21}).

It turns out that the natural inclusion functor $\iota_w \cl \catC_w \rightarrowtail R\gmod$ induces an equivalence of categories:
$$
\widetilde \iota_w \cl  \tcatC_w \buildrel \sim \over \longrightarrow \bl R\gmod \br[\dM(w\La_i,\La_i)^{\circ -1}; i\in I]=:\tRm[w].
$$ 
It should be noted 
that a bunch of objects in $R\gmod$ vanishes after the localization process. 
In \cite[Proposition~3.3]{KKOP22},
the characterization of modules which vanish under the localization functor $\Qt \cl R\gmod \to \tRm[w]$ is provided. Recall that the self-dual simple modules in $R\gmod$ correspond bijectively with the crystal $B(\infty)$ of $\Aq[\n_+]$ (\cite{LV11}). It is established that a simple module $M$ does not vanish under $\Qt$ if and only if $M$ corresponds to an element $b$ in $B_w(\infty)$, where $B_w(\infty)$ is a subset of $B(\infty)$ introduced in \cite{Kas93}, known as the Demazure crystal (see Proposition~\ref{Demazure}).

In \cite{N1}, the second author constructed a geometric crystal structure on the Schubert cell $X_w$ 
(resp.\ variety $\overline{X_w}$) associated with a Weyl group element $w$. 
It was shown that by tropicalizing this geometric crystal on the Langlands dual Schubert cell $X^\vee_w$
(resp.\ variety $\overline{X}^\vee_w$), 
one can obtain the crystal $B_{i_1} \otimes B_{i_2} \otimes \cdots \otimes B_{i_k}$, 
where $B_i \seteq \{(m)_i \mid m \in \bbZ\}$ ($i \in I$) is a crystal introduced in Example \ref{ex-cry} 
and $\uw \seteq s_{i_1}s_{i_2} \cdots s_{i_k}$ is a reduced expression of $w$. 
Note that for a Weyl group element $w$, 
a reduced expression $s_{i_1}s_{i_2} \cdots s_{i_k}$ is often identified with a reduced word $i_1i_2 \cdots i_k$. Here, the crystal $\Bw \seteq B_{i_1} \otimes B_{i_2} \otimes \cdots \otimes B_{i_k}$ is 
referred to as the {\it cellular crystal associated with $w$}.
In \cite{Kana-N}, Kanakubo and the second author proved that when $\g$ is semi-simple, 
the cellular crystal $\Bw = B_{i_1} \otimes B_{i_2} \otimes \cdots \otimes B_{i_k}$ 
associated with a reduced word $\uw = i_1i_2 \cdots i_k$ is connected as a crystal graph. It is well-known
 that the crystal $B(\infty)$ is embedded in ${\bf B}_{\underline{w_0}}$, 
where $\underline{w_0} = i_1i_2 \cdots i_N$ is a reduced word of the longest element in $\weyl$ (\cite{Kana-N,N4}). 
As a set, ${\bf B}_{\underline{w_0}}$ can be identified with $\bbZ^N$, and for any ${\bf x, \bf y} \in {\bf B}_{\underline{w_0}}$, the summation ${\bf x} + {\bf y}$ is naturally defined in $\bbZ^N$.

For $\La$ in the set $\pwtl$ of dominant integral weights, set
${\bf h}_\La := \til f_{i_1}^{m_1} \cdots \til f_{i_N}^{m_N}((0)_{i_1} \otimes \cdots (0)_{i_N})$, where $m_k \seteq \ang{h_{i_k}, s_{i_{k+1}} \cdots s_{i_N} \La}$.
Then we have 
\[
{\bf B}_{\underline{w_0}}=\st{-{\bf h}_\La+x\mid \La\in\pwtl,\;x\in B(\ify)},
\]
which implies the connectedness of ${\bf B}_{\underline{w_0}}$ and consequently $\Bw$ for any $w \in \weyl$ (\cite{Kana-N}).

\medskip
 In the meanwhile,
 Lauda and Vazirani (\cite{LV11})  showed  that the family $\Irr(R\gmod)$
 of the isomorphism classes (up to grading shifts) 
of simple modules in $R\gmod$ possesses a
crystal structure and provided an isomorphism of  the crystals 
\[
\Psi\cl \Irr(R\gmod)\,\,\overset{\sim}{\longrightarrow}\, B(\ify).
\]
For a semi-simple $\g$, we have 
$R\gmod=\catC_{w_0}$ and $\wtil{R\gmod}=\wtil{\catC_{w_0}}$. 
Let $\Phi\cl R\gmod\to \wtil{R\gmod}$ be the localization functor. In \cite{KKOP21} it  is  shown that, for any simple object, 
$Y\in \wtil{R\gmod}$ there exists $\La\in\pwtl$ and a simple object $X\in R\gmod$
such that $Y\simeq C_\La^{\circ -1}\circ \Phi(X)$, where 
$C_\La\seteq\dM(w\La, \La)=\tF_{i_1}^{m_1}\cdots \tF_{i_N}^{m_N}{\bf 1}$ is a determinantial module. Thus, it is observed that 
\eqn
\Irr(\wtil{R\gmod})\simeq{\bf B}_{\underline{w_0}}\qt{
by $C_\La^{\circ -1}\circ X\longleftrightarrow -{\bf h}_\La+\Psi(X)$,}
\eneqn
\vspace{-1pt}where $\Irr(\wtil{R\gmod})$ is
the family of the isomorphism classes
(up to grading shifts) of simple objects in  $\wtil{R\gmod}$.

Then, the second author proved (\cite{N2}) that
the family of self-dual simple objects in 
the localized category possesses a crystal 
structure and it is isomorphic to the cellular crystal ${\bf B}_{\underline{w_0}}\seteq B_{i_1}\otimes B_{i_2}\ot\cdots\otimes B_{i_N}$. It implies that the family 
$\Irr(\widetilde{R\gmod})$ also possesses a crystal 
structure and it is isomorphic to the cellular crystal ${\bf B}_{\underline{w_0}}$ since for any simple module $X$ in $R\gmod$, there exists a unique integer $n$ 
such that $q^n X$ is self-dual (see Lemma~\ref{self-n}).

In \cite{N2}   the following problem is presented: in the general Kac-Moody setting, 
$\Irr(\tC_w)$  is endowed with 
a crystal structure and it is isomorphic to the cellular 
crystal $\Bw\seteq B_{i_1}\otimes B_{i_2}\ot\cdots\otimes B_{i_k}$ associated with an arbitrary reduced word $\uw=i_1i_2\cdots i_k$ of an arbitrary $w\in\weyl$. 
The main aim of this paper is to provide a positive answer to this problem.

Let us delve into more details.
In \cite{N2}, as mentioned above, the crystal structure of $\Irr(\widetilde{R\gmod})$
is established.
However, the method employed there cannot be applied to our problem since 
 $\Irr(\Cw)$  is not  necessarily stable under the action of the 
 crystal operator $\tE_i$ (see Remark~\ref{rm:non-stable}). Additionally,  an object in the localized category $\tCw$ is no longer an $R$-module,
and thus the action of  $\tE_i$  on $\tCw$ is not properly defined. 
Therefore,  new tools are required to define the crystal structure on $\Irr(\tCw)$.
 They  are ``affinization", ``R-matrix" and ``rigidity''.

In case $R$ is symmetric, an affinization of $M\in R\gmod$ is simply given as 
${\bf k}[z]\ot_{\bf  k}M$. However, in general, the existence of an affinization in $R\gmod$ 
is not necessarily guaranteed. A general definition of affinization in
a monoidal category
introduced in \cite{ref2} (cf.\ \cite{KP18}) is provided in Definition~\ref{def:affinization}.
An object $M$ is called {\em affreal} if $M$ is real and afford an affinization
(see Definition~\ref{def:affreal}). Then by Proposition~\ref{prop:simplehd}
we find that if $M$ is affreal and $N$ is simple in a
quasi-rigid category (see \S\,\ref{subsec:simplehd}), then 
there exists R-matrix  $\rmat{M,N}\cl M\tens N\to N\tens M$.
As stated in Theorem~\ref{thm:rigidity-Cw},
the category $\tCw$ is rigid (\cite{KKOP21}, \cite{KKOP22}), 
that is, any object $X$ of $\tCw$ holds the left dual 
$\D^{-1}(X)$ and the right dual $\D(X)$ (see Definition~\ref{def:rigid}).
These tools enable us to define a crystal structure on $\Irr(\tCw)$ as in Definition~\ref{def:rootop}. Finally, we shall construct the isomorphism of crystals
${\rm K}_{\uw}\cl\Irr(\tCw)\overset{\sim}{\longrightarrow}\Bw$ inductively using the map
\eqn
&&{\rm K'}_{w',w}\cl \Irr(\tCw)\to \Irr(\tilde\catC_{w'})\times \Z\quad
\hbox{ given by }\quad X\mapsto  (\Est(X),\eps^*_i(X)),
\eneqn
where $w'=ws_i<w$ and  see Definition~\ref{def:rootop} and \S\,\ref{map-E*} for the definition of
$\Est$ and $\ep^*_i$.

We note that Woo-Seok Jung and Euiyong Park
(\cite{JP}) treats relevant topics,
including the bijectivity of $\Irr(\tCw)\to \Bw$.

\medskip
The organization of the article is as follows. 
In \S\,2, we review the basics of the crystal base theory, in particular, the explicit structure of 
the cellular crystals which will be needed for later arguments.
In \S\,3, the theory of affinization and R-matrix of monoidal category is provided following \cite{ref}, which will play a crucial role to obtain the main theorem in the subsequent sections. \S\,4 is devoted to give the definition of the quiver Hecke algebra and the useful and important properties of its modules such as 
the head simplicity of convolution product, R-matrices, the invariants $\La$, the shuffle lemmas and 
the determinantial modules.
In \S\,5, several results on the localization of the category $R\gmod$ and $\Cw$ are introduced and 
the relations of the localization functors and the convolution product are discussed.
In \S\,6, we define the crystal structure on $\Irr(\tCw)$ and 
present some technical lemmas to prove the main theorems. As a corollary  of the main theorem in \S\,7, we prove that there exists an isomorphism between $\Irr(\tCw)$ and the cellular crystal $\Bw$.
  The proof of the main theorem  
  will be presented in \S\,8. In the last section,
  we prove the connectedness of arbitrary cellular crystal
 as an application of the main theorem.


 \section{Crystal Bases and Crystals}
\subsection{Convention}
In this paper. we use the following convention.
\begin{convention}
\bnum
\item A ring is a unital ring.
\item For a ring $A$, $A^\times$ denotes the group of invertible elements of $A$.
\item For a monoidal category $\shc$ with a tensor product $\tens$,
  we denote by $\shc^\rev$ the monoidal category $\shc$ with the reversed
  tensor product $\tens_\rev$: $X\tens_\rev Y\seteq Y\tens X$.
  \ee
\end{convention}

\subsection{Crystal Bases}
A Cartan datum $ \bl\cmA,\wlP,\Pi,\Pi^\vee,(\cdot,\cdot) \br $ is a quintuple of a generalized Cartan matrix $\cartan$,    a free abelian group $\wlP$,  a set of simple roots, $\Pi = \{ \alpha_i \mid i\in I \} \subset \wlP$,  a  
 set of simple coroots $\Pi^{\vee} = \{ h_i \mid i\in I \} \subset \wlP^{\vee}\seteq\Hom( \wlP, \Z )$ ,  and a $\Q$-valued 
symmetric bilinear form  $(\cdot,\cdot)$ on $\wlP$ such that
\bna 
\item $\cmA = (\langle h_i,\alpha_j\rangle)_{i,j\in I}$,
\item  $(\alpha_i,\alpha_i)\in 2\Z_{>0}$ for any $i\in I$,
\item $\langle h_i, \lambda \rangle =\dfrac{2(\alpha_i,\lambda)}{(\alpha_i,\alpha_i)}$ for $i\in I$ and $\lambda \in \Po$,
\item for each $i\in I$, there exists $\Lambda_i \in \wlP$
such that $\langle h_j, \Lambda_i \rangle = \delta_{ij}$ for any $j\in I$.
\end{enumerate}

Let $\rlQ\seteq\soplus_{i\in I} \Z\al_i$
and $\rlQ_+\seteq\soplus_{i\in I} \Z_{\ge0} \al_i$ be the root lattice and the positive root lattice of the symmetrizable Kac-Moody algebra $\g=\g(\cartan)$,  respectively.
Let $U_q(\g)\seteq\lan e_i,f_i,q^h\mid i\in I,\,h\in P^\vee \ran$
(resp.\ $U_q^-(\g)\seteq\lan f_i\mid i\in I\ran$) be the associated quantum
group over $\Q(q)$ (resp.\ negative nilpotent subalgebra of $U_q(\g)$).

Let $\weyl$ be the Weyl group of $\g$,  the subgroup of
$\Aut(\wtl)$ generated by the simple reflections $\st{s_i}_{i\in I}$ where
$s_i(\la)=\la-\ang{h_i,\la}\al_i$ for $\la \in \wtl$.

\medskip

As in \cite{K1}, there exists the crystal basis $(L(\ify),B(\ify))$ of $\uqm$ defined by
\begin{eqnarray*}
&& L(\ify)\seteq\sum_{k\geq0,i_1,\cdots,i_k\in I}\A\til f_{i_1}\cdots\til
 f_{i_k}u_\ify,\\
&&B(\ify)=\{\til f_{i_1}\cdots\til
 f_{i_k}u_\ify\,{\rm mod}\, q L(\ify)\,|\,k\geq0, i_1,\cdots,i_k\in I\}\setminus\{0\},\\
&&\vep_i(b)={\rm max}\{k\mid\eit^kb\ne0\},\quad
\vp_i(b)=\vep_i(b)+\lan h_i,\wt(b)\ran,
\end{eqnarray*}
where $u_\ify=1\in \uq$, 
$\eit$ and $\fit\in{\rm End}_{\Q(q)}(\uqm)$ are the {\em crystal operators}
(\cite{K1}) and 
$\A$ is the local subring of $\Q(q)$ at $q=0$. We also obtain the crystal basis 
$(L(\la),B(\la))$ of the irreducible highest weight module $V(\la)$ with a highest weight vector $u_\la$ for any 
dominant weight $\la\in \wlP_+\seteq\{\la\in\wlP\mid \lan h_i,\la\ran\in\Z\}$ by a 
similar way.

\subsection{Crystals}
We define the notion of {\em crystal} as in \cite{Kas93}, which is 
a combinatorial object abstracting the properties of crystal bases:
\Def[\cite{Kas93}]\label{cryst}
A 6-tuple $\bl B,\wt, \{\vep_i,\vp_i, \eit,\fit\}_{i\in I}\br$ is 
a {\em crystal} if $B$ is a set and there exists a certain special element 
$0$ outside of $B$ and maps:
\begin{eqnarray}
&&{\rm wt}\cl B\to P,\label{wtp-c}\quad
\vep_i\cl B\to\Z\sqcup\{-\ify\},\quad\vp_i\cl B\to\Z\sqcup\{-\ify\}\quad\,(i\in I),
\\
&&\eit\cl B\sqcup\{0\}\to B\sqcup\{0\},\quad
\fit\cl B\sqcup\{0\}\to B\sqcup\{0\}\,\,(i\in I),\label{eitfit-c}
\end{eqnarray}
satisfying :
\begin{enumerate}
\item
$\vp_i(b)=\vep_i(b)+\lan h_i,\wt(b)\ran$.
\item
If  $b,\eit b\in B$, then $\wt(\eit b)=\wt(b)+\al_i$, 
$\vep_i(\eit b)=\vep_i(b)-1$, $\vp_i(\eit b)=\vp_i(b)+1$.
\item
If $b,\fit b\in B$, then $\wt(\fit b)=\wt(b)-\al_i$, 
$\vep_i(\fit b)=\vep_i(b)+1$, $\vp_i(\fit b)=\vp_i(b)-1$.
\item
For $b,b'\in B$ and $i\in I$, one has 
$\fit b=b'$ if an only if  $b=\eit b'.$
\item
If $\vp_i(b)=-\ify$ for $b\in B$, then $\eit b=\fit b=0$
and $\eit(0)=\fit(0)=0$.
\end{enumerate}
\edf

\Def\label{def:graph}
A {\em crystal graph} of a crystal $B$ is an $I$-colored 
oriented graph defined by $b\mapright{i} b'\Leftrightarrow \fit(b)=b'$ for $b,b'\in B$.
A crystal $B$ is {\em connected} if its crystal graph is a connected graph. 
\edf

\Def[\cite{Kas93}]\label{morph}\hfill 
\bnum
\item Let $B=(B,\wt,\{\varepsilon_i\},\{\varphi_i\},\{{\tilde e}_i\},\{{\tilde f}_i\})$ and $
  B'=(B', \wt', \{\varepsilon'_i\},\{\varphi'_i\},\{{\tilde e}'_i\},
  \{{\tilde f}'_i\})$ be crystals.
  A morphism of crystal from $B$ to $B'$ is a map
  $\psi\cl B\sqcup\st{0}\to B'\sqcup\st{0}$
  such that
  \bna
\item $\psi(0)=0$,
\item if $b\in B$ satisfies $\psi(b)\in B'$, then
  $\wt\bl\psi(b)\br=\wt(b)$, $\eps_i\bl\psi(b)\br=\eps_i(b)$,
  and  $\vphi_i\bl\psi(b)\br=\vphi_i(b)$, 
\item if $b_1,b_2\in B$ satisfy $\psi(b_1)\in B'$,
  $\psi(b_2)\in B'$ and $\tf_i(b_1) =b_2$, then we have
  $\tf_i'\bl\psi(b_1)\br=\psi(b_2)$.
 \ee
  \item 
    Let
    $B=(B,\wt,\{\varepsilon_i\},\{\varphi_i\},\{{\tilde e}_i,\},\{{\tilde f}_i\})$ be a crystal, and let $B'$ be a subset of
    $B$. Then $B'$ is endowed with the crystal structure
    $B'=(B', \wt', \{\varepsilon'_i\},\{\varphi'_i\},\{{\tilde e}'_i\}, \{{\tilde f}'_i\})$ induced by $B$: $\wt'$, $\varepsilon'_i$, $\varphi'_i$ are the restrictions of
    the corresponding
    functions on $B$ and
\[
\hs{3ex}\te_i'(b)=\bc \te_i(b)&\text{if $\te_i(b)\in B'$,}\\
    0&\text{otherwise,}\ec\q
    \tf_i'(b)=\bc \tf_i(b)&\text{if $\tf_i(b)\in B'$,}\\
    0&\text{otherwise}\ec\qt{for any $b\in B'$,}
\]
so that $B'\to B$ is a morphism of crystals.
The crystal $B'$ is called a {\em subcrystal} of $B$. 
\item 
We say that a morphism $\Psi\cl B_1\to B_2$ is {\em strict}
if it satisfies
$\te_i(\Psi(b))=\Psi(\te_ib)$ and 
$\tf_i(\Psi(b))=\Psi(\tf_ib)$ for any $b\in B_1$. 
\item 
We say that  a morphism
$\Psi\cl B_1\to B_2$ is an {\em embedding} if 
$B_1\sqcup\{0\}\to B_2\sqcup\{0\}$ is injective.
\end{enumerate}
\edf
\begin{prop}[\cite{Kas93}]\label{Demazure}
For any $w\in \weyl$ there exists a unique subset $B_w(\ify)$ of $B(\ify)$ satisfying the 
following properties:
\begin{enumerate}
\item[\rm(1)] $B_w(\ify)=\{u_\ify\}$ if $w=1$.
\item[\rm(2)] If $s_iw<w$, then $\displaystyle B_w(\ify)=\bigcup_{k\geq0}\fit^kB_{s_iw}(\ify)$.
  \ee
  Moreover, they satisfy
\be  
\item[\rm(3)] $\eit B_w(\ify)\subset B_{s_iw}(\ify)$.
\item[\rm(4)] If $w\geq w'$, then $B_w(\ify)\supset B_{w'}(\ify)$.
\end{enumerate}
We call $B_w(\ify)$ the {\em Demazure crystal} associated with $w\in \weyl$
with the crystal structure induced from $B(\ify)$. 
\end{prop}
We can define the tensor product of crystals as follows (\cite{K1,Kas93}):
\begin{prop}\label{tensor}
For crystals $B_1$ and $B_2$,  
set 
\[
B_1\ot B_2=\{b_1\ot b_2\mid b_1\in B_1 ,\, b_2\in B_2\}
\]
Then, $B_1\ot B_2$ becomes a crystal by defining:
\eq
&&\wt(b_1\ot b_2)=\wt(b_1)+\wt(b_2),\\
&&\vep_i(b_1\ot b_2)={\max}(\vep_i(b_1),
  \vep_i(b_2)-\lan h_i,\wt(b_1)\ran),
\label{tensor-vep}\\
&&\vp_i(b_1\ot b_2)={\max}(\vp_i(b_2),
  \vp_i(b_1)+\lan h_i,\wt(b_2)\ran),
\label{tensor-vp}\\
&&\eit(b_1\ot b_2)=
\left\{
\begin{array}{ll}
\eit b_1\ot b_2 & \text{if $\vp_i(b_1)\ge \vep_i(b_2)$,}\\
b_1\ot\eit b_2  & \text{if $\vp_i(b_1)< \vep_i(b_2)$,}
\end{array}
\right.
\label{tensor-e}
\\
&&\fit(b_1\ot b_2)=
\left\{
\begin{array}{ll}
\fit b_1\ot b_2 & \text{if $\vp_i(b_1)>\vep_i(b_2)$,}\\
b_1\ot\fit b_2  & \text{if $\vp_i(b_1)\le \vep_i(b_2)$.}
\label{tensor-f}
\end{array}
\right.
\eneq
\end{prop}


\begin{ex}\label{ex-cry}
\begin{enumerate}
\item
For $i\in I$, set
$ B_i\seteq\{(n)_i\,|\,n\in \Z\}$ and 
\begin{eqnarray*}
&&\wt((n)_i)=n\al_i,\,\,
\vep_i((n)_i)=-n,\,\,
\vp_i((n)_i)=n,\,\,\label{bidata}\\
&&
\vep_j((n)_i)=\vp_j((n)_i)=-\ify\,\,\,(i\ne j),\label{biji}\\
&&
\eit((n)_i)=(n+1)_i,\quad\fit((n)_i)=(n-1)_i,\label{bieitfit}\\
&&
\til e_j((n)_i)=\til f_j((n)_i)=0\,\,\,(i\ne j).
\label{biejfj}
\end{eqnarray*}
Then $B_i$ ($i\in I$) becomes a crystal.
Note that $\eit$ and $\fit$ are inverse to each other. 
\item
The crystal $T_\la\seteq\{t_\la\}$ is defined by 
$\eit t_\la=\fit t_\la=0$, $\vep_i(t_\la)=\vep_i(t_\la)=-\ify$ and $\wt(t_\la)=\la$.
\end{enumerate}
\end{ex}

\begin{ex}\label{ex-cry2}
  \hfill
  
\bnum
\item There exists a unique strict embedding $\Psi_i\cl B(\ify)\to B(\ify)\ot B_i$ ($u_\infty\mapsto u_\infty
\ot (0)_i$) for any $i\in I$
(\cite{Kas93}).
\item
There is a 
unique  embedding 
$\Psi_\la\cl B(\la)\to B(\ify)\ot T_\la$ ($u_\la\mapsto u_\ify\ot t_\la$)
for any $\la\in\wlP_+$.
  Here,  note that $\vep_i(b\ot t_\la)=\vep_i(b)$ and 
$\eit(b\ot t_\la)=(\eit b)\ot t_\la$. Thus, the map $\Psi_\la$ induces 
an injective map $\Psi\cl B(\la)\to B(\ify)$ such that $\Psi\eit(b) =\eit\Psi(b)$ and 
$\vep_i(b)=\vep_i(\Psi(b))$ for any $b\in B(\la)$.
\end{enumerate}
\end{ex}

\subsection{Cellular Crystal $\Bw$}\label{cellular}
For a reduced expression $\uw=s_{i_1}s_{i_2}\cdots s_{i_k}$ of some Weyl group element $w\in\weyl$, 
we call the crystal $\Bw\seteq B_{i_1}\ot\cdots\ot B_{i_k}$ 
the {\em cellular crystal}
associated with a reduced expression $\uw$. Note that
it is obtained by applying the tropicalization functor to the 
geometric crystal on the Langlands-dual Schubert cell $X_w^\vee$ (\cite{N1}).
It is immediate from the braid-type isomorphisms (\cite{N3}) that 
for any $w\in \weyl$ and its arbitrary reduced words $i_1\cdots i_m$ and 
$j_1\cdots j_m$, we get the following isomorphism of crystals:
\begin{equation}
B_{i_1}\ot\cdots\ot B_{i_m}\simeq B_{j_1}\ot\cdots\ot B_{j_m}.
\label{iso-cellular}
\end{equation}

Here, we shall investigate an explicit structure of the crystal $\Bw$.
Fix a reduced expression $\uw=s_{i_1}\cdots s_{i_m}$
and write 
\[
(x_1,\cdots,x_m)\seteq\til f_{i_1}^{x_1}(0)_{i_1}
\ot\cdots\ot \til f_{i_m}^{x_m}(0)_{i_m}=(-x_1)_{i_1}\ot\cdots\ot(-x_m)_{i_m}\in \Bw,
\]
where $\fit^n(0)_i$ means $\eit^{-n}(0)_i$ when $n<0$.

By the tensor structure of crystals 
in Proposition \ref{tensor},
we can describe 
the explicit crystal structure on $\Bw\seteq B_{i_1}\ot \cdots\ot B_{i_m}$ 
as follows: 
For $x=(x_1,\cdots,x_m)\in \Bw$ and $k\in[1,m]$, define
\begin{equation*}
 \sigma_k(x)\seteq x_k+\sum_{1\leq j<k}\lan h_{i_k},\al_{i_j}\ran x_j
\label{sigma-def}
\end{equation*}
and for $i\in I$ define
\begin{eqnarray}
&&\TY(\wtil\sigma,i,)(x)\seteq\max\{\sigma_k(x)\,|\,1\leq k\leq m\,{\rm and}\,
i_k=i\},\label{simgak}\\
&&\TY(\wtil M,i,)=\TY(\wtil M,i,)(x)\seteq\{k\,|\,1\leq k\leq m,\,i_k=i,\,
\sigma_k(x)=\TY(\wtil\sigma,i,)(x)\},
\label{Mi}
\\
&&
\TY(\wtil m,i,f)=\TY(\wtil m,i,f)(x)\seteq\max\,\TY(\wtil M,i,)(x),\q
\TY(\wtil m,i,e)=\TY(\wtil m,i,e)(x)\seteq\min\,\TY(\wtil M,i,)(x).
\label{mfme}
\end{eqnarray}
Now, the actions of the crystal operators $\eit,\fit$ and the functions
$\vep_i,\vp_i$ 
and $\wt$  are written explicitly:
\begin{eqnarray}
&& \fit(x)_k\seteq x_k+\del_{k,\TY(\wtil m,i,f)},\qq\qq
\eit(x)_k\seteq x_k-\del_{k,\TY(\wtil m,i,e)},\\
&&
\wt(x)\seteq-\sum_{k=1}^m x_k\al_{i_k},\quad
\vep_i(x)\seteq \TY(\wtil\sigma,i,)(x),\quad
\vp_i(x)\seteq\lan h_i,\wt(x)\ran+\vep_i(x).
\label{wt-vep-vp-ify}
\end{eqnarray}

\section{Generality of graded monoidal categories}
In this section, we briefly review the notions and the results
in \cite{ref} which axiomatize the contents of \cite{KKKO15}, \cite{KKKO18},
\cite{KKOP18}, etc.  For more details, consult \cite{ref}.

\subsection{Affine objects in $\Pro(\shc)$} 

\Def[see e.g.\ {\cite[\S\,1.4, 3.1, 6.1]{KS06}}]
\begin{enumerate}
\item
For a category $\shc$, let $\shc^\vee$ denote the category of functors from 
$\shc^{\rm op}$ to ${\bf  Set}^{\rm op}$. Then, we have a fully faithful functor
$\shc\to\shc^\vee$
given by $\shc\ni X\mapsto(\shc\ni Y\mapsto\Hom_\shc(X,Y)\in\mathbf{Set})$, called the 
{\it Yoneda functor}.
\item A category $I$ is {\em directed} if it satisfies the following conditions:
\begin{enumerate}
\item $I$ is non-empty,
\item for any $i,j\in I$, there exist $k\in I$ and  
morphisms $i\to k$, $j\to k$,
\item for any parallel morphisms $f,g\cl i\rightrightarrows j$, there exists a morphism $h\cl j\to k$ such that $h\circ f=h\circ g$.
\end{enumerate}
If $I^\op$ is directed, we say that $I$ is {\em co-directed}.
\item A {\em pro-object} in $\shc$ 
is an object of $\shc^\vee$ which is isomorphic to
$\prolim\beta$ for some functor $\beta\cl I\to \shc\to\shc^\vee$ with co-directed $I$.
\item Let $\Pro(\shc)$ denote the full subcategory of $\shc^\vee$ consisting 
  of pro-objects.
  \item  The co-directed projective limit in $\Pro(\shc)$ is denoted by $\proolim$.
\end{enumerate}
\edf
Let $\cor$ be a base field. Let $\shc$ be a $\cor$-linear category such that 
\eq&&\hs{2ex}\left\{
\parbox{75ex}{\be[{$\bullet$}]
\item $\shc$ is abelian,
\item $\shc$ is $\cor$-linear graded, namely $\shc$ is endowed with a $\cor$-linear auto-equivalence $q$,

\item any object has a finite length,
\item  any simple object $S$ is absolutely simple, i.e.,
  $\cor\isoto\END_\shc(S)$.
  \ee
}\right.
\label{cond:fcat}
\eneq
Here, we denote
\eq
&&\hs{3ex}\ba{l}\HOM_\shc(M,N)=\soplus_{n\in\Z}\HOM_\shc(M,N)_n\hs{1ex} \text{with $\HOM_\shc(M,N)_n=\Hom_\shc(q^nM,N)$,}\\
\text{and $\END_\shc(M)\seteq\HOM_\shc(M,M)$.}
\ea
\eneq

It follows from \eqref{cond:fcat} that $\HOM_{\shc}(M,N)$ is finite-dimensional over $\cor$ for any $M,N \in \shc$ and
that $S\not \simeq q^kS $ for any $k\in \Z\setminus\{0\}$ and any non-zero 
$S \in \shc$.

By \cite[Lemma 2.1]{ref}, 
the full subcategory $\shc$ of $\Pro(\shc)$ is stable by taking subquotients.
The grading shift functor is extended to $\Pro(\shc)$. 

Let $A=\soplus_{k\in \Z} A_k$ be a commutative graded $\cor$-algebra such that 
\eq&&\left\{
\parbox{70ex}{
\be[$\bullet$]
\item $A_0\simeq\cor$,
\item $A_k=0$ for any $k<0$,
\item $A$ is a finitely generated $\cor$-algebra.
\ee}\right.
\label{cond:gring}
\eneq
In particular, $A$ is a noetherian ring and $\dim_\cor A_k < \infty$ for all $k$.
For each $m\in \Z$, set
\eqn 
A_{\ge m} \seteq\soplus _{k\ge m} A_k, 
\eneqn
which is an ideal of $A$.

{\em Hereafter, we assume that
$\shc$ satisfies \eqref{cond:fcat} and $A$ satisfies \eqref{cond:gring}.}

We denote by $\Modg(A, \Pro(\shc))$ 
the category of graded $A$-modules in $\Pro(\shc)$,
i.e., objects $\Ma\in\Pro(\shc)$ with a graded algebra homomorphism
$A\to\END_{\Pro(\shc)}(\Ma)$.

 For $\Ma, \Na \in \Modg(A, \Pro(\shc))$, we set 
$\HOM_{A}(\Ma,\Na)\seteq\soplus_{k\in \Z} \, \HOM_A(\Ma,\Na)_k  \, $ with
$\HOM_{A}(\Ma,\Na)_k\seteq\Hom_{\Modg(A, \Pro(\shc))}(q^k \Ma,\Na)$.   
Then
$\HOM_A(\Ma,\Na)$ has a structure of a graded $A$-module.   

\begin{definition} \label{Def: Pro_coh}
We denote by $\Proc(A,\shc)$ the full subcategory of $\Modg(A, \Pro(\shc))$ consisting of $\Ma\in\Modg(A, \Pro(\shc))$ such that
\be[(a)]
\item $\Ma/A_{>0}\Ma \in \shc$,
\item $\Ma \isoto \proolim[k] \Ma / A_{\ge k} \Ma$, or equivalently $\proolim[k] A_{\ge k}\Ma  \simeq 0$.
\ee
Here $A_{\ge k}\Ma  \seteq \Im(A_{\ge k} \tens_A \Ma \to A \tens_A \Ma  ) \subset \Ma$. 
Note that $ X \tens_A  \Ma $ belongs to $\Modg(A, \Pro(\shc))$ for any 
$\Ma\in\Modg(A, \Pro(\shc))$ and any finitely generated graded $A$-module $X$.
\end{definition}

\Def \label{Def: affine objects}
Let $z$ be an indeterminate of homogeneous degree $d\in\Z_{>0}$.
An object
of $\Proc(\cor[z],\shc)$ is nothing but a pair $(\Ma,z)$ such that
\bna
\item $\Ma\in\Pro(\shc)$ and $z\in\END_{\Proc(\shc) }(\Ma)_d$,
\item $\Ma/z\Ma\in\shc$,
\item $\Ma\isoto\proolim[n]\Ma/z^{n}\Ma$.
\setcounter{myc}{\value{enumi}}
\ee
If $(\Ma,z)$ satisfies further the following condition

\bna\setcounter{enumi}{\value{myc}}
\item $z\in\END_{\Proc(\shc) }(\Ma)$ is a monomorphism,
\ee
  then we say that $(\Ma,z)$ is an {\em \afn} (in $\shc$).

\end{definition}
We denote by $\Aff(\shc)$ the category of 
\afns.

\subsection{Graded monoidal categories }

In the sequel, 
let $(\shc,\tens)$ be an abelian $\cor$-linear graded monoidal category which satisfies
\eq&&\left\{\parbox{70ex}{
\bna
\item $\shc$ satisfies \eqref{cond:fcat},
\item $\tens$ is $\cor$-bilinear and bi-exact,
\item the unit object $\one$ is simple and satisfies 
$\END_{\shc}(\one)\simeq\cor$, 
\item $\tens$ is compatible with the grading shift functor $q$, \\
i.e., $q(X\tens Y)\simeq (qX)\tens Y\simeq X\tens(qY)$ (see \cite[(1.8)]{ref2}).
\setcounter{mycc}{\value{enumi}}
\ee}\right.
\label{cond:exactmono}
\eneq
For generalities on monoidal categories, we refer the reader 
to \cite[Chapter 4]{KS06}. 

By identifying $q$ with the invertible central object $q\one\in\shc$, we have
$qX\simeq q\tens X$.

The category $\Pro(\shc)$ has also a structure of monoidal category in which the tensor product $\tens$ is bi-exact.

\Def\label{def:rigid}
An object $X$ \ro resp.\ $Y$\rf of a monoidal category $\sha$ has a
{\em right  dual \ro resp.\ left dual\/\rf} $Y$ (resp.\ $X$) if there exists an object $Y\in\sha$ \ro resp.\ $X\in\sha$\rf 
and morphisms
$\ep\cl X\tens Y\to{\bf 1}$ and $\eta\cl{\bf 1}\to Y\tens X$ such that
the compositions
\eqn
&&X\simeq X\tens \one
\To[\id_X\tens \eta]X\tens Y\tens X
\To[\ep\tens\id_X]
\one\tens X\simeq X,\\
&&Y\simeq \one\tens Y
\To[\eta\tens\id_Y]Y\tens X\tens Y
\To[\id_Y\tens \eps]Y\tens\one\simeq Y
\eneqn
are the identities. We denote the right dual of $X$ by $\D(X)$ and the left dual of 
$Y$ by $\D^{-1}(Y)$. If any object in $\sha$ has a right dual and a left dual, we say that
$\sha$ is a {\em rigid} category.
\edf

\Def[{\cite[Definition 4.12]{ref}}] \label{Def: Raff}
We denote by $\Rat(\shc)$ the category 
with $\Ob\bl\Proc(\cor[z],\shc)\br$  
as the set of objects and with the morphisms defined as follows. 
For $\Ma,\Na\in\Ob\bl\Rat(\shc)\br$,
\eqn\Hom_{\,\Rat(\shc)}(\Ma,\Na)&&=\Indlim_{k\in\Z_{\ge0}}
\Hom_{\,\Proc(\cor[z],\shc)}\bl(z^k\cor[z])\tens[{\cor[z]}]\Ma,\Na\br\\
&&\simeq
\Indlim_{k\in\Z_{\ge0}}
\Hom_{\,\Proc(\cor[z],\shc)}(\Ma,\cor[z]z^{-k}\tens_{\cor[z]}\Na).\eneqn
\edf
Note that we have
\[
\HOM_{\Rat(\shc)}(\Ma,\Na)\simeq
\cor[z,z^{-1}]\tens_{\cor[z]}\HOM_{\Proc(\cor[z],\shc)}(\Ma,\Na)
\qt{for $\Ma,\Na\in\Proc(\cor[z],\shc)$.}
\]
We call $\Rat(\shc)$ the category of rational affine objects of $\shc$.

The category $\Rat(\shc)$ is an abelian category with $\cor$-finite-dimensional hom.
We have a faithful functor $\Aff(\shc)\to\Rat(\shc)$, which is essentially surjective
\cite[Lemma 4.13]{ref}.

\subsection{\Afns}
Let $\shc$ be an abelian graded $\cor$-linear monoidal category which satisfies
\eqref{cond:exactmono}.
The category  $\Aff(\shc)$ has a monoidal category structure with $\tensz$ as a tensor product, where 
$\Ma\tensz \Na\seteq\Coker(\Ma\tens \Na\To[z\tens1-1\tens z]\Ma\tens \Na)$.

\Prop[{\cite[Proposition 5.6]{ref}}]\label{prop:rigid-duality}
Assume that $\shc$ is rigid.
Then $\Aff(\shc)$ is a rigid monoidal category.
We denote by
$\Da^{\pm1}$ their right and left duality functors.
\enprop

\Def \label{def:rmat}
Let $M,N\in\shc$ be simple objects.
\bnum
\item
  If $\dim\HOM(M\tens N,N\tens M)=1$,
then we say that the pair $(M,N)$ is {\em $\La$-definable},
and a non-zero
morphism $\rmat{}\in \HOM(M\tens N,N\tens M)$ is called the {\em R-matrix}
between $M$ and $N$, denoted by $\rmat{M,N}$. It is well-defined up to a constant multiple.
We set
$$\La(M,N)\seteq\deg(\rmat{}).$$
\item
  If $(N,M)$ is also $\La$-definable in addition, we say that $(M,N)$ is
  {\em $\de$-definable} and  write
  $$\de(M,N)\seteq\dfrac{1}{2}\bl\La(M,N)+\La(N,M)\br\in\dfrac{1}{2}\Z.$$
  \ee
\edf


\subsection{Rational centers and affinizations}
Let $\shc$ be a $\cor$-linear graded monoidal category which satisfies
the following conditions:
\eq
&&\hs{1ex}\left\{\parbox{70ex}{
$\shc$ satisfies \eqref{cond:exactmono} and the following additional condition:
\bna
\setcounter{enumi}{\value{mycc}}
\item $\shc$ has a decomposition $\shc=\soplus_{\la\in\Idx}\shc_\la$
compatible with a monoidal structure where $\Idx$ is an abelian monoid,
and $\one\in \shc_\la$ with $\la=0$. 
\ee
}\right.\label{cond:g}
\eneq

Let $z$ be a homogeneous indeterminate with degree $d\in\Z_{>0}$.

Remark that we have also bi-exact bifunctors:
$$\tens\;\cl \shc\times\Rat(\shc)\to \Rat(\shc)
\qtq \tens\;\cl \Rat(\shc)\times\shc\to \Rat(\shc).$$

\Def\label{def:ratcent}
A {\em rational center} in $\shc$ is a triple $(\Ma,\phi,\Rmat_\Ma)$
of $\Ma\in\Aff(\shc)$, an additive map $\phi\cl\Idx\to\Z$ and an isomorphism 
$$\Rmat_\Ma(X)\cl q^{\tens\,\phi(\la)}\tens\Ma\tens X\isoto X\tens\Ma$$
in $\Rat(\shc)$ functorial in $X\in\shc_\la$ such that the diagrams
$$
\xymatrix@C=15ex
{q^{\tens\,\phi(\la+\mu)}\tens\Ma\tens X\tens Y\ar[r]^{\Rmat_\Ma(X)\tens Y}\ar[dr]_{\Rmat_\Ma(X\tens Y)}&q^{\tens\,\phi(\mu)}\tens X\tens\Ma\tens Y\ar[d]^{X\tens\Rmat_\Ma(Y)}\\
&X\tens Y\tens\Ma}$$
and $$\xymatrix@C=12ex{\Ma\tens\one\ar[r]^{\Rmat_\Ma(\one)}\ar[dr]^-{\sim}&\one\tens \Ma\ar[d]^\bwr\\
&\Ma}
$$
commute in $\Rat(\shc)$ for any $X\in\shc_\la$ and $Y\in\shc_\mu$ ($\la,\mu\in\Idx$).
\edf
Note that the commutativity of the bottom diagram is a consequence of the one of the top.

{\em In the sequel, we neglect grading shifts if there is no afraid of confusion.}
For example if we say that $f\cl M\to N$ is a morphism,
it means that $f\in\HOM(M,N)_d$ for some $d\in\Z$.
If we want to emphasize $d$, we say that $f$ is a morphism of homogeneous degree $d$.

\Def[{\cite[Definition 6.7]{ref}}]\label{def:affinization}
An {\em affinization} $\Ma$ of $M\in\shc$ (of degree $d\in\Z_{>0}$) is an \afn 
$(\Ma,\z)$ endowed with a rational center $(\Ma,\Rmat_\Ma)$ and an isomorphism
$\Ma/\z\Ma\simeq M$ such that $\deg(\z)=d$. 
\edf

We sometimes simply write $\Ma$ for affinization if no confusion arises.  

\Prop[{\cite[Proposition 6.9]{ref}}]\label{pro:rren}
Let $(\Ma,\Rmat_{\Ma})$ be a rational center in $\shc$, and let $L\in\shc$.
Assume that $\Ma$ and $L$ do not vanish.
Then there exist $m\in\Z$ and  a morphism
$\Rre_{\Ma,L}\cl \Ma\tens L\to L\tens \Ma$ in $\Proc(\cor[z],\shc)$ such that
\bna
\item $\Rre_{\Ma,L}$ is equal to $z^m\Rmat_{\Ma}(L)\cl \Ma\tens L\to L\tens \Ma$ in $\Rat(\shc)$,
\setcounter{myc}{\value{enumi}}
\item $\Rre_{\Ma,L}\vert_{z=0}\cl (\Ma/z\Ma)\tens L\to L\tens (\Ma/z\Ma)$
does not vanish.
\ee
\addtocounter{myc}{1}
Moreover such an integer $m$ and an $\Rre_{\Ma,L}$ are unique.

Similarly,  there exist $m\in\Z$ and  a morphism $\Rre_{L,\Ma}\cl L\tens\Ma\to \Ma\tens L$ in $\Proc(\cor[z],\shc)$ such that  
\bna
\setcounter{enumi}{\value{myc}}
\item $\Rre_{L,\Ma}$ gives $z^m\Rmat_{\Ma}(L)^{-1}\cl L\tens \Ma\to \Ma\tens L$ in $\Rat(\shc)$,
\item $\Rre_{L,\Ma}\vert_{z=0}\cl L\tens (\Ma/z\Ma)\to (\Ma/z\Ma)\tens L$
does not vanish.
\ee

\enprop

When $\deg(\z)=\deg(\zN)$, since $\Rre_{\Ma,\Na}$ commutes with $\zM$ and $\zN$, it induces a morphism in $\Aff(\shc)$
\eq \label{eq: induced_r}
 \Ma \tens_z \Na \To[\bRre_{\Ma,\Na}] \Na\tens_z \Ma\,,
\eneq
which is  denoted  by $\bRre_{\Ma,\Na}$. Here, $z$ acts on $\Ma$ and $\Na$ by $\z$ and $\z[\Na]$, respectively.

\smallskip
We say that a simple object $M$ in a monoidal abelian category is called \emph{real} if $M\tens M$ is simple. 

\Def\label{def:affreal}
We say that a simple object $M\in \shc$ is \emph{\afr} if $M$ is real and there is an affinization $(\Ma,z)$ of $M$.
If $\deg(z)=d$, we say that $M$ is \emph{\afr of degree $d$}.
\edf

\subsection{Quasi-rigid Axiom}\label{subsec:simplehd}

Recall the duality in monoidal categories.

\Lemma\label{lem:MNDM}
Assume that $\shc$ is rigid.
If $(M,N)$ is a $\La$-definable pair of simple modules, then
$(N,\D M)$ is also $\La$-definable and
$\La(M,N)=\La(N,\D M)$.
Moreover $\rmat{N,\D M}$ is obtained as the composition
$$N\tens \D M\to \D M\tens M\tens N\tens\D M\To[\rmat{M,N}]
\D M\tens N\tens M\tens \D M \to\D M\tens N.$$
\enlemma
\Proof
It follows from
$$\HOM(M\tens N,N\tens M)\simeq\HOM(N, \D M\tens N\tens M)
\simeq\HOM(N\tens \D M,\D M\tens N).$$
\QED
\Def[{\cite{ref}}]
Let $(\sha,\tens)$ be a monoidal category.
We say that $\sha$ is a {\em \KO} monoidal category if it satisfies:
\bna
\item $\sha$ is abelian and $\tens$ is bi-exact,
\item for any $L,M,N\in\sha$, $X\subset L\tens M$ and
$Y\subset M\tens N$ such that $X\tens N\subset L\tens Y\subset L\tens M\tens N$, there exists $K\subset M$ such that
$X\subset L\tens K$ and $K\tens N\subset Y$,
\item for any $L,M,N\in\sha$, $X\subset M\tens N$ and
$Y\subset L\tens M$ such that $L\tens X\subset Y\tens N\subset L\tens M\tens N$, there exists $K\subset M$ such that
$X\subset K\tens N$ and $L\tens K\subset Y$.
\ee
\edf

As for the following statements, see \cite{KKOP23}, \cite[Sect.6]{ref},
which will be needed in \S\,\ref{sec:cat-C}:

\Lemma\label{lem:monoepi}
Let $M_k$ be non-zero object $( k=1,2,3)$ of a quasi-rigid monoidal category
$\sha$, and let $\vphi_1\cl  L \to  M_1\tens M_2$ and
$\vphi_2\cl  M_2 \tens M_3 \to L'$ be non-zero morphisms.  Assume further that $M_2$ is simple.
Then the composition
\eqn
L \tens M_3 \To[\vphi_1 \tens M_3] M_1 \tens M_2 \tens M_3 \To[M_1 \tens\vphi_2] M_1 \tens L'
\eneqn
does not vanish.
\enlemma

\Lemma[{\cite[Lemma~6.16]{ref}}]\label{lem:q-rigid>rigid}
An abelian rigid monoidal category is \KO.
\enlemma
Indeed, the notion ``quasi-rigid" is rather natural since the monoidal categories appearing in this article  
are all quasi-rigid as seen in the example below.
 
\Ex\label{ex:q-rigid}
By \cite[Lemma~3.1]{KKKO15} we find that
the monoidal categories $R\gMod$, $R\gmod$ and $\Cw$ introduced later are quasi-rigid.
By Lemma~\ref{lem:q-rigid>rigid} above, we see that the category $\tCw$ is also quasi-rigid
since it is rigid \cite{KKOP22}.
\enex

Recall that for objects $M,~N$ in a monoidal category $(\shc,\tens)$, we denote by $M\hconv N$(resp.\ $M\sconv N$) the
head (resp.\ socle) of $M\tens N$.

\Prop[{\cite{KKKO15}, \cite[Proposition~6.18, 6.19]{ref}}] \label{prop:simplehd}
Assume that $\shc$ is \KO.
Let $M$ be an \afr object of $\shc$ and
$N$ a simple object of $\shc$.
Then, we obtain 
\bnum
\item
$M\tens N$ and $N\tens M$ have simple heads and simple socles,
\item the pair $(M,N)$ is $\de$-definable \ro see Definition~\ref{def:rmat}\/\rf
  and
  \eqn
  &&\HOM_\shc(M\tens N,N\tens M)={\bf k}\rmat{M,N}, \q
\HOM_\shc(N\tens M,M\tens N)={\bf k}\rmat{N,M},\\
&&M\hconv N\simeq\Im(\rmat{M,N})\simeq N\sconv M\qtq
N\hconv M\simeq\Im(\rmat{N,M})\simeq M\sconv N
\eneqn
up to grading shifts,
\item
$\La(M,M\htens N)=\La(M,N)$ and $\La(N\htens M,M)=\La(N,M)$,
\item for any simple subquotient $S$ of 
the radical $\Ker(M\tens N\to M\htens N)$,
we have $\La(M,S)<\La(M,N)$,\label{item:srad}
\item for any simple subquotient $S$ of $(M\tens N)/(M\stens N)$,
  we have $\La(S,M)<\La(N,M)$, 
\item for any simple subquotient $S$ of 
the radical $\Ker(N\tens M\to N\htens M)$,
we have $\La(S,M)<\La(N,M)$,
\item for any simple subquotient $S$ of $(N\tens M)/(N\stens M)$,
we have $\La(M,S)<\La(M,N)$
\item 
we have $\END_\shc(M\tens N)=\cor\id_{M\tens N}$
and $\END_\shc(N\tens M)=\cor\id_{N\tens M}$,
\item \label{it:de tens}
$\de(M, M^{\tens n}\hconv N)\le\max(\de(M,N)-n,0)$ for any $n\in\Z_{\ge0}$.
\ee
In particular,
$M\htens N$, as well as $M\stens N$,  appears only once
in the composition series of $M\tens N$
\ro up to a grading\rf.

\enprop

The analogous statement to the following lemma appears in 
\cite[Corollary 3.13]{KKKO15}, \cite[Lemma 2.7]{PBWQA}.

\Lemma\label{lem:invhvonv} Assume that $\shc$ is quasi-rigid.
Let $L$ be an \afr simple object of $\shc$.
Then for any simple $M\in\shc$, we have
$$(L\htens M)\htens \D L\simeq M\qtq
\D^{-1}L\htens (M\htens L)\simeq M.$$
\enlemma
Remark that this lemma is important to define the crystal structure on $\Irr(\tCw)$
later.

\Prop\label{prop:ddd}
Assume that $\shc$ is quasi-rigid.
Let $L$ be a real simple object in $\shc$ with an affinization of degree
$2d\in \Z_{>0}$.
Then, we have
\bnum
\item
$\de(L,X)\in d\Z_{\ge0}$ for any simple $X\in\shc$,
\item
$\La(L,M)+\La(L,N)-\La(L,S)\in 2d\Z_{\ge0}$
for any simple $M$, $N$ and any simple subquotient $S$ of $M\tens N$.
\ee
\enprop
\Proof
Let $(\Laa,\Rmat_\Laa)$ be an affinization with a spectral parameter $z$ with $\deg(z)=2d$.
For a simple $X\in\shc$, there exists $n_X\in\Z$ such that
$z^{n_X}\Rmat_\Laa(X)\cl \Laa\tens X\to X\tens\Laa$ is a morphism in $\Aff(\shc)$
and does not vanish at $z=0$. Its degree is $\La(L,X)$.
Similarly, there exists $m_X\in\Z$ such that
$z^{m_X}\Rmat_\Laa(X)^{-1}\cl X\tens\Laa\to\Laa\tens X$  is a morphism in $\Aff(\shc)$
and does not vanish at $z=0$. Its degree coincides with $\La(X,L)$.

\snoi(i)\ 
The composition
$z^{m_X+n_X}\id_{\Laa\tens X}$ is an endomorphism of $\La\tens X$ in $\Aff(\shc)$
of degree $2(m_X+n_X)d\ge0$,
Since it coincides with $\La(L,X)+\La(X,L)=2\de(L,X)$.
Hence $\de(L,X)\in\Z_{\ge0}d$.

\snoi
(ii)\
$z^{n_M+n_N}\Rmat_\Laa(M\tens N)\cl \Laa\tens (M\tens N)\to
(M\tens N)\tens \Laa$ is a morphism in $\Aff(\shc)$.
Hence it induces a morphism
$z^{n_M+n_N}\Rmat_\Laa(S)\cl \Laa\tens S\to
S\tens \Laa$ in $\Aff(\shc)$.
Hence there exists $s\in\Z_{\ge0}$ such that
$z^{n_M+n_N-s}\Rmat_\Laa(S)\cl \Laa\tens S\to
S\tens \Laa$ is a morphism in $\Aff(\shc)$ and does not vanish.
Since we have
\eqn
\La(L,S)&&=\deg\bl z^{n_M+n_N-s}\Rmat_\Laa(S)\br
=\deg\bl z^{n_M+n_N-s}\Rmat_\Laa(M\tens N)\br\\
&&=\La(L,M)+\La(L,N)-\deg(z^s),\eneqn
which implies the desired result.
\QED

 \Rem
It follows from Example~\ref{ex:q-rigid}
that we can apply the results of this subsection
to the categories $R\gmod$, $\Cw$ and $\tCw$.
\enrem
%
\section{Quiver Hecke algebras and their modules}
\subsection{Quiver Hecke algebras } In the sequel, $\bR$ is a base field.
Let $ \bl\cmA,\wlP,\Pi,\Pi^\vee,(\cdot,\cdot) \br $ be a Cartan datum.
A { \em quiver Hecke datum} is a family of polynomials
$(\qQ_{i,j}(u,v) \in \bR[u,v])_{i,j\in I}$ of the form

\begin{align}
\qQ_{i,j}(u,v) =\bc
                   \sum\limits
_{p(\alpha_i , \alpha_i) + q(\alpha_j , \alpha_j) = -2(\alpha_i , \alpha_j) } t_{i,j;p,q} u^pv^q &
\text{if $i \ne j$,}\\[3ex]
0 & \text{if $i=j$,}
\ec\label{eq:Q}
\end{align}
such that $t_{i,j;-a_{ij},0} \in  \bR^{\times}$ and
$\qQ_{i,j}(u,v)= \qQ_{j,i}(v,u) \quad \text{for all} \ i,j\in I.$
We set
\begin{align*}
\bQ_{i,j}(u,v,w)\seteq\dfrac{ \qQ_{i,j}(u,v)- \qQ_{i,j}(w,v)}{u-w}\in \bR[u,v,w].
\end{align*}

For $\beta\in \rlQ_+$,  the set
$I^\beta\seteq  \Bigl\{\nu=(\nu_1, \ldots, \nu_n ) \in I^n \bigm| \sum_{k=1}^n\alpha_{\nu_k} = \beta \Bigr\}$ is stable under the action
of the symmetric group $\mathfrak{S}_n = \langle s_k \mid k=1, \ldots, n-1 \rangle$ given by place permutations.

The height of $\beta=\sum_{i\in I} b_i \al_i\in\rtl$ is given by $\height{\beta}\seteq\sum_{i\in I} |b_i|$.
\Def\label{df:basic}
Let $\beta\in\rlQ_+$ with $\height{\beta}=n$.
The {\em quiver Hecke algebra} $R(\beta)$ associated with a Cartan datum $ \bl\cmA,\Pi,\wlP,\Pi^\vee,(\cdot,\cdot) \br $ and a quiver Hecke datum $(\qQ_{i,j}(u,v))_{i,j\in I}$
is the $\bR$-algebra generated by
$$
\{e(\nu) \mid \nu \in I^\beta \}, \; \{x_k \mid 1 \le k \le n \},
 \; \text{and} \; \{\tau_l \mid 1 \le l \le n-1 \}
$$
subject to the defining relations:
\eqn
&& e(\nu) e(\nu') = \delta_{\nu,\nu'} e(\nu),\ \sum_{\nu \in I^{\beta}} e(\nu)=1,\
x_k e(\nu) =  e(\nu) x_k, \  x_k x_l = x_l x_k,\\
&& \tau_l e(\nu) = e(s_l(\nu)) \tau_l,\  \tau_k \tau_l = \tau_l \tau_k \text{ if } |k - l| > 1, \\[1ex]
&&  \tau_k^2 = \sum_{\nu\in I^\beta}\qQ_{\nu_k, \nu_{k+1}}(x_k, x_{k+1})e(\nu), \\[5pt]
&& \tau_k x_l - x_{s_k(l)} \tau_k =
\bl\delta(l=k+1)-\delta(l=k)\br
\sum_{\nu\in I^\beta,\ \nu_k=\nu_{k+1}}e(\nu),\\
&&\tau_{k+1} \tau_{k} \tau_{k+1} - \tau_{k} \tau_{k+1} \tau_{k}
=\sum_{\nu\in I^\beta,\ \nu_k=\nu_{k+2}}
\bQ_{\,\nu_k,\nu_{k+1}}(x_k,x_{k+1},x_{k+2}) e(\nu).
\eneqn
\edf

The algebra $R(\beta)$ is equipped with the $\Z$-grading given by
\eq
&&\deg(e(\nu))=0, \quad \deg(x_k e(\nu))= ( \alpha_{\nu_k} ,\alpha_{\nu_k}), \quad  \deg(\tau_l e(\nu))= -(\alpha_{\nu_{l}} , \alpha_{\nu_{l+1}}).
\eneq

We denote by 
$R(\beta) \gMod$  the category of graded $R(\beta)$-modules with degree preserving homomorphisms.
We  write $R(\beta)\gmod$ for the full subcategory of $R(\beta)\gMod$ consisting of the graded modules which are  finite-dimensional over $\bR $, and 
 $R(\beta)\proj$ for  the full subcategory of $R(\beta)\gMod$ consisting of finitely generated  projective graded $R(\beta)$-modules.
We set $R\gMod \seteq \bigoplus_{\beta \in \rlQ_+} R(\beta)\gMod$, $R\proj \seteq \bigoplus_{\beta \in \rlQ_+} R(\beta)\proj$, and  $R\gmod \seteq \bigoplus_{\beta \in \rlQ_+} R(\beta)\gmod$.
The trivial $R(0)$-module $\cor$ with degree 0 is denoted by $\trivialM$.
For simplicity, we write ``a module" instead of ``a graded module''.
We define the grading shift functor $q$
by $(qM)_k = M_{k-1}$ for a graded module $M = \bigoplus_{k \in \Z} M_k $.
For $M, N \in R(\beta)\gMod $, $\Hom_{R(\beta)}(M,N)$ denotes the space of degree preserving module homomorphisms.
We define
\[
\HOM_{R(\beta)}( M,N ) \seteq \bigoplus_{k \in \Z} \HOM_{R(\beta)}(M, N)_k
\qt{with $\HOM_{R(\beta)}(M, N)_k=\Hom_{R(\beta)}(q^{k}M, N)$,}\]
and set $ \deg(f) \seteq k$ for $f \in \HOM_{R(\beta)}(M, N)_k$.
When $M=N$, we write $\END_{R(\beta)}( M ) = \HOM_{R(\beta)}( M,M)$.
We sometimes write $R$ for $R(\beta)$ in $\HOM_{R(\beta)}( M,N )$ for simplicity.

 For $\beta,\beta'\in\prtl$, set
$
e(\beta, \beta') \seteq \sum_{\nu \in I^\beta, \nu' \in I^{\beta'}} e(\nu\ast\nu'),
$
where $\nu\ast\nu'$ is the concatenation of $\nu$ and $\nu'$. 
Then there is an injective ring homomorphism
$$R(\beta)\tens R(\beta')\to e(\beta,\beta')R(\beta+\beta')e(\beta,\beta')$$
given by
$e(\nu)\tens e(\nu')\mapsto e(\nu,\nu')$,
$x_ke(\beta)\tens 1\mapsto x_ke(\beta,\beta')$,
$1\tens x_ke(\beta')\mapsto x_{k+\height{\beta}}e(\beta,\beta')$,
$\tau_ke(\beta)\tens 1\mapsto \tau_ke(\beta,\beta')$ and
$1\tens \tau_ke(\beta')\mapsto \tau_{k+\height{\beta}}e(\beta,\beta')$.
For $a\in R(\beta)$ and $a'\in R(\beta')$, the image of $a\tens a'$ is sometimes denoted by $a\etens a'$.

For $M \in R(\beta)\gMod$ and $N \in R(\beta')\gMod$, we set
$$
M \conv N \seteq R(\beta+\beta') e(\beta, \beta') \otimes_{R(\beta) \otimes R(\beta')} (M \otimes N).
$$
For $u\in M$ and $v\in N$, the image of
$u\tens v$ by the map $M\tens N\to M\conv N$ is sometimes denoted by
$u\etens v$. 
We also write $M\etens N\subset M\conv N$ for the image of $M\tens N$ in $M\conv N$.

For $\al,\beta\in\prtl$, let $X$ be an $R(\al+\beta)$-module.
Then $e(\al,\beta)X$ is an $R(\al)\tens R(\beta)$-module.
We denote it
by $\Res_{\al,\beta}X.$
We have
\eq \label{eq:adjoints}
&&\ba{rl}
\Hom_{R(\al)\tens R(\beta)}\bl M\tens N,\Res_{\al,\beta}(X)\br
&\simeq\Hom_{R(\al+\beta)}(M\conv N,X),\\
\Hom_{R(\al)\tens R(\beta)}\bl\Res_{\al,\beta}(X), M\tens N\br
&\simeq\Hom_{R(\al+\beta)}(X, q^{(\al,\beta)}N\conv M)\ea
\eneq
for any $R(\al)$-module $M$, any $R(\beta)$-module $N$ and
any $R(\al+\beta)$-module $X$.

\medskip
Recall that we denote by $M \hconv N$ the head of $M \conv N$ and by $M \sconv N$ the socle of $M \conv N$.
We say that simple $R$-modules $M$ and $N$ \emph{strongly commute} if $M \conv N$ is simple.  Recall that a simple $R$-module
$L$ is \emph{real} if $L$ strongly commutes with itself, that is, 
 $L\conv L$ is simple. 
Note that if $M$ and $N$ strongly commute,  then $M$ and $N$ commute, i.e., $M\conv N \simeq N \conv M$ up to a grading shift.

For $i\in I$, the functors $\E_i$ and $\F_i$ are defined by
\eqn&&
\ba{rl}\E_i(M)& = e(\alpha_i, \beta-\al_i) M \in  R(\beta-\al_i) \gMod\\
 \F_i(M) &= R(\alpha_i) \conv M \in  R(\beta+\al_i) \gMod 
\ea\hs{5ex}\qt{for any $R(\beta)$-module $M$.}
\eneqn
For $i\in I $ and $n\in \Z_{>0}$, let $L(i)$ be the simple $R(\alpha_i)$-module concentrated on  degree 0 and
 $P(i^{n})$  the indecomposable  projective $R(n \alpha_i)$-module
 whose head is isomorphic to $L(i^n) =\ang{i^n}
 \seteq q_i^{\frac{n(n-1)}{2}} L(i)^{\conv n}$,
where $q_i\seteq q^{(\al_i,\al_i)/2}$. 
We use the both notations $L(i^n)$ and $\ang{i^n}$ in the sequel.

Then, for $M\in  R(\beta) \gMod $, we define
\eqn&&
\ba{ll}
\E_i^{(n)} M& \seteq \HOM_{R(n\alpha_i)} \bl P(i^{n}),\, e(n\alpha_i, \beta - n\alpha_i) M\br \in  R(\beta-n\al_i) \gMod,\\
\F_i^{(n)} M& \seteq  P(i^{n}) \conv M \in  R(\beta+n\al_i) \gMod. 
\ea
\eneqn
For $i \in I$ and a non-zero  $M \in R(\beta)\gMod$, we define
\eq
&&\ba{ll}\wt(M) = - \beta, \akew[3ex]&\ep_i(M) = \max \{ k \ge 0 \mid E_i^k M \not\simeq 0 \}, \\&
\ph_i(M) = \ep_i(M) + \langle h_i, \wt(M) \rangle.\ea
\eneq
We can also define $\E_i^*$, $\F_i^*$, $\ep^*_i$, etc.\ in the same manner as above if we replace  $e(\alpha_i, \beta-\al_i)$, $R(\alpha_i)\conv -$,
etc.\  with 
$e(\beta-\al_i, \alpha_i)$, $- \conv R(\alpha_i)$, etc.

We denote by $K(R\proj)$ and $K(R\gmod)$ the Grothendieck groups of $R\proj$ and $R\gmod$, respectively.

Let $\Irr(R\gmod)$ be the set of the equivalence classes of simple
modules in $R\gmod$ up to grading shifts.
Namely, for simple modules $M$ and $N$,
we identify $M$ and $N$ in $\Irr(R\gmod)$ if $M\simeq q^nN$ for some $n\in\Z$.

For a simple module $M\in R\gmod$, define
\eq
&&\hs{5ex}\tE_i(M)\seteq q_i^{\eps_i(M)-1}{\rm hd}\,(\E_iM)\simeq q_i^{1-\eps_i(M)}{\rm soc}\,(\E_iM)\qtq
\tF_i(M)\seteq q_i^{\eps_i(M)}\ang{i}\hconv M.
\label{eq:eit-fit}
\eneq
We set
$$\tEm_i(M)\seteq\tE_i^{\ms{4mu}\eps_i(M)}M\simeq\E_i^{(\eps_i(M)}M.$$
The following theorem is due to Lauda-Vazirani. 
\Th[{\cite{LV11}}]\label{thm:L-V} 
The 6-tuple $\bl\Irr(R\gmod), \wt,\{\ep_i\}_{i\in I},\{\ph_i\}_{i\in I},
\{\tE_i\}_{i\in I},\{\tF_i\}_{i\in I}\br$ is a crystal and 
there exists an isomorphism of crystals:
\[
\Psi\cl\Irr(R\gmod)\mapright{\sim} B(\ify)\quad (\triv\mapsto u_\ify).
\]
\end{thm}

\Rem[{cf.\ \cite{KR10}}]\label{rem:modulenot}
Let $\beta\in\prtl$ and let $S$ be a subset of $I^\beta$ satisfying the following conditions:
\bnam
\item
  if $\mu=(\mu_1,\ldots,\mu_m)\in S$, then $\mu_k\not=\mu_{k+1}$
  for any $k$ such that $1\le k <n$,
\item for any $\mu=(\mu_1,\ldots,\mu_m)\in S$
  and $k$ such that $1\le k <n$,
$s_k\mu\in S$ if $(\al_{\mu_k},\al_{\mu_{k+1}})=0$,

\item
    if $\mu=(\mu_1,\ldots,\mu_m)\in S$ and $\mu_k=\mu_{k+2}$
    for $k$ such that $1\le k\le n-2$, we have
    $\ang{h_{\mu_k},\al_{\mu_{k+1}}}<-1$.
    \ee
    Then we can easily check that
    the vector space $M\seteq\soplus_{\mu\in S}\cor\,u_\mu$ with basis $\st{u_\mu}_{\mu\in S}$
    has a structure of
    $R(\beta)$-module such that $\deg(u_\mu)=0$, $e(\mu)u_\mu=u_\mu$
    $x_k\vert_M=0$ and
    $\tau_k u_\mu=\bc u_{s_k\mu}&\text{if $(\al_{\mu_k},\al_{\mu_{k+1}})=0$,}\\
    0&\text{otherwise.}\ec$
    
\noi
    Here we assume that $Q_{i,j}(u,v)=1$ when $(\al_i,\al_j)=0$ without loss of generality.
    \noi
    If $S$ is a minimal one among the non-empty subsets satisfying
    the above conditions,
    then $M$ is a simple module.
    Sometimes, $M$ is denoted by $\ang{\nu_1,\ldots,\nu_n}$ taking a member
    $\ang{\nu_1,\ldots,\nu_n}$ of $S$.
    With this terminology,
   $\ang{i}$ and $\ang{i,j}$ ($i\not=j$) are such examples.
We have $\ang{i}=L(i)$, Thus, in the sequel, we use frequently $\ang{i}$ instead of $L(i)$ for an abbreviation.
        \enrem


\subsection{Duality in $R\gmod$}\label{subsec:dual}

For $M \in R(\beta)\gmod$, we set $M^\star \seteq  \HOM_{\bR}(M, \bR)$ with the $R(\beta)$-action given by
$$
(a \cdot f) (u) \seteq  f(\rho(a)u), \quad \text{for  $f\in M^\star$, $a \in R(\beta)$ and $u\in M$,}
$$
where $\rho$ is the antiautomorphism of $R(\beta)$ which fixes the generators
 $e(\nu)$, $x_k$, $\tau_l$. 
 We say that $M\in R\gmod$ is \emph{self-dual}
 if $M \simeq M^\star$ in $R\gmod$.
 The following lemma is well-known.
 \Lemma[{cf.\ \cite{KKOP21}}]\label{self-n}
 For any simple $M\in R\gmod$, there exists a unique
 $n\in\Z$ such that
 $q^nM$ is self-dual.
 \enlemma
 It means that any simple module in $R\gmod$ is isomorphic to its dual 
 up to a grading shift.

We have
$$(M\conv N)^\star\simeq q^{(\wt M,\wt N)}N^\star\conv M^\star
\qt{functorially in $M,N\in R\gmod$.}$$

\Lemma\label{lem:Laspl}
Let $M$ and $N$ be simple modules in $R\gmod$.
\bnum
\item Any simple quotient of $M\conv N$ appears as a simple submodule of $N\conv M$ \ro up to grading shifts\/\rf. Conversely, any simple submodule of $M\conv N$
  appears as a simple quotient of $N\conv M$.
\item
  There exists a morphism $f\cl M\conv N\to N\conv M$ such that
  $\Im(f)$ is simple.
\item
  If $(M,N)$ is a $\La$-definable pair \ro see {\rm Definition~\ref{def:rmat}}\/\rf,
  then $M\conv N$ has a simple head and $N\conv M$ has a simple socle.
  Moreover, we have
  $$M\hconv N\simeq\Im(\rmat{M,N})\simeq N\sconv M.$$  
\ee
\enlemma
\Proof
If $M\conv N\epito S$ is a simple quotient of $M\conv N$, then,
applying ${}^\star$, we obtain
$$S\simeq S^\star\monoto (M\conv N)^\star\simeq N^\star\conv M^\star\simeq N\conv M.$$
Hence $S$ appears as a simple submodule of $N\conv M$.
Other statements are now obvious.
\QED

\subsection{R-matrices} \label{subSec: R-matrices}
Let $\beta \in \rlQ_+$ with $m =  \Ht(\beta)$. For  $k=1, \ldots, m-1$ and $\nu \in I^\beta$, the \emph{intertwiner} $\varphi_k \in R(\beta) $ is defined by 
\eq
\varphi_k e(\nu) =
\bc
 \bl\tau_k(x_k-x_{k+1})+1\br e(\nu) 
& \text{ if } \nu_k = \nu_{k+1}, 
 \\
 \tau_k e(\nu) & \text{ otherwise.}
\ec\label{def:intertwiner}
\eneq
Since $\st{\vphi_k}_{1\le k<m}$ satisfies the braid relation, we can define
$\vphi_w$ for any element $w$ of the symmetric group $\sym_m$ of degree $m$.

For $m,n \in \Z_{\ge 0}$, we set $w[m,n]$ to be the element of $\sg_{m+n}$ such that
$$
w[m,n](k) \seteq
\left\{
\begin{array}{ll}
 k+n & \text{ if } 1 \le k \le m,  \\
 k-m & \text{ if } m < k \le m+n.
\end{array}
\right.
$$

For $M\in \Modg\bl R(\beta)\br $ and $N\in \Modg\bl R(\gamma)\br$,
the map $M \otimes N \rightarrow N \conv M$ given by
 $$u \otimes v \mapsto \varphi_{w[\height{\gamma},\height{\beta}]}(v \etens u)$$
is
$R(\beta)\otimes R(\gamma)$-linear, and hence
it extends to an $R(\beta+\gamma)$-module homomorphism (neglecting a grading shift)
$$
\Runi_{M,N}\cl  M\conv N \longrightarrow N \conv M.
$$
We call it the {\em universal R-matrix}.

\medskip
For $\beta \in \rlQ_+$ and $i\in I$,  let $\mathfrak{p}_{i, \beta}$ be an element in the center $Z(R(\beta))$ of $R(\beta)$ given by
\begin{align} \label{Eq: def of p}
\mathfrak{p}_{i, \beta}  \seteq \sum_{\nu \in I^\beta} \Bigl(\hs{1ex}  \prod_{a \in \{1, \ldots, \Ht(\beta) \},\ \nu_a=i} x_a \Bigr) e(\nu)\in Z(R(\beta)).
\end{align}

\Def[\cite{KP18}]\label{def:affinization1}
For a simple module $M\in R(\beta)\gmod$,
an  \emph{affinization of $M$} is a pair  $(\Ma,z_{\Ma})$
of  an $R(\beta)$-module $\Ma$ and an endomorphism $z_{\Ma}$ of $\Ma$
with degree $d_{\Ma} \in \Z_{>0}$ such that 
\eq \label{eq:oldaff}
&&\hs{2ex}\parbox{75ex}{
\begin{enumerate}[\rm (i)]
\item $\Ma / z_{\Ma} \Ma \simeq M$,
\item $\Ma$ is a finitely generated free module over the polynomial ring $\bR[z_{\Ma}]$,
\item $\mathfrak{p}_{i,\beta} \Ma \ne 0$ for all $i\in I$.\label{it:nonzeroP}
\end{enumerate}
}\eneq

We recall that a real simple module is {\em \afr} if it admits an affinization.
\edf

Note that, if $(\Ma,z_{\Ma})$ is an affinization of $M$,
then $(\Ma,\Runi_{\Ma})$ is an affinization of $M$ in the sense of
Definition~\ref{def:affinization} (\cite{KP18}). Here, $\Runi_{\Ma}(X)\cl \Ma\conv X\to X\conv\Ma$ is the universal R-matrix,
and 
we identify $\Ma$ with $\proolim[{k}]\Ma/z_\Ma^k\Ma\in\Aff(R\gmod)$. 
Hence we can apply the results in \S\,\ref{subsec:simplehd}, in particular Proposition~\ref{prop:simplehd}.

\Ex\label{ex:i-affreal}
The simple $R$-module $L(i)=\ang{i}$ is \afr.
\enex

\subsection{Normal sequences} 

Recall the definition of $\La(M,N)$ and $\de(M,N)$
in  Definition~\ref{def:rmat}. 
We define
\eqn
&&\tLa(M,N) \seteq   \frac{1}{2} \bl \La(M,N) + (\wt(M), \wt(N)) \br
\eneqn
for a $\La$-definable pair $(M,N)$ of simple modules in $R\gmod$.

Note that  
$ \de(M,N)$ and $ \tLa(M,N)$ are non-negative integers (\cite[Lemma 3.11]{KKOP21}) as soon as one of $M$ and $N$ is \afr.

\bigskip
Let $\shc$ be a quasi-rigid monoidal category satisfying \eqref{cond:exactmono}.

\Def\label{def:normal-seq}
Let $(M_1,\ldots,M_r)$ be a sequence  of simple objects
in $\shc$.
\bnum
\item
  We say that $(M_1,\ldots,M_r)$ is {\em almost \afr} if
  $M_k $  is \afr except for  possibly one $k$.
  \item
A sequence $(M_1,\ldots,M_r)$
of simple objects is called a \emph{normal sequence} 
if  $(M_j,M_k)$ is $\La$-definable for any $1\le j<k\le r$ and
the composition of r-matrices 
\eqn \rmat{M_1,\ldots,M_r}\seteq
&&(\rmat{M_{r-1},M_r})  \circ \cdots \circ (\rmat{M_2,M_r}\circ \cdots \circ \rmat{M_2,M_3})  \circ (\rmat{M_1,M_r} \circ \cdots  \circ \rmat{M_1,M_2}) 
\\
  &&\cl M_1\conv \cdots \conv M_r \longrightarrow M_r \conv \cdots \conv  M_1
\eneqn
does not vanish.
\ee
\edf

\Lemma[{\cite[\S\,2.2]{KKOP22}}]
If $(M_1,\ldots,M_r)$ be an almost \afr normal sequence in $\shc$,
then
$\Im (\rmat{M_1,\ldots,M_r})$ is simple, and it is isomorphic to the head of $M_1\conv \cdots\conv  M_r$
and to the socle of $M_r\conv \cdots\conv  M_1$.
\enlemma

\Lemma [{\cite[Lemma 2.7, Lemma 2.8]{KK19}}] \label{lem:normal}
 Let $(L_1,\ldots,L_r)$ be an almost \afr sequence  of  simple modules  in 
$\shc$.
Then  we have
\bnum
\item
Assume that $L_1$ is \afr.
Then 
$(L_1,\ldots,L_r)$    is a normal sequence
if and only if
$(L_2,\ldots,L_r)$ is a normal sequence and 
$$\La(L_1, \hd(L_2\conv\cdots \conv L_r)) = \sum_{2\le j\le r} \La(L_1,L_j).$$
\item
  Assume that $L_r$ is \afr.
  Then, $(L_1,\ldots,L_{r})$ is a normal sequence if and only if
    $(L_1,\ldots,L_{r-1})$ is a normal sequence  and  
$$\La(\hd(L_1\conv\cdots \conv L_{r-1}), L_r) = \sum_{1\le j\le r-1} \La(L_j,L_r).$$
\end{enumerate}
\enlemma

\Lemma \label{lem:normalD}
Let $(L,M,N)$ be an almost \afr sequence of simples in $\shc$.
\bnum
\item If $L$ and $M$ commute, then $(L,M,N)$ is a normal sequence.
\item Assume that $\shc$ is rigid.
  Then, $(L,M,N)$ is a normal sequence if and only if $(M,N,\D L)$
  is normal.
\item Assume that $\shc$ is rigid and that $\D L$ and $N$ commute.
  Then, $(L,M,N)$ is a normal sequence.
  \ee
  \enlemma
  \Proof
  (i) follows from $\La(L\hconv M,N)=\La(L,N)+\La(M,N)$ and Lemma~\ref{lem:normal} (ii).

    \snoi
    (ii)
    Let $$\rmat{L,M,N}\cl L\tens M\tens N\To[\rmat{M,N}]
    L\tens N\tens M\To[\rmat{L,N}]N\tens L\tens M\To[\rmat{L,M}]
    N\tens M\tens L$$
    be the R-matrix.
    We see by Lemma~\ref{lem:MNDM} that
    $\rmat{M,N,\D L}$ is obtained as the composition 
        $$\xymatrix@C=15ex{
      M\tens N\tens \D L\ar[d]\ar[r]^{\rmat{M,N,\D L}}&{\D L \tens N\tens M} \\
    \D L\tens L\tens M\tens N\tens \D L\ar[r]^{\D L\tens\rmat{L,M,N}\tens\D L}&
    \D L \tens N\tens M\tens L\tens\D L.\ar[u]}$$
  Hence if $\rmat{M,N,\D L}$ does not vanish, then $\rmat{L,M,N}$ does not vanish.
Similarly
    $\rmat{L,M,N}$ is obtained as the composition
    $$\xymatrix@C=15ex{
      L\tens M\tens N\ar[d]\ar[r]^{\rmat{L,M,N}}&N\tens M\tens L\\
    L\tens M\tens N\tens \D L\tens L\ar[r]^{L\tens\rmat{M,N,\D L}\tens L}&
     L\tens \D L\tens N\tens M\tens L.\ar[u]}$$
  Hence if $\rmat{L,M,N}$ does not vanish, then $\rmat{M,N,\D L}$ does not vanish.

  \snoi
  (iii) follows from (i) and (ii).

    \QED
 
 \Prop[{\cite[Proposition~2.9]{KKOP22}}]\label{prop:Normal} 
Let $(L,M,N)$  be an almost \afr sequence of simple modules in $R\gmod$. 
If $\tLa(L,N)=0$ and one of $L$ and $N$ is \afr, then
$(L,M,N)$ is a normal sequence.
\enprop

\subsection{Head simplicity of convolutions with $L(i)$}
\

\Def
For $i\in I$, $\beta\in\rtlp$ and a simple $R(\beta)$-module $M$, define
$$\de_i(M)\seteq \eps_i(M)+\eps^*_i(M)+\ang{h_i,\wt(M)}.$$
\edf
Recall the following lemma.
\Lemma [{\cite[Corollary 3.8]{KKOP18}}]\label{LmLmtilde}
For $i\in I$, $\beta\in\rtlp$ and a simple module $R(\beta)$-module $M$,
we have
\eqn
\La(L(i), M) &&= (\alpha_i, \alpha_i) \ep_i(M) + (\alpha_i, \wt(M)), \\
\La( M, L(i))&&= (\alpha_i, \alpha_i) \ep^*_i(M) + (\alpha_i, \wt(M)),\\
\tLa(L(i),M)&&=\dfrac{(\al_i,\al_i)}{2}\eps_i(M),\qquad
\tLa(M,L(i))=\dfrac{(\al_i,\al_i)}{2}\eps^*_i(M),\\
\de(L(i),M)&&=\dfrac{(\al_i,\al_i)}{2}\de_i(M).\\
\eneqn
\enlemma

\Prop[{cf.\ \cite{LV11}}]\label{prop:epsBi}
Let $i\in I$, $n\in\Z_{\ge0}$ and let $M$ be a simple module.
\bnum
\item 
If $\de_i(M)=0$,
then we have $L(i)\hconv M\simeq L(i)\conv M \simeq  M\conv L(i) \simeq  M\hconv L(i)$ \ro up to grading shifts\/\rf  
and $\de_i(L(i)\conv M)=0$.

\noi
If $\de_i(M)>0$,
then we have
\eqn
\de_i(L(i)\hconv M)&=&\de_i(M\hconv L(i))=\de_i(M)-1,\\
\eps_i(M\hconv L(i))&=&\eps_i(M),\quad\eps^*_i(L(i)\hconv M)=\eps^*_i(M).
\eneqn
\item We have
\eq
&&\ba{l}
\de_i(L(i^n)\hconv M)=\max\bl\de_i(M)-n, 0\br,\\
\de_i(M\hconv L(i^n))=\max\bl\de_i(M)-n, 0\br.
\ea\label{eq:dei}
\eneq
\item We have
\eqn
\eps_i(M\hconv L(i^n))&&=\max\bl\eps_i(M),\; n-\wt_i (M)-\eps^*_i (M) \br,\\
\eps^*_i(L(i^n)\hconv M)&&=\max\bl\eps^*_i(M),\; n-\wt_i(M)-\eps_i(M)\br.\\
\eneqn
\ee
\enprop

\subsection{Determinantial modules }\label{subsec:determinantial}

We define the partial order $\ble$ on $\wtl$ as follows:
$\la \ble \mu$  for $\la,\mu\in\wtl$ if there exists a sequence  of positive real roots $\beta_1,\ldots,\beta_r$ such that $(\beta_k, s_{\beta_{k+1}} s_{\beta_{k+2}} \cdots s_{\beta_r}\mu) >0$ for all $1\le k\le r$ and $\la=s_{\beta_1} s_{\beta_2} \cdots s_{\beta_r} \mu$. 
We have $\mu-\la\in\prtl$ if $\la \ble \mu$. Hence $\ble$ is a partial order on $\wtl$.

Assume that $\lambda, \mu \in \weyl \Lambda$ for some $\La \in \pwtl$.
Then $\lambda \ble \mu$ if and only if there exist $w,v\in\weyl$ such that
$\la=v\La$, $\mu=v\La$ and $w\ge v$.

For $\La\in\pwtl$ and $\la,\mu\in\weyl\La$ such that $\la\ble \mu$,
there exists a self-dual simple module  $\dM(\lambda, \mu)$ in $R(\mu-\la)\gmod$,  called  the \emph{determinantial module}.  (See \cite[Section 3.3]{KKOP21} for the precise definition and more properties of them.)

\begin{prop}
[{\cite[Lemma 1.7, Proposition 4.2]{KKOP18}, \cite[Lemma 3.23, Theorem 3.26]{KKOP21}}] \label{Prop: dM properties}
Let $\Lambda \in \wlP_+$, and $\lambda, \mu \in \weyl \Lambda$ with $\lambda \wle \mu$.
\begin{enumerate} [\rm (i)]
\item  $\dM(\lambda, \mu)$ is an affreal simple module.
\item If $ \langle h_i, \lambda \rangle \le 0$ and $s_i \lambda \preceq \mu $, then
$$
\ep_i( \dM(\lambda, \mu)) = - \langle h_i, \lambda \rangle \quad \text{and} \quad  E_i^{(- \langle h_i, \lambda \rangle)} \dM(\lambda, \mu) \simeq \dM(s_i\lambda, \mu).
$$
\item If $ \langle h_i, \mu \rangle \ge 0$ and $ \lambda \preceq s_i \mu $, then
$$
\ep^*_i( \dM(\lambda, \mu)) =  \langle h_i, \mu \rangle \quad \text{and} \quad  E_i^{* (\langle h_i, \mu \rangle)} \dM(\lambda, \mu) \simeq \dM(\lambda, s_i\mu).
$$
\end{enumerate}
\end{prop}

For $w\in\weyl$, we say that
$\la \in \wtl$ is {\em $w$-dominant} if 
$(\beta,\la)\ge0$ for any $\beta\in\prD\cap w^{-1}\nrD$.
In this case, we have $\la-w\la\in\prtl$. 

\Th[{\cite[Theorem 2.18]{KKOP22}}] \label{thm:gdm}
Let $w\in W$ and let $\la\in\wtl$.
Assume that $\la$ is $w$-dominant.
Then there exists a self-dual simple $R(\la-w\la)$-module
$\Mm_w(w\la,\la)$ which satisfies the following conditions.
\bna
\item
If $i\in I$ satisfies $\ang{h_i,w\la}\ge0$,
then $\eps_i\bl\Mm_w(w\la,\la)\br=0$.\label{item eps}

\item
If $i\in I$ satisfies $\ang{h_i,\la}\le0$,
then $\eps^*_i\bl\Mm_w(w\la,\la)\br=0$.\label{item epsstar}
\item 
If $i\in I$ satisfies $s_iw  \prec w$, then
$\Mm_w(w\la,\la)\simeq L(i^m)\hconv \Mm_{s_iw}(s_iw\la,\la)$
where $m=\ang{h_i,s_iw\la} = \eps_i(\Mm_w(w\la,\la)) \in\Z_{\ge0}$.\label{item a}
\item
If $i\in I$ satisfies $ws_i\prec w$, then
$\Mm_w(w\la,\la)\simeq \Mm_{ws_i}(w\la,s_i\la)\hconv L(i^m) $
where $m=\ang{h_i,\la} = \eps^*_i(\Mm_w(w\la,\la))\in\Z_{\ge0}  $.\label{item: sla} \

\item
For any $\mu\in\wlP_+$ such that
$\la+\mu\in\wlP_+$,  we have
$$\Mm(w\mu,\mu)\hconv \Mm_w(w\la,\la)\simeq\Mm\bl w(\la+\mu),\la+\mu\br$$
up to a grading shift.\label{item:wla}
Moreover such an $\Mm_w(w\la,\la)$ is unique up to an isomorphism.
\ee
In particular,
$\Mm_w(w\La,\La)\simeq\Mm(w\La,\La)$ for $\La\in\pwtl$.
\enth

We call $\Mm_w(w\la,\la)$ a {\em generalized determinantial module}.

\Lemma[{\cite[Lemma 2.21]{KKOP22}}] \label{lem:MhconvMw}
Let $w\in W$ and let $\la\in \wtl$
be a $w$-dominant weight.
Then, for any $\mu\in\wlP_+$,  we have
$$\Mm(w\mu,\mu)\hconv \Mm_w(w\la,\la)\simeq\Mm_w\bl w(\la+\mu),\la+\mu\br$$
up to a grading shift.
\enlemma


\subsection{Shuffle lemma}

The shuffle lemma is useful to calculate explicitly 
the actions of $\E_i$ and $\Es$ on a convolution product 
of $R$-modules.
Let $\rt$ be the set of roots, and $\rert$ the set of real roots.
For $\al\in\rert$, set $\sfd_\al\seteq(\al,\al)/2$, $\al^\vee\seteq(\sfd_\al)^{-1}\al$.
We set also $\sfd_i\seteq\sfd_{\al_i}$ for $i\in I$.

\Lemma[{\cite[Proposition 10.1.5]{KKKO18}}]\label{lem:shuffleE}
Let $M\in\Mod(R(\beta))$, $N\in\Mod(R(\gamma))$, $i\in I$ and $n\in\Z_{\ge0}$.
Then $\E_i^{(n)}(M\conv N)$ has an increasing filtration $\st{F_k}_{k\in\Z}$ of
$R(\beta+\gamma-n\al_i)$-submodules such that
\bna
\item $F_k=0$ for $k<0$,
\item $F_k=\E_i^{(n)}(M\conv N)$ for $k\ge n$,
\item $F_k/F_{k-1}\simeq q^{k(n-k)\sfd_i-k(\al_i,\beta)}\E_i^{(n-k)}(M)\conv \E_i^{(k)}(N)$.
\ee
\enlemma
It is a categorification of the formula of the comultiplication:
$$\Delta(e_i^{(n)})=\sum_{k=0}^n q^{k(n-k)\sfd_i}e_i^{(n-k)}t_i^{k}\tens e_i^{(k)}.
$$
\Lemma[{\cite[Proposition 10.1.5]{KKKO18}}]\label{lem:shuffleEs}
Let $M\in\Mod(R(\beta))$, $N\in\Mod(R(\gamma))$, $i\in I$ and $n\in\Z_{\ge0}$.
Then $\Es^{(n)}(M\conv N)$ has an increasing filtration $\st{F_k}_{k\in\Z}$ of
$R(\beta+\gamma-n\al_i)$-submodules such that
\bna
\item $F_k=0$ for $k<0$,
\item $F_k=\Es^{(n)}(M\conv N)$ for $k\ge n$,
\item $F_k/F_{k-1}\simeq q^{k(n-k)\sfd_i-k(\al_i,\gamma)}\Es^{(k)}(M)\conv \Es^{(n-k)}(N)$.
\ee
\enlemma

\Lemma[{\cite[Proposition 10.1.5]{KKKO18}}]
Let $Y,Z\in R\gmod$, and $y,z\in\Z_{\ge0}$.
Then we have
\bnum
\item 
If $\E_i^{y+1}Y\simeq0$, then 
there exists a monomorphism of degree $yz\sfd_i+z\bl\al_i,\wt(Y)\br$:
\eq \E_i^{(y)}Y\conv \E_i^{(z)}Z\monoto \E_i^{(y+z)}(Y\conv Z).\label{mor:YZ1}
\eneq

\item 
If $\E_i^{z+1}Z\simeq 0$, then 
there exists an epimorphism of degree $-yz\sfd_i-z\bl\al_i,\wt(Y)\br$:
\eq \E_i^{(y+z)}(Y\conv Z)\epito \E_i^{(y)}Y\conv E_i^{(z)}Z.
\label{mor:YZ2}\eneq
\item
If $E_i^{y+1}Y\simeq0$ and $E_i^{z+1}Z\simeq 0$, then
\eqref{mor:YZ1} and \eqref{mor:YZ2} are isomorphisms, inverse to each other.
\item 
If $\E_i^{y+1}Y\simeq0$ and $\E_i^{z+1}Z\simeq 0$, 
Then there exists an exact sequence
$$\xymatrix{
0\ar[r]& E_i^{(y)}Y\conv E_i^{(z-1)}Z
\ar[r]^f&E_i^{(y+z-1)}(Y\conv Z)\ar[r]^g&
E_i^{(y-1)}Y\conv E_i^{(z)}Z\ar[r]&0.}$$
Here $f$ is of degree $y(z-1)\sfd_i+(z-1)\bl\al_i,\wt(Y)\br$
and $g$ is of degree $-(y-1)z\sfd_i-z\bl\al_i,\wt(Y)\br$.
\ee
\enlemma
\Proof It follows from Lemma~\ref{lem:shuffleE}. 
\QED

\Lemma[{\cite[Proposition 10.1.5]{KKKO18}}]\label{lem:E*}
Let $Y,Z\in R\gmod$, and $y,z\in\Z_{\ge0}$.
Then we have
\bnum

\item 
If $\Es^{z+1}Z\simeq0$, then 
there exists a monomorphism of degree $yz\sfd_i+y\bl\al_i,\wt(Z)\br$:
\eq \Es^{(y)}Y\conv \Es^{(z)}Z\monoto \Es^{(y+z)}(Y\conv Z).\label{mor:YZs1}
\eneq

\item 
If $\Es^{y+1}Y\simeq 0$, then 
there exists an epimorphism of degree $-yz\sfd_i-y\bl\al_i,\wt(Z)\br$:
\eq \Es^{(y+z)}(Y\conv Z)\epito \Es^{(y)}Y\conv \Es^{(z)}Z.
\label{mor:YZs2}\eneq

\item
If $\Es^{y+1}Y\simeq0$ and $\Es^{z+1}Z\simeq 0$, then
\eqref{mor:YZs1} and \eqref{mor:YZs2} are isomorphisms, inverse to each other.
\item 
If $\Es^{y+1}Y\simeq0$ and $\Es^{z+1}Z\simeq 0$, 
Then there exists an exact sequence
$$\xymatrix{
0\ar[r]& \Es^{(y-1)}Y\conv \Es^{(z)}Z
\ar[r]^f&\Es^{(y+z-1)}(Y\conv Z)\ar[r]^g&
\Es^{(y)}Y\conv \Es^{(z-1)}Z\ar[r]&0.}$$
Here, $f$ is of degree $(y-1)z\sfd_i+(y-1)\bl\al_i,\wt(Z)\br$,
and $g$ is of degree $-y(z-1)\sfd_i-y\bl\al_i,\wt(Z)\br$.
\ee
\enlemma

\Prop\label{prop:E*}
Let $M\in R\gmod$ and $L\in R\gmod$ be simple modules, and assume that one of them is \afr.
Set $L'=\Esm(L)$ and $M'=\Esm(M)$.
Then we have
\bnum
\item $\eps_{i}^*(L\hconv M)\le\eps_{i}^*(L)+\eps_{i}^*(M)$.
\item
$\La(L', M')\le\La(L,M)+\eps_{i}^*(L)\bl\al_{i},\wt(M)\br-\eps_{i}^*(M)\bl\al_{i},\wt(L)\br$
and $$\tLa(L', M')\le\tLa(L,M)+\eps_{i}^*(L)\bl\al_{i},\wt(M)\br
+\eps_{i}^*(L)\eps_{i}^*(M)\sfd_i.$$
If $\eps_{i}^*(L\hconv M)=\eps_{i}^*(L)+\eps_{i}^*(M)$, then the equality holds.
\item
Assume either 
$\eps_{i}^*(L\hconv M)=\eps_{i}^*(L)+\eps_{i}^*(M)$ or
$\La(L', M')=\La(L,M)+\eps_{i}^*(L)\bl\al_{i},\wt(M)\br-\eps_{i}^*(M)\bl\al_{i},\wt(L)\br$.
Then we have $\Esm(L\hconv M)\simeq\ L'\hconv M'$ up to grading shifts,
\ee
\enprop

\Proof
The statement (i) follows immediately from Proposition~\ref{prop:ddd}.

\snoi
(ii) \ 
Let $(\Laa,z)$ be an affinization of $L$ and set $\eps^*_i(L)=n$ and $\eps^*_i(M)=m$.
Then $\Laa'\seteq\Es^{(n)}\Laa$ is an affinization of $L'$.
Let $R\cl \Laa\conv M\to M\conv \Laa$ be
a renormalized R-matrix.
By applying $\Es^{(n+m)}$, we obtain a commutative diagram:
\eqn
\xymatrix@C=10ex{
\Es^{(n)}(\Laa)\conv\Es^{(m)}M\ar[r]^{R'}\ar[d]^\bwr_f&\Es^{(m)}M\conv\Es^{(n)}\Laa
\ar[d]^\bwr_g\\
\Es^{(n+m)}(\Laa\conv M)\ar[r]^{\Es^{(n+m)}(R)}&\Es^{(n+m)}(M\conv \Laa),
}
\eneqn
where the isomorphisms $f$ and $g$ are obtained by Lemma~\ref{lem:E*},  
$\deg(f)=nm\sfd_i+n(\al_i, \wt\,M)$
and $\deg(g)=nm\sfd_i+m(\al_i, \wt\,L)${, which follow from Lemma~\ref{lem:shuffleE}}.
Hence $R'\cl \Laa'\conv M'\to M'\conv \Laa'$ has degree
$\La(L,M)+n(\al_i, \wt\,M)-m(\al_i, \wt\,L)$.
Since $R'=z^s \Rre_{\Laa.M'}$, we have
$\La(L',M')\le\deg(R')$. If $\eps_{i}^*(L\hconv M)=\eps_{i}^*(L)+\eps_{i}^*(M)$, 
then $\Es^{(n+m)}(L\conv M)\to\Es^{(n+m)}(L\hconv M)$ is an isomorphism, and hence
$R'\vert_{z=0}$ does not vanish, which implies that $\La(L',M')=\deg(R')$.
Thus we obtain (ii). (iii) is obvious. 
\QED

\Lemma[{\cite[Lemma 3.13]{KL09}}]
Let $M\in R\gmod$ be a simple module.
Then for any $n\in\Z_{\ge0}$ such that $n\le\eps_i(M)$,
$\E_i^{(n)}(M)$ has a simple socle and a simple head.
Moreover they are isomorphic to $\tE_i^nM$ \ro up to grading shifts\rf.
\enlemma

\section{The categories $\Cw$ and $\tCw[{w}]$}\label{sec:cat-C}

We shall see several basic properties of the category $\catC_{w}$ and its localization $\tCw$, which is 
the main target of this article.
\subsection{Categories $\catC_{w}$} \label{Sec: Cwv}
In this subsection, we recall  the category $\catC_w$
associated with $w\in\weyl$ (see \cite{KKOP18} ). For more details, see
\cite{KKOP21}, \cite{KKOP22}, \cite{KKOP23}.

For $M\in R(\beta) \gMod$ we define
\eqn
&&\ba{l}
\gW(M) \seteq  \{  \gamma \in  \rlQ_+ \cap (\beta - \rlQ_+)  \mid  e(\gamma, \beta-\gamma) M \ne 0  \}, \\
\gW^*(M) \seteq  \{  \gamma \in  \rlQ_+ \cap (\beta - \rlQ_+)  \mid  e(\beta-\gamma, \gamma) M \ne 0  \}.
\ea
\eneqn

For any simple module $M$, we have (\cite{TW16})
\eq
\gW(M)\subset \Sp(\gW(M)\cap\prt).\eneq
Here, for a subset $A$ of $\R\tens\rtl$, we denote by $\Sp(A)$
the smallest convex cone containing $A\cup\st{0}$.

For $w\in\weyl$,  we define the full monoidal subcategory of $R\gmod$
by
\eq
&&
\catC_w\seteq\st{M\in R\gmod\Mid \W\,(M)\subset\prtl\cap\, w\nrtl} .
\label{def:Cw}
\eneq

An ordered pair $(M,N)$ of $R$-modules is called {\em unmixed}
if $$\sgW(M)\cap\gW(N)\subset\{0\}.$$

\Lemma[{\cite[Corollary 2.3]{KKOP22}}]\label{lem:unm}
If $(M,N)$ is unmixed, then $(M,N)$ is $\La$-definable and
$\tLa(M,N)=0$,
\enlemma
For $\La \in \pwtl$ and $w, v\in \weyl$ such that $v\le w$, we have
\eq
\dM(w\La,v\La) \in \catC_w.
\label{det-in-C}
\eneq

\Rem
Even if $\la$ is $w$-dominant, $\Mm_w(w\la,\la)$ may not belong to 
$\Cw$.
For example, take $\g=A_2$, $w=s_1s_2$ and $
\la=s_1\La_1$.
Then $\Mm_w(w\la,\la)\simeq \ang{2}$ does not belong to $\Cw$,
since
$\al_2\not\in\prD\cap w\nrD$.
\enrem

\subsection{Localizations of $\catC_{w}$} 
\label{subsec:locCw}
In this subsection we briefly recall  the localizations of the categories $\catC_{w}$ via left braiders studied in \cite{KKOP21, KKOP22}.

Let $w \in \weyl$.
Then there exists a real commuting  family of graded left braiders in $R\gmod$
\eq
\st{\bl\dM(w\La_i,\La_i),\coRl_{\dM(w\La_i,\La_i)}, \dphi_{\dM(w\La_i,\La_i)}\br}_{i\in I} \label{eq:leftbraider}
\eneq
(see \cite[Proposition 5.1]{KKOP21} for notations).
Here $\coRl_{\dM(w\La_i,\La_i)}(X)\cl \dM(w\La_i,\La_i)\conv X\to
X\conv\dM(w\La_i,\La_i)$ is a morphism of degree
$\dphi_{\dM(w\La_i,\La_i)}(\wt(X))$ functorially in $X\in R\gmod$.
Moreover it is a family of central objects in the category $\catC_w$,
i.e. $\coRl_{\dM(w\La_i,\La_i)}(X)$ is an isomorphism for any $X\in\Cw$. 
Note that
\eq \label{eq:phi_w}
\dphi_{\dM(w\La_i,\La_i)}(\beta)= -(w\La_i+\La_i, \beta) \qt{for any } \ \beta\in \rtl. 
\eneq

We denote by $\bl R\gmod\br[\dM(w\La_i,\La_i)^{\circ -1}; i\in I]$
the localization of $R\gmod$ via the real commuting  family of graded left braiders \eqref{eq:leftbraider},
and by $\tCw$ the localization of $\Cw$ by the same left braiders.
We have a quasi-commutative diagram
\eqn
\xymatrix@C=10em{
\catC_w \akew\ar@{>->}[r]_{\iota_w} \ar[d]_{\Phi_w}& R\gmod \ar[d]_{\Qt} \\
\tcatC_w \akew\ar@{>-->}[r]_{\widetilde\iota_w}&\bl R\gmod \br[\dM(w\La_i,\La_i)^{\circ -1}; i\in I]
}
\eneqn
where $\Phi_w$ and $\Qt$ denote the localization functors,   and $\widetilde \iota_w$ is the induced functor from  the inclusion functor $\iota_w$. 

\begin{thm} [{\cite[Theorem 5.9, Theorem 5.11]{KKOP21}, \cite[Theorem 3.9]{KKOP22}}] \label{thm:rigidity-Cw}
  \hfill
  
\bnum
\item The functor $\widetilde \iota_w \cl  \tcatC_w \to \bl R\gmod \br[\dM(w\La_i,\La_i)^{\circ -1}; i\in I]$ is an equivalence of categories.
\item For any simple $M\in R\gmod$,
  $\Qt(M)\in\tCw$ is simple or zero.
\item The monoidal category $\tcatC_w$ is rigid, that is, every object of $\tcatC_w$ has a left dual and a right dual. Moreover $\tCw$ satisfies condition
  \eqref{cond:fcat}.
\ee
\end{thm}

\Th[{\cite[Theorem 8.3, Corollary 8.4]{ref} }]\label{th:ff}
\hfill
 \bnum
 \item
$\Phi_w\cl \Cw\to\tCw$ is fully faithful.
\item
The full subcategory $\catC_w$ of $\tcatC_w$ is stable by taking subquotients.
\ee
\enth

For $w\in\weyl$, let us take a reduced expression $\uw=s_{i_1}\cdots s_{i_\ell}$ of $w$.
We set
$$I_w\seteq\st{i\in I\mid w\La_i\not=\La_i}=\st{i_1,\ldots,i_\ell}.$$

Then
the localization functor $\Qt$.
factors as
$$R\gmod\To R_{I_w}\gmod\to\Cw.$$
Here, $R_{I_w}$ denotes the quiver Hecke algebra defined by using $I_w$ instead of $I$,
and $R\gmod\To R_{I_w}\gmod$ is the exact monoidal functor given by
$$\text{$R(\al_i)\mapsto R_{I_w}(\al_i)$ or $R(\al_i)\mapsto 0$ 
whether $i\in I_w$ or not.}$$

\Rem
{\em Although the composition $\Cw\to R\gmod\to\tCw$ is fully faithful
{\rm (\cite[Theorem 8.3]{ref})},
we do not regard an object of\/ $\Cw$ as an object of\/ $\tCw$.}
For $M\in\Cw$, we distinguish $M\in\Cw$ and $\Qt(M)\in\tCw$.
\enrem
For $\La\in\pwtl$, we set
\eqn
\dC_\La\seteq\dM(w\La_,\La)\in\Cw,\qtq
\tdC_\La\seteq\Qt(\dC_\La)\in\tCw.
\eneqn

For $\la\in\wtl$, define
\eq\tdC_\la\seteq(\tdC_{\mu})^{-1}\conv\tdC_{\la+\mu}\in\tCw
\label{eq:center}
\eneq
with $\mu\in\pwtl$ such that $\la+\mu\in\pwtl$ (up to grading shifts).
Hence, we have
$\tdC_\la\conv\tdC_\mu\simeq\tdC_{\la+\mu}$ for any $\la,\mu\in\wtl$.
If we want to emphasize $w$, we write
$\tdC^w_\la$ for $\tdC_\la$.

Note that 
\eq\Qt\bl\dM_w(w\la,\la)\br\simeq\tdC_\la\qt{for any $w$-dominant $\la\in\wtl$.}
\label{eq:gendet}
\eneq
Indeed, by Lemma~\ref{lem:MhconvMw} we have
$$\dM(w\mu,\mu)\hconv\dM_w(w\la,\la)\simeq\dM(w(\la+\mu),\la+\mu)$$
for any $\mu\in\pwtl$ such that $\la+\mu\in\pwtl$.

\Prop\label{prop:aff-simple}
Let $L\in R\gmod$ be a simple module,
and let $(\Laa,z)$ be an affinization of $L$ in the sense of
{\rm \/Definition \ref{def:affinization1}}.
Then $(\Qt(\Laa),z)$ is an affinization of\/ $\Qt(L)$
in the sense of {\rm \/ Definition \ref{def:affinization}}.
\enprop

\Proof
Let $\sha$ be the abelian monoidal category
consisting of pairs $(X,R_X)$
where $X\in\tCw$ and $R_X\cl \Qt(\Laa)\conv X\isoto X\conv\Qt(\Laa)$ is an isomorphism
in $\Rat(\tCw)$ (see Definition~\ref{Def: Raff}).
The morphisms are obvious ones.
The monoidal category structure of $\sha$ is given by:
$(X,R_X)\tens (Y,R_Y)=(X\conv Y,R_{X\conv Y})$ where
$R_{X\conv Y}$ is the composition
$\Qt(\Laa)\conv X\conv Y\To[R_X\circ Y]
X\conv\Qt(\Laa)\conv Y \To[X\circ R_Y] X\conv Y\conv\Qt(\Laa)$.

Then we can define
the monoidal functor $\Psi\cl\Cw\to \sha$ by
$M\mapsto \bl\Qt(M),\Runi_{\Laa,M}\br$.
In order to see that $\Psi$ extends to $\widetilde\Psi\cl \tCw\to\sha$,
it is enough to show that $\Psi(\dC_\La)$ is invertible for any $\La\in\pwtl$.
Set $\ttC=\Qt(\dC_\La)$ and $\tLaa=\Qt(\Laa)$.
Then, we have an isomorphism
$\tLaa\conv\ttC\isoto\ttC\conv\tLaa$ in $\Rat(\tCw)$.
By applying $\ttC^{-1}\conv\;\scbul\;\conv\ttC^{-1}$, we obtain
an isomorphism
$\ttC^{-1}\conv\tLaa\isoto\tLaa\conv\ttC^{-1}$.
Let $R_{\ttC^{-1}}$ be its inverse.
Then we can see easily that $(\ttC^{-1},R_{\ttC^{-1}})$ is an inverse object
of $\Psi(\dC_\La)$.

Thus, we see that  $\Psi$ extends to $\widetilde\Psi\cl \tCw\to\sha$.
For $X\in\tCw$, set $(X,R_X)=\widetilde\Psi(X)$, and
let us define $R_{\Qt(\Laa)}(X)$ by  $R_X\cl\Qt(\Laa)\tens X\to X\conv\Qt(\Laa)$.
Then $\bl\Qt(\Laa),R_{\Qt(\Laa)}\br$ is a rational center
(see Definition~\ref{def:ratcent}).
Thus $\bl\Qt(\Laa),R_{\Qt(\Laa)}\br$ is an affinization of $\Qt(L)$
in the sense of
Definition \ref{def:affinization}.
\QED

\Prop\label{prop:shsw}
Let $M$ and $N$ be simple modules in $\Cw$.
Assume that $\bl\Qt(M),\Qt(N)\br$ is a $\La$-definable pair of simple objects in $\tCw$.
Then $(M,N)$ is a $\La$-definable pair.
In particular, $M\conv N$ has a simple head, $N\conv M$ has a simple socle,
and we have
$$M\hconv N\simeq\Im(\rmat{M,N})\simeq N\sconv M.$$
\enprop
\Proof
It is an immediate consequence of
 Theorem~\ref{th:ff} and  Lemma~\ref{lem:Laspl}.
\QED

\subsection{Equivalence between $(\tCw)^{\rev}$ and $\tCw[{w^{-1}}]$}

Recall that, for a monoidal category $(\shc,\tens)$, we denote by $\shc^\rev$
the monoidal category $\shc$ with the tensor product $\tens_\rev$: $M\tens_\rev N\seteq N\tens M$.

Let $\psi\cl R(\beta)\isoto R(\beta)$ be the ring automorphism
\eq \label{eq:psi}
\psi&\;:\hs{2ex}&\ba{l}
e(\nu_1,\ldots ,\nu_n) \mapsto e(\nu_n,\ldots ,\nu_1), \\
 x_k \mapsto x_{n+1-k} \qt{($1\le k\le n$),}\\
\tau_l\mapsto -  \tau_{n-l}\qt{($1\le l<n$),}
\ea\label{def:antipsi}
\eneq
 where $n=\height{\beta}$.

Then $\psi$ induces an equivalence of monoidal categories (\cite[\S\,3.2.]{KKOP22})
\begin{equation}
\psi_*\cl (R\gmod)^\rev\simeq R\gmod.
\label{rev-iso}
\end{equation}

\Lemma[{\cite[Lemma2.23]{KKOP22}}]\label{lem:*-inv}
We have
$$\psi_*\bl\Mm_w(w\la,\la)\br
\simeq\Mm_{w^{-1}}(-\la,-w\la).$$
for any $w\in W$ and any $w$-dominant $\la\in\wtl$
\enlemma

\begin{thm}[{\cite[Theorem 3.7]{KKOP22}}] \label{th:main1}
There is an equivalence of  monoidal categories
$(\tCw)^\rev$ and $\tCw[{w^{-1}}]$.
More precisely, we have a quasi-commutative diagram 
$$\xymatrix@C=10ex@R=5ex{
(R\gmod)^\rev\ar[r]^-\ssim_-{\psi_*}\ar[d]_\Qt& R\gmod\ar[d]_{\Qt[w^{-1}]}\\
(\tCw)^\rev\ar[r]^-\ssim_-{\psi_*}&\tCw[w^{-1}].}
$$
\end{thm}

\subsection{The category $\CBw$}\label{loc-Dema}

In this subsection, we shall see the basic properties of the localization functor $\Qt$ and the category 
$\CBw$.

Recall that $S(b)$ denotes
the self-dual simple $R$-module corresponding to $b\in B(\infty)$.

In Proposition \ref{Demazure}, we introduced
a subset $B_{w}(\infty)$ of $B(\infty)$ for each $w\in\weyl$.

\begin{prop} [{\cite[Proposition 3.1.]{KKOP22}}]\label{prop:binBw}
  Let $w\in\weyl$ and let $\underline{w} = s_{i_1}s_{i_2}\cdots s_{i_\ell}$ be a reduced expression of $w$.
Let $M$ be a simple $R$-module.
Then the following conditions are equivalent.
\bna
\item 
$e(i_1^{n_1},\ldots,i_l^{n_l})M\not\simeq0$ for some
$(n_1,\ldots, n_\ell)\in\Z_{\ge0}^\ell$,
\item $\tE_{i_\ell}^{\max} \tE_{i_{\ell-1}}^{\max} \cdots \tE_{i_2}^{\max} \tE_{i_1}^{\max} M \simeq \one$,
\item $\tEsm_{i_1}\cdots\tEsm_{i_\ell}M\simeq\one$,
  \item $M$ is a quotient of 
$\ang{i_1^{n_1}}\conv\cdots\conv\ang{i_\ell^{n_\ell}}$ for some
$(n_1,\ldots, n_\ell)\in\Z_{\ge0}^\ell$,
\item $M\simeq S(b)$ up to a grading shift for some  $b \in B_w(\infty)$,
\item $\Qt(M)$ is a simple module,
\item $\Qt(M)\not\simeq0$,
\item $\La(\dC^w_\La, M)=-\bl w\La+\La,\wt(M)\br$ for any
  $\La\in\pwtl$,
\item\label{item:i}
$\dC^w_\La\hconv M\in\Cw$ for some $\La\in\pwtl$.
\end{enumerate}
\end{prop}

Let $\CBw$ be the full subcategory of $R\gmod$ consisting of modules $M$
such that any simple subquotient $S$ of $M$ satisfies
the equivalent conditions above.

Let $\Irr(\Cw)$ 
be the set of equivalence classes of simple
modules in $\Cw$ up to grading shifts.
Namely, for simple modules $M$ and $N$ in $\Cw$,
we identify $M$ and $N$ in $\Irr(\Cw)$ if $M\simeq q^nN$ for some $n\in\Z$.
Similarly, let $\Irr(\CBw)$ 
be the set of equivalence classes of simple
modules in $\CBw$ up to grading shifts.
Hence, we have
$\Irr(\Cw)\subset \Irr(\CBw)$.

We have the following isomorphism of crystals 
in the sense of Definition~\ref{morph} (ii)
\[
\Irr(\CBw)\simeq B_w(\infty).
\]
\Lemma 
$\Irr(\CBw)\sqcup\st{0}$ is stable by $\tE_i$ and $\tE^*_i$.
\enlemma

\Proof
Since the case of $\tE^*_i$ is shown similarly to the case of $\tE_i$, we shall show the stability only for $\tE_i$. 
Take  $S\in\Irr(\CBw)$.
The assertion is trivial when $\tE_iS\simeq0$, we may assume that
$\tE_i S$ is simple. 
By the condition (d) above, there exists an epimorphism
$\ang{i_1^{n_1}}\conv\cdots\conv\ang{i_\ell^{n_\ell}}\epito S$ for some
$(n_1,\ldots, n_\ell)\in\Z_{\ge0}^\ell$.
Since $\E_i$ is an exact functor, we have an epimorphism
\[
\pi\cl \E_i(\ang{i_1^{n_1}}\conv\cdots\conv\ang{i_\ell^{n_\ell}})\epito \E_i S\epito \tE_i S.
\]
By applying Lemma~\ref{lem:shuffleE} (i),
there exists a filtration
$$0=F_{0}\subset F_1\subset \cdots\subset F_{\ell-1}\subset F_\ell
=\E_i(\ang{i_1^{n_1}}\conv\cdots\conv\ang{i_\ell^{n_\ell}})$$
such that $F_k/F_{k-1}$ is isomorphic to 
$\ang{i_1^{n_1}}\conv\cdots\conv\E_i(\ang{i_k^{n_k}})\conv\cdots\conv\ang{i_\ell^{n_\ell}}$.
Let $k_0\in\{1,2,\cdots,\ell\}$ be the minimum among $k's$ such that $\pi_{|F_k}\cl F_k\epito \tE_iS$ does not vanish.
Then, $\pi\vert_{F_{k_0}}\cl F_{k_0}\to \tE_iS$ is an epimorphism
since $\tE_iS$ is simple.
Hence we obtain $\E_i(\ang{i_{k_0}^{n_{k_0}}})\not\simeq0$, which implies
$i=i_{k_0}$, $n_{k_0}>0$, and obtain
\[
\til\pi\cl F_{k_0}/F_{k_0-1}\simeq
\ang{i_1^{n_1}}\conv\cdots\conv\ang{i_{k_0}^{{n_{k_0}}-1}}\conv\cdots\conv\ang{i_\ell^{n_\ell}}\epito\tE_iS.
\]
\QED

\Lemma\label{lm:stable}
$\Irr(\Cw)\sqcup\st{0}$ is stable by $\tE^*_i$.
\enlemma
\Proof
This follows from
$
\gW(\tE^*_i S)\subset \gW(S)$.
\QED

\Rem\label{rm:non-stable}
$\Irr(\Cw)\sqcup\st{0}$ is not necessarily stable by $\tE_i$. 
For example, for $\g=A_2$, $w=s_1s_2$, we have $\ang{12}\in\Cw$ but
$\tE_1\ang{12}\simeq\ang{2}\not\in\Cw$. However $\ang{2}\in\CBw$.
\enrem

The following lemma immediately follows from the definitions and Theorem~\ref{thm:gdm}.
\Lemma
For any $w$-dominant $\la\in\wtl$, we have
$$\dM_w(w\la,\la)\in\CBw.$$
\enlemma

\Lemma\label{lem:CBww'}
Let $w\in\weyl$ and $i\in I$.
Assume that $w'\seteq ws_i<w$.
\bnum\item
If a simple $S\in\CBw$ satisfies $\eps^*_i(S)=0$, then
$S\in\CBw[w']$.
\item If $S\in\CBw$ is simple, then $({{\tFs_i}})^m(S)\simeq S\hconv\ang{i^m}\in\CBw$
  for any $m\in\Z_{\ge0}$
  \ee
\enlemma
\Proof
Take a reduced expression $s_{i_1}\cdots s_{i_\ell}$ of $w$ with $i=i_\ell$.
By the assumption there exists an epimorphism
$$\ang{i_1^{m_1}}\conv\cdots\conv \ang{i_\ell^{m_\ell}}\epito S.$$

\snoi
(i)\ Since $\eps^*_i(S)=0$, we have $m_\ell=0$, which implies that $S\in\CBw[w']$.

\snoi
(ii) We have an epimorphism
$$\ang{i_1^{m_1}}\conv\cdots\conv \ang{i_{\ell-1}^{m_{\ell-1}}}
\conv \ang{i_\ell^{m+m_\ell}}\epito S\hconv\ang{i^m}.$$
\QED

\Lemma\label{lem:lasimp}
Let $M, N\in R\gmod$ be simple modules.
If there exits a simple subquotient $S$ of $M\conv N$ such that
$S\in\CBw$, then $M$ and $N$ belong to $\CBw$.
\enlemma
\Proof
Since $\Qt(S)\not\simeq0$ is a subquotient of
$\Qt(M\conv N)\simeq\Qt(M)\conv\Qt(N)$,
we have
$\Qt(M)\not\simeq0$ and $\Qt(N)\not\simeq0$.
\QED

\Lemma\label{lem:MNla}
Let $M$ be a simple module in $\Cw$
and let $N$ be a simple module in $\CBw$.
Then any simple quotient of $M\conv N$ belongs to $\CBw$.
\enlemma
\Proof
Let $S$ be a simple quotient of $M\conv N$.
We shall show that $S$ belongs to $\CBw$.

By Lemma~\ref{lem:Laspl}, $S$ is a simple submodule of $N\conv M$. 
Take $\La\in\pwtl$ such that $\dC_\La\hconv N$ belongs to $\Cw$ (see Proposition~\ref{prop:MCw}).
Then the composition
$$\dC_\La\conv S\monoto\dC_\La\conv N\conv M\epito (\dC_\La\hconv N)\conv M$$
does not vanish by Lemma~\ref{lem:monoepi}. Since $(\dC_\La\hconv N)\conv M\in\Cw$, we conclude that
$\dC_\La\conv S$ has a simple quotient in $\Cw$,  
which implies the desired result
by Proposition~\ref{prop:binBw}.
\QED

\Lemma\label{lem:wLadf}
Let $M$ and $N$ be simple modules in $R\gmod$, and
let $f\cl M\conv N\to N\conv M$ be a morphism.
Assume that
\bna
\item
$S\seteq\Im(f)$
is simple and belongs to $\CBw$,
\item
$\bl \Qt(M),\Qt(N)\br$ is a $\La$-definable pair \ro see {\rm Definition~\ref{def:rmat}}\rf.
\ee
Then we have
\bnum
\item $M,N\in\CBw$,
\item
  $\La\bl\Qt(M),\Qt(N)\br=\deg(f)$.
  \ee
\enlemma
\Proof
We have a diagram in $\tCw$:
$$\Qt(f)\cl \Qt(M)\conv\Qt(N)\epito\Qt(S)\monoto\Qt(N)\conv\Qt(M).$$
Since $\Qt(S)$ is simple, $\Qt(M)$ and $\Qt(N)$ do not vanish, and hence
we obtain (i).
Since $\Qt(f)$ does not vanish, it coincides with the R-matrix
$\rmat{\Qt(M),\Qt(N)}$, and hence we obtain (ii).
\QED

The following lemma is an immediate consequence of the lemma above.

\Lemma\label{lem:LMN}
Let $M,N$ be simple modules in $R\gmod$. Assume that $(M,N)$ is $\La$-definable.
If $M\hconv N\in\CBw$, then we have
\bnum
\item $M,N\in\CBw$,
\item$\Qt(M)\hconv\Qt(N)\simeq \Qt(M\hconv N)$,
\item
$\La(M,N)=\La\bl\Qt(M),\Qt(N)\br$.
\ee
\enlemma

\Prop\label{prop:CB}
Let $L\in\Cw$ and $M\in\CBw$ be simple modules.
Assume that  $(L,M)$ is $\La$-definable.
Then we have
\bnum
\item $L\hconv M\in\CBw$,
\item $\Qt\bl L\hconv M\br \simeq \Qt(L)\hconv \Qt(M)$,
\item $\La(L,M)=\La\bl \Qt(L), \Qt(M)\br$.
\ee
\enprop
\Proof
(i) follows from Lemma~\ref{lem:MNla}.
The remaining statements follow from Lemma~\ref{lem:LMN}.
\QED

\Lemma\label{lem:LaQt}
Let $M$ and $N$ be simple modules in $\CBw$.
Assume that one of them is \afr. Then, we have
$$\La\bl\Qt(M),\Qt(N)\br\le\La(M,N).$$
\enlemma
\Proof
The proof is similar to the one of Lemma~\ref{lem:LMN}.
Assume that $M$ is \afr.
Let $(\Ma,z)$ be an affinization of $M$, and let
$R\cl \Ma\conv N\to N\conv \Ma$
be the renormalized R-matrix.
Then $\La(M,N)$ is the homogeneous degree of $R$.
Applying $\Qt$ we obtain a morphism in $\Aff(\tCw)$:
$$\Qt(R)\cl \Qt(\Ma)\conv \Qt(N)\to \Qt(N)\conv \Qt(\Ma).$$
Take the largest $c\in\Z_{\ge0}$ such that
the image of $\Qt(R)$ is contained in $\Qt(N)\conv z^c\Qt(\Ma)$.
Then $\La\bl\Qt(M),\Qt(N)\br$ is the homogeneous degree of $z^{-c}\Qt(R)$.
Hence we obtain the desired result.

We can similarly prove the assertion when $N$ is \afr.
\QED

\Rem\label{rem:12}
For \afr $L\in\Cw$ and $M\in\CBw$, the simple module $M\hconv L$ may not belong to $\CBw$.
For example, take $\g=A_2$, $w=s_1s_2$,
$L=\ang{1}$, $M=\ang{2}$.
Then $L\hconv M\simeq\ang{12}\in\Cw$ and
$M\hconv L\simeq\ang{21}\not\in\CBw$.
\enrem

\smallskip
For a $\La$-definable pair  $(X,Y)$ of simple objects of $\tCw$,
 we set
\[
\tLa(X,Y)=\Bigl(\La(X,Y)+\bl\wt(X),\wt(Y)\br\Bigr)/2.
\]
If one of $X$ and $Y$ is \afr, then $\tLa(X,Y)\in\Z$ (cf.\ \cite[Lemma~3.11]{KKOP21}). 

The following lemma is obtained from \cite[5.1]{KKOP21}.
\Lemma\label{lem:Lacenter}
For any $\la\in\wtl$ and any simple $X\in\tCw$, we have
\eq
&&\ba{l}
\La(\tdC_\la, X)=-\bl w\la+\la,\wt(X)\br\qtq\La(X,\tdC_\la)=\bl\wt(X),w\la+\la\br,\\[.5ex]
\tLa(\tdC_\la, X)=-\bl\la,\wt(X)\br\qtq\tLa(X,\tdC_\la)=\bl\wt(X),w\la\br.
\ea
\eneq
\enlemma
Note that they are possibly negative integers.
  \Proof
  Since $\La(\tdC_{\la-\mu}, X)=\La(\tdC_\la, X)-\La(\tdC_\mu, X)$, etc.,
  we may assume that $\la\in\pwtl$.
For any $X\in\tCw$, we can find 
$\mu\in\pwtl$ and a simple $M\in\Cw$ such that
$X\simeq \tdC_{-\mu}\conv\Qt(M)$. 
We also get $\La(\dC_\la,M)=-(w\la+\la,\wt(M))$ by \cite[Proposition~5.1]{KKOP21}.
Thus, by \cite[Lemma 3.16 (ii),\,(5.1)]{KKOP21} and Proposition~\ref{prop:CB} (iii), we obtain
\eqn
\La(\tdC_\la, X)&=&-\La(\tdC_\la, \tdC_\mu)+\La(\tdC_\la, \Qt(M))
=-\La(\dC_\la,\dC_\mu)+\La(\dC_\la, M)\\
&=&-(w\la+\la,\mu-w\mu)-
(w\la+\la,\wt(M))
=-(w\la+\la,\wt(X)),
\eneqn
where also note that $\La(\dC_\la,\dC_\mu)=(\la+w\la,\mu-w\mu)$ as in 
\cite[(5.1)]{KKOP21}.
\QED

\Cor \label{cor:epseps}
For any $i\in I_w$ and any simple $M\in\Cw$, we have
\[
\tLa\bl \Qt(M),\Qt(\ang{i})\br=\tLa(M,\ang{i})=\sfd_i\eps_i^*(M).
\]
\encor
\Proof
By Lemma \ref{LmLmtilde},
$\tLa(M,\lan i\ran)=\sfd_i\vep^*_i(M)$, and by Example~\ref{ex:i-affreal}, 
Proposition~\ref{prop:CB}
implies that
\[
  \tLa\bl \Qt(M),\Qt(\ang{i})\br=\tLa(M,\ang{i}).\qedhere
\]
\QED

\Prop\label{prop:eit-fit}
Let $i\in I_w$ and let $M\in\CBw$ be a simple module. 
\bnum
\item
Assume that $\tF_i(M)\in\CBw$.
Then we have
$$\Qt(\tF_i(M))\simeq\Qt(\ang{i})\hconv\Qt(M).$$
\item Assume that $\tE_i(M)\not\simeq0$.
Then we have
$$\Qt(\tE_i(M))\simeq\Qt(M)\hconv\D\bl\Qti\br.$$
\ee
\enprop
\Proof

\snoi
(i) follows from Lemma~\ref{lem:LMN}.

\snoi
(ii)\ 
{Since $M\simeq \tF_i\tE_i M = \ang{i}\hconv\tE_i M$, }
it follows from (i) that 
$\Qt(M)\simeq\Qti\hconv\Qt(\tE_i( M))$, which implies that
$\Qt(\tE_i(M))\simeq\Qt(M)\hconv\D\bl\Qti\br$ by
Lemma~\ref{lem:invhvonv}.
\QED
Note that this proposition motivates us to define the crystal operators on $\Irr(\tCw)$ in \S\;\ref{sec:crystal}.

\subsection{Convolutions with $\dC_\La$}
In this subsection, we shall prove Proposition~\ref{prop:MCw},
which is a more precise statement than Proposition~\ref{prop:binBw}\;\eqref{item:i}.

\Lemma\label{lem:Mvan}
Let $M$ be a simple module in $R\gmod$ such that $\Qt(M)$ is simple and 
$\gW^*(M)\in\prtl\cap w\prtl$.
Assume that $M$ commutes with $\dC_{\La_i}$ for any $i\in I_w$ such that $ws_i>w$.
Then $M\simeq\one$.
\enlemma

\Proof
We may assume that $I_w=I$.
Set $K=\st{i\in I\mid ws_i>w}$.
Let $\beta=-\wt(M)\in\prtl$.
For any $\La\in\pwtl$, we have
$\La(\dC_\La,M)=(w\La+\La,\beta)$.
On the other hand, since $\dC_\La\in\Cw$, the pair $(M,\dC_\La)$ is unmixed, and hence Lemma~\ref{lem:unm} implies that
$$\La(M,\dC_\La)=-(w\La-\La,-\beta)=(w\La-\La,\beta).$$
Therefore, we obtain
$$\de(\dC_\La,M)=(w\La,\beta).$$
Hence we obtain
$$(w\La_i,\beta)=0\qt{for any $i\in K$,}$$
which means that
$w^{-1}\beta\in \sum_{i\in I\setminus K}\Z\al_i$.
Since $\beta\in w\prtl$, we have
$\beta\in \sum_{i\in I\setminus K}\Z_{\ge0}w\al_i$.
Since $w\al_i\in\nrtl$ for any $i\in I\setminus K$, we obtain
$\beta\in\nrtl$. Since $\beta\in \prtl$, we conclude that $\beta=0$.
\QED

\Prop\label{prop:ML}
Let $M\in R\gmod$ be a simple module $R\gmod$ such that $\Qt(M)$ is simple.
If $M$ commutes with $\dC_{\La_i}$ for any $i\in I_w$ such that $ws_i>w$,
then $M\in\Cw$.
\enprop
\Proof
We may assume that $I_w=I$.
We can write (\cite{TW16})
$$M\simeq X\hconv Y$$
with a simple $X$ such that $\gW^*(X)\subset\prtl\cap w\prtl$
and a simple $Y\in\Cw$.
Since $\Qt(M)$ is simple, the objects
$\Qt(X)$ and $\Qt(Y)$ are also simple.
For any $\La\in\pwtl$, we have (see Proposition~\ref{prop:binBw})
\eqn
\La(\dC_\La, M)=
-\bl w\La+\La,\wt(M)\br&&=-\bl w\La+\La,\wt(X)\br-\bl w\La+\La,\wt(Y)\br\\
&&=\La(\dC_\La, X)+\La(\dC_\La, Y).\eneqn
On the other hand, since $(X,\dC_\La)$ is unmixed, $(X,Y,\dC_\La)$ is a normal sequence by Proposition~\ref{prop:Normal} and Lemma~\ref{lem:unm}, and hence we have
$$\La(M,\dC_\La)=\La(X,\dC_\La)+\La(Y,\dC_\La).$$
Hence we have
$$\de(\dC_\La,M)=\de(\dC_\La,X)+\de(\dC_\La,Y)=\de(\dC_\La,X).$$
Hence $X$ commutes with $\dC_{\La_i}$ for any $i\in I_w$ such that $ws_i>w$.
Then $X\simeq\one$ by Lemma~\ref{lem:Mvan},
which implies that $M\in\Cw$.
\QED
We need the following proposition in \S\,\ref{sec:proof} to prove the main theorem.
\Prop\label{prop:MCw}
Let $M$ be a simple module in $R\gmod$ such that $\Qt(M)$ is simple.
Set $K=\st{i\in I_w\mid ws_i>w}$.
Then there exists
$\La\in\sum_{i\in K}\Z_{\ge0}\La_i$ such that
$\dC_\La\hconv M\in\Cw$.
\enprop
\Proof
Set $\la=\sum_{i\in K}\La_i$.
Proposition~\ref{prop:simplehd} \eqref{it:de tens} 
implies that $\de(\dC_{\la},(\dC_{\la})^{\conv n}\htens M)=0$ for some $n>0$. 
Hence
setting $\La=n\la\in\sum_{i\in K}\Z_{\ge0}\La_i$,
the module $\dC_\La\hconv M$ commutes with $\dC_{\La_i}$
for any $i\in K$.
Since $\Qt(\dC_\La\hconv M)$ is simple, Proposition~\ref{prop:ML} implies that $\dC_\La\hconv M\in\Cw$.
\QED

\section{Crystal structure on $\tCw$}\label{sec:crystal}
Let $\Irr(\tCw)$ be the
set of the equivalence classes of simple
objects in $\tCw$ up to grading shifts. 
In this section, we introduce a crystal structure on $\Irr(\tCw)$, which is 
one of our main aims.

\subsection{Root objects in the localization $\tCw$}
Recall that an object $X$ in a monoidal category $\shc$ is {\em invertible}
if there exists $X'\in\shc$ such that $X\tens X'\simeq X'\tens X\simeq\one$.
If $\shc$ is rigid, then $X'\simeq\D X\simeq \D^{-1}X$.
\Lemma\label{lem:inv} Let $L\in\tCw$ be a simple object
and assume that $L$ is invertible.
Then any simple $X\in\tCw$ strongly commutes with $L$.
\enlemma
\Proof
We have
$X\simeq X\conv(L\conv\D L)
\simeq(X\conv L)\conv \D L$.
Hence, $(X\conv L)\conv \D L$ is simple, which implies 
that $X\conv L$ is simple.
\QED

\Lemma\label{lem:inv-d=0}
An \afr simple $L\in\tCw$ is invertible if and only if $\de(L,\D L)=0$.
\enlemma
\Proof
It is obvious that the invertibility of $L$ implies $\de(L,\D L)=0$.
Assume that $\de(L,\D L)=0$. 
Then $L\conv\D L\simeq L\conv\D L$ is simple,
and hence 
$L\conv \D L\to 1$ is an isomorphism, which implies that $L$ is invertible.
\QED

\Def\label{def:rootob}
Let $L\in\tCw$ be a real simple object.
We say that $L$ is a {\em root object} \ro of degree $d_L\in\Z_{>0}$\rf if
\bna
\item $L$ has an affinization $(\Laa,z)$ with $\deg(z)=2d_L$ \ro see
  Definition~\ref{def:affinization}\rf,
\item $\de(L,\D^{-1} L)=d_L$.\label{it:dual}
\ee
\edf

\Lemma If $L\in\tCw$ is a root object, then
$\D^{\pm1}L$ is also a root object.
\enlemma
\Proof
By \cite[Proposition 5.6]{ref}, $\D^{\pm1}L$ has an affinization of degree $2d_L$.
\QED
\Lemma
If $L\in\tCw$ is a root object, then
$\La(L,\D L)=0$ and $\La(\D L,L)=2d_L$.
\enlemma
\Proof
We have
$2d_L=2\de(L,\D L)=\La(L,\D L)+\La(\D L,L)$
and $\La(L,\D L)=\La(L,L)=0$ by Lemma~\ref{lem:MNDM}. Hence  we obtain
$\La(\D L,L)=2d_L$.
\QED
\Prop[{cf.\ Proposition \ref{prop:simplehd}}]
\label{prop:root}
Let $L$ be a root object in $\tCw$ of degree $d_L$, and $M$ a simple object in $\tCw$.
If $\de(L,M)>0$, then we have
\bnum
\item
\eqn
&&\La(L,L\hconv M)=\La(L,M)\qtq\La(L\hconv M,L)=\La(M,L)-2d_L,\\
&&\La(L,M\hconv L)=\La(L,M)-2d_L\qtq\La(M\hconv L,L)=\La(M,L),\\
&&\de(L\hconv M,L)=\de(M\hconv L,L)=\de(M,L)-\sfd_L.
\eneqn
\item
Assume further that $\bl\wt(L),\wt(L)\br=2d_L$.
Then, we have 
$$\tLa(L,M\hconv L)=\tLa(L,M)\qtq\tLa(L\hconv M,L)=\tLa(M,L).$$
\ee
\enprop
\Proof

Since $\La(L\hconv M,L)<\La(M,L)$ by Proposition~\ref{prop:simplehd},
we have $\La(L\hconv M,L)\le\La(M,L)-2d_L$ by Proposition~\ref{prop:ddd}.
On the other hand, the same proposition implies that 
$$\La(L\hconv M,L)+\La(\D L,L)-\La(M,L)\ge0$$
since $(L\hconv M)\hconv\D L\simeq M$.
Since $\La(\D L,L)=2\sfd_L$,  we have
$\La(M,L)-2d_L\le \La(L\hconv M,L)$. Therefore, we obtain $\La(L\hconv M,L)=\La(M,L)-2d_L$.
The other assertions are obtained similarly.
\QED

\Prop\label{prop:LL*}
Let $L$ be a root object of $\tCw$, and $M$ a simple object in $\tCw$.
If \/ $\de(L,M)\ne d_L$, then we have
$$(L\hconv M)\hconv L\simeq L\hconv(M\hconv L).$$
\enprop
\Proof
If $\de(L,M)=0$, then $L$ and $M$ commute, and hence
$L\conv M\conv L$ is simple.
Thus, we may assume that $\de(L,M)\ge 2d_L$
by Proposition~\ref{prop:ddd}.
Let $K$ be the kernel of $M\conv L\to M\hconv L$, 
i.e., the radical of $M\conv L$.
Then we have a commutative diagram with an exact row:
$$\xymatrix{
0\ar[r]&L\conv K\ar[r]^{\iota\,\,\,}\ar[dr]_f&L\conv M\conv L\ar[r]^{p}\ar@{->>}[d]_{\pi'}& L\conv (M\hconv L)\ar[r]
\ar@{-->>}[ld]^{\pi}&0\\
&&(L\hconv M)\hconv L}
$$
Let us show that
the composition $f\cl L\conv K\to(L\hconv M)\hconv L$ vanishes. Indeed, if it holds,
we obtain an epimorphism $\pi\cl L\circ(M\hconv L)\epito(L\hconv M)\hconv L$, which 
shows the claim. 

In order to see the vanishing of $f$, it is enough to show that
$L\hconv S\not\simeq (L\hconv M)\hconv L$
for any simple subquotient $S$ of $K$ since  $L\hconv S$ covers the composition series of 
$L\circ K$.
We have 
$$\La(S,L)<\La(M,L)=\La\bl L\hconv M, L)+2d_L=
\La\bl(L\hconv M)\hconv L,L\br+2d_L.$$
Here, the first inequality follows from Proposition~\ref{prop:simplehd},
and the equalities follow from
Proposition~\ref{prop:root}.
If $\de(L,S)>0$, then we have
$\La(L\hconv S,L)=\La(S,L)-2d_L$ by
Proposition~\ref{prop:root}.
we have $L\hconv S\not\simeq (L\hconv M)\hconv L$. 
Hence we may assume that $\de(L,S)=0$ or equivalently, $L$ commutes with $S$.
Then, $S\hconv L\simeq (L\hconv M)\hconv L$
implies that
$S\simeq L\hconv M$, which derives a contradiction $
0=\de(L,S)=\de(L,L\hconv M)=\de(L,M)-d_L\geq 2d_L-d_L>0$.
\QED

The following proposition will be used in  
this section and \S\;\ref{sec:proof}.

\Prop\label{prop:simpleroot}
Assume that $i\in I$ satisfies either $s_iw<w$ or $ws_i<w$. Then,
$\Qti$ is a root object or invertible.
\enprop
\Proof
Since $\Qti$ has an affinization of degree $2\sfd_i$ by
Proposition~\ref{prop:aff-simple},
it is enough to show that
\eq\de\bl\Qti,\D^{-1}\Qti\br\le\sfd_i.
\label{eq:rootin}
\eneq
Set $M=\dM(w\La_i,\La_i)$. Then
$\D^{-1}\Qti\simeq\tdC_{-\La_i}\conv\Qt(\Es M)$ 
by Theorem~\ref{thm:gdm}\eqref{item: sla}.
Hence we obtain
$$\de\bl\Qti,\D^{-1}\Qti\br\le\de\bl\ang{i},\Es M\br=\sfd_i\de_i(\Es M)$$
by Lemma~\ref{lem:LaQt}.

By  Proposition~\ref{Prop: dM properties},
we have $\eps^*_i(\Es M)=0$ and $\eps_i(\Es M)\le\max\bl0,-\ang{h_i,w\La_i}\br$.
Hence, we obtain 
\eqn
\de_i(\Es M)&=&\eps_i(\Es M)+\epss_i(\Es M)
 +\ang{h_i, \wt(\Es M)}\\
 &\le&\max\bl0,-\ang{h_i, w\La_i}\br+\ang{h_i,w\La_i-s_i\La_i}
 =\max\bl\ang{h_i,w\La_i-s_i\La_i}, 1\br.
 \eneqn

If $ws_i<w$, then we have $w\La_i-s_i\La_i\in\nrtl$ and hence
$ \ang{h_i,w\La_i-s_i\La_i}\le0$.
Hence we obtain \eqref{eq:rootin}.

If $s_iw<w$, then we have
$w\al_i\in \rt_-$ and hence
$ \ang{h_i,w\La_i-s_i\La_i}=\ang{w^{-1}h_i,\La_i}+1\le1$,
which implies \eqref{eq:rootin}.
\QED

\Lemma\label{lem:Iw'}
Assume that $i\in I$ satisfies $w'\seteq ws_i<w$ and $i\not\in I_{w'}$.
Then $\Qt[w](\ang{i})$ is invertible.
\enlemma
\Proof
Set $\la=\La_i+s_i\La_i=2\La_i-\al_i\in\pwtl$.
Then, $\ang{h_j,s_i\La_i}=\ang{h_j,\la}\in\Z_{\ge0}$ for any $j\in I_{w'}$, and we have
$$\dM(w's_i\La_i,s_i\La_i)\simeq \dM(w'\la,\la)\simeq\dM(w\la,\la).$$
Hence, we have $\dM(w\la,\la)\hconv\ang{i}\simeq\dM(w\La_i,\La_i)$, which implies that
$\Qt(\ang{i})\simeq \tdC_\la^{\circ-1}\conv\tdC_{\La_i}\simeq\tdC_{-s_i\La_i}$
is invertible.
\QED

\Rem
In general, $\Qti$ may be neither a root object nor an invertible object.
For example, take $\g=A_3$, $w=s_2w_0=s_1s_3s_2s_1s_3$, and $i=2$.
Then $\Qti\simeq\ang{132}^{-1}\conv \ang{12}\conv\ang{32}$ and
$\D^{-1}\Qti\simeq \ang{1}\conv\ang{3}\conv\ang{132}^{-1}$.
Hence $\de\bl\Qti,\D^{-1}\Qti\br=2$.
\enrem

\subsection{Simple root operators}\label{root-op}
Now we shall define the simple root operators on the localized category.
In this case, the simple root operators are bijective.

\Def\label{def:rootop} For any $i\in I_w$ and a simple $X\in\tCw$, we define
\eqn
\eps_i(X)&&=\sfd_i^{-1}\tLa(\Qti,X)\in\Z,\\
\eps^*_i(X)&&=\sfd_i^{-1}\tLa(X,\Qti)\in\Z,\\
\vphi_i(X)&&=\eps_i(X)+\ang{h_i, \wt X}\in\Z,\\
\vphi^*_i(X)&&=\eps^*_i(X)+\ang{h_i, \wt X}\in\Z,\\
\de_i(X)&&=\sfd_i^{-1}\de(\Qti,X)=\eps_i(X)+\eps_i^*(X)+\ang{h_i,\wt X}\in\Z_{\ge0},\\
\tF_i X&&=q_i^{\eps_i(X)}\Qti\hconv X\in\tCw,\\
\tF^*_i X&&=q_i^{\eps^*_i(X)}X\hconv\Qti\in\tCw,\\
\tE_i X&&=q_i^{\vphi_i(X)+1}X\hconv\D\Qti\in\tCw,\\
\tE^*_i X&&=q_i^{\vphi^*_i(X)+1}\D^{-1}\Qti\hconv X\in\tCw.
\eneqn
For $i\in I\setminus I_w$, we understand
$\tF_i X=\tE_i X=0$, $\eps_i(X)=\eps_i^*(X)=\de_i(X)=-\infty$.
\edf

The following proposition is immediate from the definition above and
Lemma~\ref{lem:invhvonv}. 
\Prop\label{pro:inverse}
$\tF_i$ and $\tE_i$ are inverse to each other.
Similarly, $\tF^*_i$ and $\tE^*_i$ are inverse to each other.
\enprop
One of the main aims is to show that these data define a crystal structure on 
$\tCw$ the set of equivalence classes of  simples in $\tCw$ up to grading shifts.

\Th\label{thm:crystal}
The data in Definition~\ref{def:rootop} defines a crystal structure on $\Irr(\tCw)$.
\enth
\Proof
We have all the conditions in Definition~\ref{cryst}  except (2) and (3). 
Thus, it suffices to show that 
if $X\in \Irr(\tCw)$, then $\wt(\tF_i X)=\wt(X)-\al_i$, 
$\vep_i(\tF_i X)=\vep_i(X)+1$, $\vp_i(\tF_i X)=\vp_i(X)-1$, which implies Definition~\ref{cryst} (3).
The condition (2) is obtained by the fact $\tE_i=\tF_i^{-1}$.

\medskip
Now, by Proposition~\ref{prop:root}(i) we have
$\La\bl\Qt(\ang{i}),\Qt(\ang{i})\hconv X\br=\La\bl\Qt(\ang{i}),X\br$ and then
\eqn
&&2\tLa\bl\Qt(\ang{i}),\Qt(\ang{i})\hconv X\br\\ 
 &&\hs{10ex}=\La(\Qt(\ang{i},X)
+(\wt(\Qt(\ang{i}),\wt(\Qt(\ang{i})+\wt(X))\\
&& \hs{10ex}=2\tLa(\Qt(\ang{i},X)
+(\wt(\Qt(\ang{i})),\wt(\Qt(\ang{i})))
=2\tLa(\Qt(\ang{i},X)+2\sfd_i,
\eneqn
which means $\vep_i(\tF_i X)=\vep_i(X)+1$ and
the other assertions are obtained by the definitions.
\QED

%
\Prop
For any $i\in I_w$ and any simple $X\in\tCw$, we have
\bnum
\item
  $\vphi_i(X)=\sfd_i^{-1}\tLa\bl X,\D\Qt(\ang{i})\br$,
  \item
  $\vphi^*_i(X)=\sfd_i^{-1}\tLa\bl \D^{-1}\Qt(\ang{i}),X\br$.
\ee
\enprop
\Proof
This follows from
$\tLa(M,N)=\tLa(N,\D M)+(\wt\, M,\wt\, N)=\tLa(\D^{-1}N,M)+(\wt\, M,\wt\, N)$ 
for simple $M$, $N\in\tCw$ such that one of them is \afr (see Lemma~\ref{lem:MNDM}).
\QED


\Lemma\label{lem:tEcomp}
Let $M\in\CBw$ be a simple module and $i\in I_w$. Then we have
\bnum
\item $\tE_i M$, $\tE^*_i M\in\CBw$,
\item If $\tE_i M\not\simeq0$, then $\Qt(\tE_i M)\simeq
\tE_i\bl\Qt(M)\br$,
\item If $\tEs_i M\not\simeq0$, then $\Qt(\tEs_i M)\simeq
\tEs_i\bl\Qt(M)\br$.
\ee
\enlemma
\Proof
We prove only the statements for $\tE_i$.
We may assume  $\tE_i M\not\simeq0$.
Then we have
$$\ang{i}\conv \tE_iM\epito M\monoto \tE_iM\conv \ang{i}.$$
Applying  $\Qt$, we obtain
$$\Qti\conv \Qt(\tE_iM)\epito \Qt(M)\monoto \Qt(\tE_iM)\conv \Qti.$$
Hence $\Qt(\tE_iM)\not\simeq0$, and then we have
$\Qti\hconv \Qt(\tE_iM)\simeq\Qt(M)$, or equivalently,
$\tF_i\bl \Qt(\tE_iM)\br\simeq \Qt(M)$.
It implies that
$\tE_i\Qt(M)\simeq\Qt(\tE_i M)$.
\QED

Now, we can describe $\D\bl\Qti$ explicitly by using the generalized determinantial modules.
\Lemma
For any $i\in I_w$, we have
$$\D\bl\Qti\br\simeq \Qt\bl \E_i\dM_{w}(-\La_i,-w^{-1}\La_i)\br
\conv\tdC_{w^{-1}\La_i}.$$
\enlemma
\Proof
Note that Lemma~\ref{lem:*-inv} implies
\[
\E_i\dM_{w}(-\La_i,-w^{-1}\La_i)\simeq\psi_*\Bigl(\Es
\dM(w^{-1}\La_i,\La_i)\Bigr)
\simeq\psi_*\Bigl(\dM(w^{-1}\La_i,s_i\La_i)\Bigr).
\]
The assertion follows from
$\ang{i}\conv\bl \E_i\dM_{w}(-\La_i,-w^{-1}\La_i)\br
\epito \dM_{w}(-\La_i,-w^{-1}\La_i)$ and
$\Qt\bl\dM_{w}(-\La_i,-w^{-1}\La_i)\br\simeq\tdC_{-w^{-1}\La_i}$
by \eqref{eq:gendet}.
\QED
\Lemma\label{lem:normalij}
Let $i,j\in I_w$ such that $i\not=j$. 
Assume further that either $s_iw<w$ or $ws_j<w$. Then 
$\D\Qti$ and $\Qti[j]$ commute.
\enlemma
\Proof
Assume first $s_iw<w$. 
Then, we have $\gW(\ang{i})\setminus\st{0}=\{\al_i\}\subset\Delta_+\cap w\nrt$ and hence
$\ang{i}\in\Cw$. 
Then $\dM(w\La_j,\La_j)$ commutes with $\ang{i}$.
Hence $\Esm[j]\dM(w\La_j,\La_j)\simeq\dM(w\La_j,s_j\La_j)$
commutes with $\Esm[j]\ang{i}\simeq\ang{i}$,
which implies that
$\D^{-1}\Qti[j]\simeq \tdC_j^{-1}\conv\Qt\bl\Esm[j]\dM(w\La_j,\La_j)\br$
commutes with $\Qti$.

\mnoi
Now assume that $ws_j<w$.
Then $s_jw^{-1}<w^{-1}$.
Hence
$\D\bl\Qt[w^{-1}](\ang{j})\br$ and $\Qt[w^{-1}](\ang{i})$ commute.
By applying $\psi_*$, we obtain that
$\D^{-1}\bl\Qt\ang{j}\br$ and $\Qti$ commute by Theorem~\ref{th:main1}.
\QED

The following lemma is an immediate consequence of
Proposition~\ref{prop:root}, Proposition~\ref{prop:LL*}
and Proposition~\ref{prop:simpleroot}.
\Lemma\label{lem:iistar} Assume that $i\in I$ satisfies $ws_i<w$. 
Let $X$ be a simple object of $\tCw$.
Then, we have
\bnum
\item if $\de_i(X)>0$, then
  $\eps_i\bl(\tEs_i)^n(X)\br=\eps_i(X)$ for any $n\in\Z_{\ge0}$
  and $\tF_i\tEs_i(X)\simeq\tEs_i\tF_i (X)$,

\item
  $\eps_i(\tEs_i X)=\eps_i(X)$ or $\eps_i(\tEs_i X)=\eps_i(X)-1$.
\ee
\enlemma
Note that in (i), $\Qti$ is a root object.
\Lemma\label{lem:ij*}
Let $i,j\in I_w$ such that $i\not=j$. 
Assume further that either $s_iw<w$ or $ws_j<w$. Then we have
\bnum
\item
$\tF_i$ and $\tF^*_j$ commute as operators on $\Irr(\tCw)$,
\item 
$\eps_i(\tF^*_jX)=\eps_i(X)$ and
$\eps^*_j(\tF_iX)=\eps^*_j(X)$ for any $X\in\Irr(\tCw)$.\label{it:ij*}
\ee
\enlemma
\Proof 
By Lemma~\ref{lem:normalij} and Lemma~\ref{lem:normalD},
$(\Qti,X,\Qti[j])$ is an almost \afr normal sequence. 
Hence we have
$$(\Qti\hconv X)\hconv\Qti[j]\simeq \Qti\hconv(X\hconv\Qti[j])$$
and $$\tLa(\Qti,X\hconv\Qti[j])=\tLa(\Qti,X)+\tLa(\Qti,\Qti[j]).$$
Note that
$\tLa(\Qti,\Qti[j])=\tLa(\ang{i},\ang{j})=0$ since $\ang{i\,j}\in\CBw$.
\QED

\Lemma
Let $M\in\CBw$ be a simple module.
Then we have
\bnum
\item if $s_iw<w$, then we have
$\tF_i(\Qt(M))\simeq\Qt(\tF_i(M))$ and $\eps_i(\Qt(M))=\eps_i(M)$,
\item if $ws_i<w$, then we have
$\tF^*_i(\Qt(M))\simeq\Qt(\tF^*_i(M))$
and $\eps^*_i(\Qt(M))=\eps^*_i(M)$.
\ee
\enlemma
\Proof
(i) Since $\ang{i}\in\Cw$, it follows from Proposition~\ref{prop:CB}.

\snoi
(ii) follows from (i) for $w^{-1}$ by applying $\psi_*$ (see
Theorem~\ref{th:main1}).
\QED

\Rem
In general, $\eps_i(\Qt(M))$ can be negative for $M\in\Cw$
(cf.\ Corollary~\ref{cor:epsCw}).
The equality $\eps_i\bl\Qt(M)\br=\eps_i(M)$ {\em does not hold} for every
$i\in I_w$ and every simple $M\in\Cw$. For example, take
$\g=A_2$, $w=s_1s_2$, $i=2$ and $M=\ang{1}$.
Then, $\Cw$ is generated by its central objects $\ang{1}$ and $\ang{12}$.
We have $\ang{1}\hconv \ang{2}\simeq \ang{12}$ in $R\gmod$ as in 
Remark~\ref{rem:12}
and then we 
have
$\Qt(\ang{1})\conv\Qt(\ang{2})\simeq\Qt(\ang{12})$ in $R\gmod[C_i^{\circ-1};i\in I]\simeq\tCw$
(see Theorem~\ref{thm:rigidity-Cw}). Thus, we have 
$\Qt(\ang{2})=\Qt(\ang{1})^{-1}\conv\Qt(\ang{12})$ and then
\eqn
\eps_2\bl\Qt(\ang{1})\br&=&
\sfd_2^{-1}\tLa(\Qt(\ang{1})^{-1}\conv\Qt(\ang{12}),\, \Qt(\ang{1}))\\
&=&-\sfd_2^{-1}\tLa(\Qt(\ang{1}),\Qt(\ang{1}))+\sfd_2^{-1}\tLa(\Qt(\ang{12}),\Qt(\ang{1}))\\
&=&-\tLa(\ang{1},\ang{1})+\tLa(\ang{12},\ang{1})=-1+0=-1,
\eneqn
where $\sfd_2=1$. However, we have $\eps_2(\ang{1})=0$.

Another example is $\g=A_3$, $w=s_1w_0=s_2s_1s_3s_2s_1$,
$i=1$ and $M=\ang{23}$.
Then we have $\Cw=\st{Z\in R\gmod\mid \E_1(M)\simeq0}$,
$\Qt(\ang{1})=\ang{23}^{\circ-1}\conv \ang{231}$.
We have
$\eps_1(M)=0$ and $\eps_1\bl\Qt(M)\br=\tLa(\ang{231},\ang{23})-\tLa(
\ang{23},\ang{23})=-1$.
Indeed, since the pair $(\ang{231},\ang{23})$ is unmixed,  
\cite[Corollary~2.3]{KKOP22}) implies that $\tLa(\ang{231},\ang{23})=0$.
\enrem 

\subsection{The map $\Est$}\label{map-E*}
{\em In this subsection, let $w\in\weyl$ and
  assume that $i\in I$ satisfies $w'\seteq ws_i<w$.}
By Lemma~\ref{lem:CBww'} we can define
$$\tEsm_i\cl\Irr(\CBw)\to\Irr(\CBw[w'])$$
by 
$$
M\mapsto \tEsm_{i}(M).
$$
Then it induces
$$\Irr(\Cw)\to\Irr(\tCw[w'])$$
by
$$M\longmapsto \Qt[{w'}]\bl\tEsm_{i}(M)\br.$$
Note that by the shuffle lemma (Lemma~\ref{lem:E*}) we have
for any simple $M\in\Cw$ and $\La\in\pwtl$,
\eq
\tEsm_{i}(\dC_\La\conv M)\simeq\tEsm_{i}(\dC_\La)\conv
\tEsm_{i}(M)
\label{La-M}
\eneq

Hence,
it extends to 
$$\Est\cl\Irr(\tCw)\to\Irr(\tCw[w'])$$
 by
$$(\tdC_\La)^{-1}\conv \Qt(M)\longmapsto \Qt[{w'}]\bl\tEsm_{i}(\dC_\La)\br^{-1}\conv\Qt[{w'}]\bl\tEsm_{i}(M)\br \qt{for $M\in\Irr(\Cw)$ and $\La\in\pwtl$.}$$
Note that $\Qt[{w'}](\tEsm_{i}\dC_\La)$ is invertible by \eqref{eq:center}:
\eqn
\Qt[w']\bl\tEsm_{i}\dC_\La\br&&\simeq\Qt[w']\bl\dM(w\La,s_{i}\La)\br
\simeq\Qt[w']\bl\dM_{w'}(w\La,s_{i}\La)\br\simeq\tdC^{w'}_{s_i\La}.
\eneqn

Together with the map
\[
\eps^*_{i}\cl \Irr(\tCw)\to\Z,
\]
we obtain 
\eq
\Irr(\tCw)\to \Irr(\tCw[w'])\times\Z 
\qquad
(X\mapsto (\Est(X),\eps^*_i(X)).
\eneq

\Rem
For $M\in\Cw$ and $i\in I$ such that $w'\seteq ws_i<w$,
the statement $\Esm M\in\Cw[w']$ may not hold.
For example, take $\g=A_2$, $w=s_1s_2s_1$. Then, we have $\Cw=R\gmod$.
Take $i=1$ and
$M=\ang{2}$. Then $\Esm M\simeq M\not\in\Cw[w']=\st{M\mid E_2M\simeq0}$.
\enrem

\Prop\label{prop:ELM}
Let $L\in\Cw$ be an \afr simple module and $M\in\CBw$ a simple module.
Set $L'=\Est(L)$ and $M'=\Est(M)$.
Assume that $\eps_{i}^*(L\hconv M)=\eps_{i}^*(L)+\eps_{i}^*(M)$.
Then we have
\bnum
\item 
$\Est(L\hconv M)\simeq\ L'\hconv M'$ up to grading shifts,
\item
$\La(L', M')=\La(L,M)+\eps_{i}^*(L)\bl\al_{i},\wt(M)\br-\eps_{i}^*(M)\bl\al_{i},\wt(L)\br$.
\ee
\enprop
\Proof
It follows from Proposition~\ref{prop:E*}.
\QED
Recall that we have assumed that $w'\seteq ws_i<w$ in this subsection.
\Prop\label{prop:CL}
Let $M\in\CBw$ and $\La\in\pwtl$.
Set $\dC''_\La\seteq\Esm(\dC_\La)\simeq\dM(w\La,s_{i}\La))
\simeq\dC^{w'}_{s_i\La}\in\Cw\cap\CBw[w']$ and $M'\seteq\Esm(M)\in\CBw[w']$.
Then we have
\bnum
\item $\eps_{i}^*(\dC_\La\hconv M)=\eps_{i}^*(\dC_\La)+\eps_{i}^*(M)$,
\item $\dC''_\La\hconv M'\in\CBw[w']$,
\item 
$\Esm(\dC_\La\hconv M)\simeq\dC''_\La\hconv M'$,
\item
$\La(\dC''_\La, M')=-\bl w\La+s_i\La,\wt(M')\br$.\label{it:it4}
\ee
\enprop
\Proof
Since $M\hconv\ang{i_{}}\in\CBw$, Proposition~\ref{prop:binBw} implies that
$$\La(\dC_\La,M\hconv\ang{i_{}})
=-\bl w\La+\La, \wt(M\hconv\ang{i_{}})\br
=\La(\dC_\La,M)+\La(\dC_\La,\ang{i_{}}).$$
Hence $(\dC_\La,M,\ang{i_{}})$ is a normal sequence.
Hence we have
$$\La(\dC_\La\hconv M,\ang{i_{}})=
\La(\dC_\La,\ang{i_{}})+\La(M,\ang{i_{}}),$$
which is nothing but (i).
Hence Proposition~\ref{prop:E*} (iii) implies (iii).
Since $\dC_\La\hconv M\in\CBw$, (iii) implies (ii) by Lemma~\ref{lem:CBww'}.

It remains to prove (iv).
We have
\eqn
\La(\dC''_\La,M')&&=\La\bl\Qt[w'](\dC''_\La),\Qt[w'](M'))
=\La\bl\tdC^{w'}_{s_i\La},\Qt[w'](M'))
  =-(w's_i\La+s_i\La,\wt(M')\br.
  \eneqn
  Here the first equality follows from (ii) and Lemma~\ref {lem:LMN}.
\QED

\Cor\label{cor:CBE}
For any simple $M\in\CBw$, we have
$$\Qt[w'](\Esm(M))=\Est(\Qt(M)).$$
That is, the following diagram commutes;
$$\xymatrix@C=10ex{
\Irr(\CBw)\ar[dr]^{\Qt[w']\,\circ\,\tEsm_i}\ar[d]_{\Qt\vert_{\Irr(\CBw)}}\\
\Irr(\tCw)\ar[r]_{\Est}&\Irr(\tCw[w']).}$$
\encor
\Proof
Take $\La\in\pwtl$ such that $N\seteq\dC_\La\hconv M\in\Cw$.
Since $\Qt(M)\simeq\Qt(\dC_\La)^{-1}\conv\Qt(N)$, we have
\eqn
\Est(\Qt(M))&&\simeq\Est(\Qt(\dC_\La))^{-1}\conv\Est(\Qt(N))
\simeq\Qt[w'](\dC''_\La)^{-1}\conv\Qt[w']\bl\Esm(N)\br\\
&&\simeq\Qt[w'](\dC''_\La)^{-1}\conv\Qt[w'](\dC''_\La\hconv\Esm(M)\br\\
&&\simeq\Qt[w'](\dC''_\La)^{-1}\conv\Qt[w'](\dC''_\La)\conv\Qt[w']\bl\Esm(M)\br
\simeq\Qt[w']\bl\Esm(N)\br.
\eneqn
\QED

\Lemma\label{lem:EFX=EX}
Let  $M\in\CBw$ be a simple module and let $X$ be a simple object of $\tCw$. Then, we have
\bnum
\item $\eps^*_i\bl\Qt(M)\br=\eps^*_i(M)$,
\item$\tF^*_i\Qt(M)\simeq \Qt(M\hconv\ang{i})$,
\item $\Est(\tF^*_i(X))\simeq\Est(X)$.
\ee
\enlemma
\Proof
Since $M\hconv\ang{i}\in\CBw$, Lemma~\ref{lem:LMN}
implies that
$\La\bl\Qt(M),\Qti\br=\La(M,\ang{i})$ and
$\Qt(M)\hconv\Qti\simeq\Qt(M\hconv\ang{i})$,
which reads as (i) and (ii).

\mnoi
(iii)\
Let $X=(\tdC_\La)^{\circ-1}\conv \Qt(M)$ for a simple $M\in\Cw$
and $\La\in\pwtl$.
Then we have
$\tFs_i(X)\simeq(\tdC_\La)^{\circ-1}\conv \Qt(\tFs_iM)$ and
\eqn\Est(\tF^*_i(X))
&&\simeq(\tdC^{w'}_{s_i\La})^{\circ-1}\conv \Qt(\tEsm_i\tFs_iM)\\
&&\simeq(\tdC^{w'}_{s_i\La})^{\circ-1}\conv \Qt(\tEsm_iM)
\simeq\Est(X).
\eneqn\
\QED

\subsection{Crystals}
Let us take a reduced expression $\uw=s_{i_1}\cdots s_{i_\ell}$ of $w\in\weyl$.
We set
$$
w_{\le k}=s_{i_1}\cdots s_{i_{k}}\qtq
\beta_k=w_{\le k-1}\al_{i_k}=w_{<k}\La_{i_k}-w_{\le k}\La_{i_k}\in\prt.
$$
It is well-known that $\prt\cap w\nrt=\st{\beta_1,\ldots,\beta_\ell}$.

Recall the crystal $\Bw$ introduced in \S\,\ref{cellular}:
$$\Bw\seteq B_{i_1}\tens\cdots\tens B_{i_\ell}.$$

For a simple $M\in \CBw$, we define
$\CP(M)\in\Z^\ell$ as follows:
\eqn
&&M_\ell=M,\\
&&c_k=\eps^*_{i_k}(M_k)\qt{for $1\le k\le \ell$,}\\
&&M_{k-1}=(\tEs_{i_k})^{c_k}(M_k), \\
&&\CP(M)=(c_1,\ldots,c_\ell).
\eneqn
We also regard $\CP(M)$ as an element of $\Bw$
by $(c_1,\ldots,c_\ell)\mapsto \tf_{i_1}^{c_1}(0)_{i_1}\tens\cdots\tens\tf_{i_\ell}^{c_\ell}(0)_{i_\ell}\in\Bw$.

Then $\CP$ induces
$$\Irr(\Cw)\monoto \Irr(\CBw)\monoto \Bw.$$

Since by Proposition~\ref{prop:CL}\;(i) (see also \cite[Sect.8]{N2}) we have
\eq
\CP(\dC_\La\hconv M)=\CP(\dC_\La)+\CP(M)\qt{for any $M\in\Irr(\Cw)$ and $\La\in\pwtl$,}\eneq
we can extend  $\Irr(\Cw)\monoto \Bw$ to $\Irr(\tCw)\monoto \Bw$ by
\eqn
\tdC_\La^{-1}\conv \Qt(M)\longmapsto \CP(M)-\CP(\dC_\La)\qt{for $\La\in\pwtl$,}
\eneqn
which we write by the same letter $\CP$.
Then we obtain the commutative diagram:
\eqn
\xymatrix@C=10ex{
\Irr(\Cw)\ar@{^{(}->}[d]\ar@{^{(}->}[r]&\Irr(\CBw)\ake\ar@{^{(}->}[d]\ar[dl]|{\ake\ake[2.5ex]\CP\circ\Qt}\\
\Irr(\tCw)\ar[r]_{\CP}&\Bw.
}
\eneqn
By Proposition~\ref{Prop: dM properties}, we have that
$$\CP(\tdC_\la)=\st{-\bl \beta_k^\vee, w\la\br}_{1\le k\le \ell}\qt{for any $\la\in\wtl$.}$$
The following lemma is obtained in \cite{Kana-N} in the case when $\g$ is semi-simple and $w$ is the longest element.
\Lemma\label{lem:fC}
For any $\la\in\wtl$, $i\in I$ and $b\in\Bw$, we have
\bnum
\item $\tf_i\bl\CP(\tdC_\la)+b\br=\CP(\tdC_\la)+\tf_i(b)$,
\item
$\eps_i(\CP(\tdC_\la))=-\ang{h_i,w\la}$.
\ee
\enlemma

\Proof
Set $\CP(\dC_\la)=(c_1,\ldots,c_\ell)$.
Then we have
$c_k=-(\beta^\vee_{i_k}, w\la)$.
It is enough to show that we have
$a_k=-\ang{h_{i_k}, w\la}$, where
\eqn
a_k&&\seteq c_k+\sum_{1\le s<k}\ang{h_{i},\al_{i_s}}c_s=
-(\beta_k^\vee, w\la)-\sum_{1\le s<k}\ang{h_{i},\al_{i_s}}(\beta_s^\vee, w\la).
\eneqn
Indeed, we have seen in \S\;\ref{cellular} that 
\[
\eps_i(\CP(\tdC_\la))=\max\{a_k\mid i_k=i\},
\]
which implies (ii).
Since
\[
  \TY(\til m, i,f)(\CP(\tdC_\la)+b)=\TY(\til m, i,f)(b),
  \]
we obtain the assertion (i).

\medskip
Set $i=i_k$, $A_k=\beta_k+\sum_{1\le s<k}(\al_{i},\al_{i_s})\beta_s^\vee$.
Then, we have
$a_k=-(A_k,w\la)/\sfd_{i}$.
Hence it is enough to show that
$A_k=\al_i$, that is, $A_k$ does not depend on $k$.

We have
\eqn
A_k&&=\beta_k+\sum_{1\le s<k}(\al_{i},\al_{i_s}^\vee)\beta_s
=\beta_k+\sum_{1\le s<k}w_{\le s-1}(\al_{i},\al_{i_s}^\vee)\al_s\\
&&=\beta_k+\sum_{1\le s<k}w_{\le s-1}(\al_{i}-s_{i_s}\al_{i})
=\beta_k+\sum_{1\le s<k}(w_{\le s-1}\al_{i}-w_{\le s}\al_{i})\\
&&=\beta_k+\al_{i}-w_{\le k-1}\al_{i}=\al_{i}.
\eneqn
\QED
Note that $\eps_i(\dC_\la)=\max(-\ang{h_i,w\la},0)$ for
any $\la\in\pwtl$ and $i\in I_w$,
while we will see that
$\eps_i(\tdC_\la)=-\ang{h_i,w\la}$.

The following proposition immediately follows from the definitions of $\Est(X)$ and $\CP(X)$.
\Prop\label{prop:CPww'}
Let $X\in\Irr(\tCw)$.
Set $X'=\Est(X)\in\Irr(\tCw[w'])$,
$\CP(X)=(c_1,\ldots,c_\ell)$.
Then, we have
\bnum
\item $c_\ell=\eps_{i_\ell}^*(X)$,
\item$\CP[\uw'](X')=(c_1,\ldots,c_{\ell-1})$.
\ee
\enprop

\section{Main Theorem}
Let $\uw=s_{i_1}\cdots s_{i_\ell}$ (resp.\ $\uw'=s_{i_1}\cdots s_{i_{\ell-1}}$) 
be a reduced expression of $w\in\weyl$ (resp.\ $w'\seteq ws_{i_{\ell}}\in\weyl$). 
In this section, we shall prove the following Main Theorem~\ref{th:Main},
which says that $\Irr(\tCw)\to\Bw$ is a morphism of crystals.

\MTh\label{th:Main}
For any simple $X\in\tCw$ and $i\in I$,  we have
\bnum
\item 
$$\eps_i(X)=\bc
\eps_i(\Est(X))&\text{if $i\not=i_\ell$,}\\
\max\bl \eps_i(\Est(X)), -\eps_{i_\ell}^*(X)-\ang{h_i,\wt(X)}\br&
\text{if $i=i_\ell$.}
\ec$$\label{MT:1}
\item $\Est\bl\tF_i(X)\br\simeq\tF_i\bl\Est(X)\br$ if $i\not=i_\ell$.
  \label{MT:2}
\item if $i=i_\ell$,  then
\eqn
\Est\bl\tF_i(X)\br\simeq 
\bc\tF_i\bl\Est(X)\br&\text{if $\de_i(X)>0$ and $i\in I_{w'}$,}\\
\Est(X)&\text{otherwise,}
\ec
\eneqn\label{MT:3}
\item$\eps^*_{i_\ell}(\tF_i(X))=
\bc\eps^*_{i_\ell}(X)+1&\text{if $i=i_\ell$ and $\de_{i_\ell}(X)=0$,}\\
\eps^*_{i_\ell}(X)&\text{otherwise.}
\ec$
\label{MT:4}
\ee
\enmth
Note that 
\eq
&& -\eps_{i_\ell}^*(X)-\ang{h_{i_\ell},\wt(X)}=
\eps_{i_\ell}^*(X)-\ang{h_{i_\ell},\wt(\Est X)}.
\label{eq:e*E}
\eneq
The proof of the main theorem will be given in the following section.

Main Theorem implies the following result.
\Cor 
$\CP\cl\Irr(\tCw)\to\Bw$ is a morphism of crystals.
\encor
\Proof
By induction on the length of $w$, it is enough to show that
$\Irr(\tCw)\to \Irr(\tCw[w'])\tens B_{i_\ell}$ is a morphism of crystals.
It  is evident $\wt(X)=\wt(\CP(X))$. Thus, by Main Theorem~\ref{th:Main} (i), we have 
$\varepsilon_i(X)=\varepsilon_i(\CP(X))$ by \eqref{tensor-vep} and \eqref{eq:e*E}, 
and  also get $\varphi_i(X)=\varphi_i(\CP(X))$.

Next, let us see the compatibility of the actions by $\tf_i$.
By the definition, we have
\[
\CP(X)=\CPP(\Est(X))\otimes \tf_{i_\ell}^{\vep^*_{i_\ell}(X)}(0)_{i_\ell}.
\]
 Then, if $i\ne i_{\ell}$, by induction on the length of $w$ and 
Main Theorem~\ref{th:Main}\;\eqref{MT:2}, we obtain
\eqn
\tf_i\CP(X)&=&\tf_i\CPP(\bl\Est(X)\br)\otimes \tf_{i_\ell}^{\vep^*_{i_\ell}(X)}(0)_{i_\ell}
=\CPP(\tF_i\bl\Est(X)\br)\otimes \tf_{i_\ell}^{\vep^*_{i_\ell}(X)}(0)_{i_\ell}\\
&=&\CPP(\Est(\tF_i(X)))\otimes \tf_{i_\ell}^{\vep^*_{i_\ell}(\tF_i(X))}(0)_{i_\ell}
=\CP(\tF_i(X)),
\eneqn
where $\vep^*_{i_\ell}(\tF_i(X))=\vep^*_{i_\ell}(X)$ follows from
Main Theorem~\ref{th:Main}\;\eqref{MT:4} and
$i\ne i_\ell$. 

Assume that $i=i_\ell$. We have  
\eqn
\de_{i_\ell}(X)&&=\eps_{i_\ell}(X)
+\eps^*_{i_\ell}(X)+\ang{h_{i_\ell},\wt(X)}\\
&&=\max\bl\eps_{i_\ell}(\Est(X)),\ -\eps_{i_\ell}^*(X)-
\ang{h_{i_\ell},\wt(X)}\br
+\eps^*_{i_\ell}(X)+\ang{h_{i_\ell},\wt(X)}\\
&&=\max\bl\eps_{i_\ell}(\Est(X))-\eps_{i_\ell}^*(X)+\ang{h_{i_\ell},\ \wt(\Est X)},0\br.
\eneqn
Hence $\de_{i_\ell}(X)=0$ if and only if 
$\eps_{i_\ell}(\Est(X))\le\eps_{i_\ell}^*(X)-\ang{h_{i_\ell},\wt(\Est X)}$, which is equivalent to  
\[
\varphi_{i_\ell}(\CPP(\Est(X))-\vep_{i_\ell}(\tf_{i_\ell}^{\vep^*_i(X)}(0)_{i_\ell})
=\eps_{i_\ell}(\Est(X))-\eps_{i_\ell}^*(X)+\ang{h_{i_\ell},\wt(\Est X)}\leq0.
\]
Therefore, under the condition $\de_{i_\ell}(X)=0$, \eqref{tensor-f} and Theorem~\ref{th:Main}\;\eqref{MT:4} imply that
\eqn
\tf_{i_\ell}\CP(X)&=&\tf_{i_\ell}(\CPP(\Est(X))\ot \tf_{i_\ell}^{\vep^*_{i_\ell}(X)}(0)_{i_\ell})
=\CPP(\Est(X))\ot \tf_{i_\ell}^{\vep^*_{i_\ell}(X)+1}(0)_{i_\ell}\\
&=&\CPP(\Est(\tF_i(X)))\ot \tf_{i_\ell}^{\vep^*_{i_\ell}(\tF_i(X))}(0)_{i_\ell}
=\CP(\tF_{i_\ell}(X)).
\eneqn
The condition $\de_{i_\ell}(X)>0$ is equivalent to the condition
$\varphi_{i_\ell}\bl\CPP(\Est(X)\br>\vep_{i_\ell}(\tf_{i_\ell}^{\vep^*_i(X)}(0)_{i_\ell})$. Hence, 
under the condition $\de_{i_\ell}(X)>0$,
Main Theorem~\ref{th:Main}\;\eqref{MT:4} implies that
\eqn
\tf_{i_\ell}\CP(X)&=&\tf_{i_\ell}(\CPP(\Est(X))\ot \tf_{i_\ell}^{\vep^*_{i_\ell}(X)}(0)_{i_\ell})
=\tf_{i_\ell}\CPP(\Est(X))\ot \tf_{i_\ell}^{\vep^*_{i_\ell}(X)}(0)_{i_\ell}\\
&=&\CPP(\Est(\tF_{i_\ell}(X)))\ot \tf_{i_\ell}^{\vep^*_{i_\ell}(\tF_i(X))}(0)_{i_\ell}
=\CP(\tF_{i_\ell}(X)).
\eneqn

Finally, since $\tE_i$ is an inverse of $\tF_i$ on $\Irr(\tCw)$ and $\te_i$ is an inverse of 
$\tf_i$ on $\Bw$,  $\CP \bl\tF_i(X)\br=\tf_i\CP(X)$ implies that
$\CP \tE_i(X)=\te_i\CP(X)$.
\QED
\Cor \label{cor:epsCw}
For any $i\in I_w$ and any simple $M\in\CBw$, we have
$$\eps_i(M)=\max\bl\eps_i\bl\Qt(M)\br,0\br.$$
\encor

\section{Proof of Main Theorem}\label{sec:proof}
In this section, we give the proof of Main Theorem~\ref{th:Main}.
Recall that $w'\seteq ws_{i_\ell}$.

\subsection{Preparation}
Note that $\Qti[i_\ell]$ is a root object or invertible by Proposition~\ref{prop:simpleroot}.

In order to prove Main Theorem, we may assume that $i\in I_w$
without loss of generality.
If $i\not\in I_{w'}$,
then $\eps_i=0$ and $\tF_i=0$ on $\Irr(\tCw[w'])$ and
$\Qti$ is invertible by Lemma~\ref{lem:Iw'}, and hence Main Theorem is obvious.
Hence we may assume that
\eq
i\in I_{w'}.
\eneq

First note that Main Theorem~\eqref{MT:4} is an  immediate consequence
of Lemma~\ref{lem:ij*} and Proposition~\ref{prop:root}.

\medskip
For $\La\in\pwtl$, set
\eqn
\dC_\La&&\seteq\dM(w\La,\La)\in\Cw,\\
\dC_\La'&&\seteq\dC_\La^{w'}\in\Cw[w']\subset\Cw,\\
\dC''_\La&&\seteq\Esm[i_\ell]\bl\dM(w\La,\La)\br\simeq\dM(w\La,s_{i_\ell}\La)
\in\CBw[w']\cap\Cw\qtq\\
\tdC_\La&&\seteq\Qt(\dC_\La)\in\tCw.
\eneqn
Then, by Corollary~\ref{cor:CBE} we have
$\Qt[w'](\dC''_\La)\simeq\tdC^{w'}_{s_{i_\ell}\La}\simeq\Est(\tdC_\La)$.

By Lemma~\ref{lem:Lacenter}, we have
\eqn\eps_i(\tdC_\la)=\ang{h_i,w\la}
=\eps_i(\tdC^{w'}_{s_{i_\ell}\la})\qt{for any $\la\in\wtl$.}
\eneqn
Hence, we have
\eq&&\eps_i\bl\Est(\tdC_\la)\br=\eps_i(\tdC_\la)\qt{for any $\la\in\wtl$.}
\label{eq:epsE}\eneq

\Prop\label{prop:espww'}
Let $M$ be a simple module in $\Cw$ such that $\eps^*_\il(M)=0$.
Then we have
$$\eps_i\bl\Qt[w'](M)\br=\eps_i\bl\Qt(M)\br\qt{for any $i\in I_{w'}$.}$$
\enprop
\Proof
By Proposition~\ref{prop:MCw}, we can take $\La\in\pwtl$ such that $\dC_\La\hconv \ang{i}\in\Cw$ and 
$\ang{h_\il,\La}=0$.
Then, we have $\dC_\La=\dC'_\La\in\Cw[w']\subset\Cw$ and
$\dC_\La\hconv \ang{i}\in\CBw[w']\subset\CBw$ 
by Proposition~\ref{prop:CB}. Note that $\ang{i}\in\CBw[w']$. 
Let $S$ be a simple quotient of
$(\dC_\La\hconv \ang{i})\conv M$.
Since $\dC_\La\hconv \ang{i}\in\Cw$ and $M\in\CBw$, we have
$S\in\CBw$ by Lemma~\ref{lem:MNla}.

Let us first show
\eq S\in\CBw[w'].
\label{eq:SBw}
\eneq

If $\eps^*_\il(S)=0$
then $S\in\CBw[w']$ by Lemma ~\ref{lem:CBww'}.

Assume that $\eps^*_\il(S)>0$.
Since $\eps^*_\il(\dC'_\La\hconv \ang{i})=\delta_{i,\il}$ and $\eps_\il^*(M)=0$, we have
$\eps^*_\il(S)=1$, $i=\il$ and
$\Es[\il]\bl (\dC_\La\hconv \ang{i})\conv M\br\simeq\dC_\La\conv M$
by Lemma~\ref{lem:E*}.
Applying $\Es$ to $(\dC_\La\hconv \ang{i})\conv M\epito S$, we obtain
$$\dC_\La\conv M\isoto\Es[\il](S),$$
since $M\in\Cw$ implies that $\dC_\La\conv M$ is a simple module.

Hence we obtain
$$S\simeq (\dC_\La\conv M)\hconv\ang{\il},$$
which gives an epimorphism (note that $i=\il$)
$$g\cl (\dC_\La\hconv \ang{i})\conv M
\epito (\dC_\La\conv M)\hconv\ang{i}.$$
Applying $\Qt$ to $g$, we obtain 
$$\tdC_\La\conv \Qt(\ang{i})\conv\Qt(M)\epito
\Qt\bl(\dC_\La\conv M)\hconv\ang{i}\br
\simeq
\tdC_\La\conv\bl\Qt(M)\hconv\Qt(\ang{i})\br.$$
Here the last isomorphism follows from Proposition~\ref{prop:CB} and
$\dC_\La\conv M\in\Cw$, $\ang{i}\in\CBw$ and the fact that
$\Qt(\ang{i})$ is \afr.

Hence we have
$$\Qt(\ang{i})\hconv \Qt(M)
\simeq \Qt(M)\hconv\Qt(\ang{i}),$$
which implies that $\Qt(\ang{i})$ and $\Qt(M)$ commute,
and $\Qt(\ang{i})\conv\Qt(M)$ is simple.
Note that $\Qt\vert_{\Cw}\cl\Cw\to\tCw$ is fully faithful,
and $\Qt\bl(\dC_\La\hconv \ang{i})\conv M\br$ is simple.
Since $\dC_\La\hconv \ang{i}$ and $M$ belong to $\Cw$, we conclude that
$(\dC_\La\hconv\ang{i})\conv M$ is simple.
Hence we have $S\simeq (\dC_\La\hconv\ang{i})\conv M$.
Since $\dC_\La\hconv\ang{i}$ and $M$ belong to $\CBw[w']$ and
$$\Qt[w'](S)\simeq\Qt[w']\bl(\dC_\La\hconv\ang{i})\conv M\br
\simeq\Qt[w'](\dC_\La\hconv\ang{i})\conv \Qt[w'](M)$$ does not vanish,
which means that $S\in\CBw[w']$.

Thus
we have obtained \eqref{eq:SBw}.

\mnoi

Let $f\cl (\dC_\La\hconv \ang{i})\conv M
\to M\conv(\dC_\La\hconv \ang{i})$ be a morphism such that $\Im(f)\simeq S$
(see Lemma~\ref{lem:Laspl}). Since $S$ belongs to
$\CBw[w']\subset\CBw$,
Lemma~\ref{lem:wLadf} implies that
$$\La\bl\Qt[w'](\dC_\La\hconv \ang{i}),\Qt[w'](M)\br
=\deg(f)=\La\bl\Qt[w](\dC_\La\hconv \ang{i}),\Qt[w](M)\br.$$

Hence we obtain
$$\eps_i\bl\Qt[w'](M))\br=\eps_i\bl\Qt[w](M)\br.$$
\QED

\Rem
In the course of the proof of Proposition~\ref{prop:espww'}, we don't know if $\dC_\La\hconv\ang{i}$ has an affinization.
If we knew it, the proof could be much simpler.
Conjecturally, any simple module in $\Cw$ has an affinization.
\enrem
Now we divide the proof of Main Theorem into two cases:
$i\not=i_\ell$ and  $i=i_\ell$.

\subsection{Case $i\not=i_\ell$}
Assume that $i\in I_{w'}$ satisfies $i\not=i_\ell$.

\Lemma\label{lem:FiE} For any simple $X\in\tCw$, we have
\bnum
\item $\eps^*_{i_\ell}(\tF_i X)=\eps^*_{i_\ell}(X)$,\label{it;Esti}
\item $\Est(\tF_i(X))\simeq\tF_i(\Est(X))$.
\ee
\enlemma
\Proof
The assertion (i) follows from Lemma~\ref{lem:ij*}\;\eqref{it:ij*}.

\snoi
(ii) 
We may assume that $X=\Qt(M)$ for a simple $M\in\Cw$,
since for any $X\in\tCw$ there exists $\La\in\pwtl$ such that $\wtil{C_\La} \circ X\in\Qt(\Cw)$ and 
we can reduce the assertion to $\wtil{C_\La} \circ X$.

Set $n\seteq\eps^*_\il(M)$ and $M'\seteq \Es[\il]{}^{(n)}(M)$.
Take $\La\in\pwtl$ such that $\dC_\La\hconv\ang{i}\in\Cw$ and $\ang{h_\il,\La}=0$ (see Proposition~\ref{prop:MCw}).
Then Proposition~\ref{prop:CL} implies that
$\eps^*_{i_\ell}(\dC_\La\hconv\ang{i})=\eps^*_{i_\ell}(\dC_\La)+\eps^*_{i_\ell}(\ang{i})=0$.

On the other hand (i) implies that
$$\eps^*_{i_\ell}\bl(\dC_\La\hconv\ang{i})\hconv M)
=\eps^*_{i_\ell}(\dC_\La)+\eps^*_{i_\ell}(M)
=\eps^*_{i_\ell}(\dC_\La\hconv\ang{i})+\eps^*_{i_\ell}(M).$$

Hence Proposition~\ref{prop:E*} implies that
\eq
&&\ba{rl}
\Esm[{i_\ell}]\bl(\dC_\La\hconv\ang{i})\hconv M\br
&\simeq\Esm[{i_\ell}](\dC_\La\hconv\ang{i})\hconv( \Esm[{i_\ell}]M)\\
&\simeq (\dC_{\La}\hconv\ang{i})\hconv M'.
\ea
\label{eq:ELM}
\eneq

Since we have 
\eqn
\Qt[w']\bl \Esm[{i_\ell}]\bl(\dC_\La\hconv\ang{i})\hconv M\br\br
&&\simeq \Est(\tdC_\La\conv\tF_i X)\\
&&\simeq
\Qt[w'](\dC_\La)\conv\Est(\tF_iX)
\eneqn
and
\eqn
\Qt[w']\bl (\dC_{\La}\hconv\ang{i})\hconv (\Esm[{i_\ell}]M)\br
\simeq \Qt[w'](\dC_\La)\conv\tF_i\bl\Est(X)\br,
\eneqn
we obtain
$$\Est(\tF_i(X))\simeq\tF_i(\Est(X)).$$
\QED
The following lemma shows Main Theorem \eqref{MT:1} when $i\ne i_\ell$.
\Lemma For any simple $X\in\tCw$ and any $i\in I_{w}$ such that $i\not=i_\ell$, we have
$$\eps_i(\Est(X))= \eps_i(X).$$
\enlemma
\Proof
We may assume that $X=\Qt(M)$ for a simple $M\in\Cw$.
Set $M'=\tEsm_\il(M)\in\Cw\cap\CBw[w']$.
Since we have
$$\eps_i(\Est(X))=\eps_i\bl\Qt[w'](M')\br=\eps_i\bl\Qt(M')\br.$$
Here the last equality follows from Proposition~\ref{prop:espww'}.
Set $n=\epss_\il(M)$.
Then, since $\Qt(M)\simeq\tFs_i{}^n\Qt(M')$,
we obtain
\[
\eps_i\bl\Qt(M')\br=\eps_i\bl\Qt({M})\br=\eps_i(X).\qquad
\]
Here, the first equality follows from Lemma~\ref{lem:ij*}.
\QED

Thus, we find that Main Theorem~\ref{th:Main} holds for $i\not=i_\ell$.

\subsection{Case $i=i_\ell$}
 Recall that
\eq
\,\,\qquad\de_i(X)=\eps_i(X)+\eps_i^*(X)+\ang{h_i,\wt(X)}
=\eps_i(X)-\bl\eps_i^*(X)-\ang{h_i,\wt(\Est X)}\br,
\label{eq:deeps}
\eneq
where the last equality is from \eqref{eq:e*E}.

\Lemma\label{lem:eps<}
For any simple $X\in\tCw$, we have
$$\eps_i\bl\Est(X)\br\le\eps_i(X).$$
\enlemma
\Proof

By \eqref{eq:epsE}, we may assume that $X=\Qt(M)$ for a simple $M\in\Cw$.
Set $M'=\tEsm_{i}(M)\in\CBw[w']\cap\Cw$.
Then we have $\Est(X)\simeq\Qt[w'](M')$.

Let $m=\eps_i^*(M)$.
Then $M\simeq M'\hconv\ang{i^m}$.
Since $M\in\Cw$, we obtain
$$\Qt(M)\simeq\Qt(M')\hconv \Qt(\ang{i^m})\simeq(\tF^*_i)^m(\Qt(M')).$$
Since $\eps_i\bl\tF_i^*(X)\br=\eps_i(X)$ or $\eps_i(X)+1$
for any simple $X\in\Cw$ by Lemma~\ref{lem:iistar},
we obtain
$$\eps_i\bl\Qt[w](M')\br\le\eps_i\bl\Qt[w](M)\br.$$
Since  Proposition~\ref{prop:espww'}
implies that $\eps_i\bl\Qt[w'](M')\br=\eps_i\bl\Qt[w](M)\br$,
the assertion follows from
$\eps_i(\Est X)=\eps_i\bl\Qt[w'](M')\br$
and
$\eps_i(X)=\eps_i\bl\Qt[w](M)\br$.
\QED

\subsubsection{{\rm Case:} $\de_i(X)=0$.}
In this case, since $\Qt(\ang{i})$ and $X$ commute with each other, we have 
$\tF_i(X)=\tF^*_i(X)$.
Hence we obtain
$\eps_i^*(\tF_i X)=\eps_i^*(\tF^*_i X)=\eps_i^*(X)+1$,
and Lemma~\ref{lem:EFX=EX}\;(iii)
implies that $\Est(\tF_i X)=\Est (X)$.
Main Theorem\;\eqref{MT:1}
follows from Lemma~\ref{lem:eps<} and $\eps_i(X)=-\eps_\il(X)-\ang{h_i,\wt(X)}$.

Thus Main Theorem holds in this case.

\subsubsection{{\rm Case:}\;$\de_i(X)>0$.}
In this case,
$\Qti$ is a root object by
Proposition~\ref{prop:simpleroot} and Lemma~\ref{lem:inv}.

Let us show that
\eq
\Est(\tF_i X)\simeq\tF_i\bl\Est X\br\qtq\eps_i(X)=\eps_i(\Est X).
\label{eq:red}
\eneq
They, along with \eqref{eq:deeps}, imply that
$$\eps_i(\Est X)-\eps^*_i(X)+\ang{h_i,\wt (\Est\,X)}=\de_i(X)>0,$$
and we obtain Main Theorem.

\medskip
In order to see \eqref{eq:red}, we may assume that $X=\Qt(M)$ for a simple $M\in\Cw$.

Take $\La\in\pwtl$ such that $\dC_\La\hconv \ang{i}\in\Cw$
and $\ang{h_i,\La}=0$.
Then Proposition~\ref{prop:shsw} implies that
$(\dC_\La\hconv \ang{i},M)$ is a $\La$-definable pair.
Then Main Theorem~\eqref{MT:4} 
implies that
\eq
\eps_i^*\bl(\dC_\La\hconv \ang{i})\hconv M\br=\eps_i^*(M)
+\eps_i^*(\dC_\La)=\eps_i^*(M).\label{eq:epsssum}
\eneq
Set $m=\eps_i^*(M)$.

Then we have an exact sequence by Lemma~\ref{lem:E*}:
\eqn
0\To \bl\dC_\La\hconv \ang{i}\br\conv \Es^{(m)}M&&\To
\Es^{(m)}\bl(\dC_\La\hconv \ang{i})\conv M\br\\
&&\hs{7ex}\To\Es\bl\dC_\La\hconv \ang{i}\br\conv \Es^{(m-1)}M\To0.
\eneqn
Note that
$\Es\bl\dC_\La\hconv \ang{i}\br\simeq\dC_\La$.
If the composition $$f\cl\bl\dC_\La\hconv \ang{i}\br\conv \Es^{(m)}M\to
\Es^{(m)}\bl(\dC_\La\hconv \ang{i})\conv M\br\epito \Es^{(m)}\bl(\dC_\La\hconv \ang{i})\hconv M\br$$ vanishes, then we have an epimorphism
\eqn
\dC_\La%
\conv \Es^{(m-1)}M&&\simeq
\Es\bl\dC_\La\hconv \ang{i}\br\conv \Es^{(m-1)}M\\
&&\epito
\Es^{(m)}(\dC_\La\hconv \ang{i})\hconv M\br\simeq
\tEsm_i\bl(\dC_\La\hconv \ang{i})\hconv M\br,\eneqn
where the last isomorphism follows from \eqref{eq:epsssum}.

Hence it induces a non-zero $R(\beta)$-module homomorphism:
$$ \Es^{(m-1)}M\to Y\seteq\HOM\Bigr(\dC_\La''\conv R(\beta), \tEsm_i\bl(\dC_\La\hconv \ang{i})\hconv M\br\Bigl),$$
where $\beta=-\wt\bl\Es^{(m-1)}M\br\in\prtl$.
Since $\Es^{(m-1)}M$ has a simple head $(\tEs_{i}){}^{m-1}(M)$ with
$\eps^*_i\bl(\tEs_{i}){}^{m-1}(M)\br=1$,
and
$\Es(Y)=0$, it is a contradiction.

\medskip
Hence $f$ does not vanish and we obtain
\eq
(\dC_\La\hconv\ang{i})\hconv \Es^{(m)}M
\simeq\tEsm_i\bl(\dC_\La\hconv \ang{i})\hconv M\br.
\label{eq:La-i-M}
\eneq
By\eqref{eq:La-i-M}, Proposition~\ref{prop:CL} and Corollary~\ref{cor:CBE},
we have 
\eq&&\ba{rl}
\Est \bl\Qt((\dC_\La\hconv \ang{i})\hconv M)\br&\simeq
\Qt[w']\Bigl(\tEsm_i\bl(\dC_\La\hconv \ang{i})\hconv M\br\Bigr)\\
&\simeq
\hd\Bigl(\Qt[w']\bl(\dC_\La\hconv\ang{i})\conv \Es^{(m)}M\br\Bigr)\\
&\simeq 
\Qt[w']\bl\dC_\La\hconv\ang{i}\br\hconv\Est(\Qt M).\ea\label{eq:tFires}
\eneq
Hence, we obtain
$\Est(\tF_iX)\simeq\tF_i(\Est X)$, which shows Main Theorem~\ref{th:Main} \eqref{MT:2}.

\medskip
It remains to prove that
\eq\eps_i(X)=\eps_i(\Est X).\label{eq:eps=}
\eneq
The proof of \eqref{eq:eps=} is similar to the one of Lemma~\ref{lem:eps<}.
We may assume that $i\in I_{w'}$.
By {\eqref{La-M} and }\eqref{eq:epsE}, we may assume that $X=\Qt(M)$ for a simple $M\in\Cw$.
Set $M'=\tEsm_{i}(M)\in\CBw[w']\cap\Cw$.
Then we have $\Est(X)\simeq\Qt[w'](M')$ {by Corollary~\ref{cor:CBE}}.
Let $m=\eps_i^*(M)$.
Set $X'=\Qt(M')$.
Then, we have
$X'=\tEs_i{}^m(X)$.
Since $\eps_i(\Est(X))=\eps_i(X')$ by Proposition~\ref{prop:espww'}
and $\eps_i(X')=\eps_i(X)$ follows from Lemma~\ref{lem:iistar} (i)
and $\de_i(X)>0$.
Thus we obtain $\eps_i\bl\Est(X)\br=\eps_i(X)$.

Now, the proof of Main Theorem~\ref{th:Main} is complete.

\section{Bijectivity and connectedness}
\subsection{Bijectivity of $\CP$}
Let $\uw=s_{i_1}\cdots s_{i_\ell}$ be a reduced expression of $w\in \weyl$.

In this section, we shall prove the following proposition.

\Prop\label{prop:bij}
 $\CP\cl\Irr(\tCw)\to\Bw$ is bijective.
\enprop

By induction on the length of $w$
and Proposition~\ref{prop:CPww'}, we can reduce this proposition to the following proposition.

\Prop
Let $w\in \weyl$ and $i\in I$ such that $w'\seteq ws_i<w$.
Then the map
$$\CK\cl \Irr(\tCw)\to \Irr(\tCw[w'])\times\Z
\qt{given by $X\mapsto \bl\Est(X),\eps^*_i(X)\br$},$$
is bijective
\enprop

\Proof 
Injectivity is obvious.
Let us prove the surjectivity.

We have $\Est\bl(\tF^*_i)^nZ\br\simeq\Est(Z)$ by Lemma~\ref{lem:EFX=EX} 
and $\eps_i^*\bl(\tF^*_i)^nZ\br
=\eps_i^*(Z)+n$ for any $n\in\Z$ and any simple $Z\in\tCw$.
Here, if $n<0$, then $(\tF^*_i)^n$ means $(\tE^*_i)^{-n}$ 
by Proposition~\ref{pro:inverse}.
Hence, it is enough to show that
$$\text{$\Est\cl\Irr(\tCw)\to\Irr(\tCw[w'])$ is surjective.}$$

In order to prove this, let us recall results in \cite[Theorem~4.8]{KKOP23}):
\bnum
\item
there exists a quasi-commutative diagram:
$$\xymatrix{
&\Cw\ar[dl]_-{\Qt[w']\vert_{\Cw}}\ar[dr]^{\Qr}\\
\tCw[w']\ar[rr]^\sim_\Phi&&\tCw[w,s_i]
}
$$
where $\Qr\cl\Cw\to\tCw[w,s_i]$ is the localization functor 
by the commuting family of right braiders 
$\st{\dM(w\La,s_i\La)\mid\La\in\pwtl}$ in $\Cw$,
and $\Phi\cl\tCw[w']\isoto\tCw[w,s_i]$ is an equivalence of monoidal categories,
\label{it:p1}
\item for any simple $Z\in\Cw$ such that $\Qr(Z)\not\simeq0$,
there exists $\La\in\pwtl$ such that
\bna
\item
$Z\hconv\dM(w\La,s_i\La)\in\catC_{w,s_i}$,
\item $\Qr\bl Z\conv\dM(w\La,s_i\La)\br\to
\Qr\bl Z\hconv\dM(w\La,s_i\La)\br\to
\Qr\bl \dM(w\La,s_i\La)\conv Z\br$
are isomorphisms.
\ee
\label{it:p2}
\ee
Note 
that $\catC_{w,s_i}$ is the full subcategory of $\Cw$
consisting of objects $M\in\Cw$
such that $\Es(M)\simeq0$.

Now let us show the surjectivity of $\Est\cl\Irr(\tCw)\to\Irr(\tCw[w'])$.
Recall  $\tdC^{w'}_\la$ as in \eqref{eq:center}. 
Since Proposition~\ref{prop:ELM} implies that 
\[
\Est(\Qt(M)\conv\tdC_{s_i\la})
\simeq\Est(\Qt(M))\conv\tdC^{w'}_\la
\]
 for any simple $M\in\Cw$ and $\la\in\wtl$,
it is enough to show that
\eq\hs{3ex}&&\text{for any simple $M\in\Cw[w']$, there exists
$X\in\tCw$ such that $\Est(X)\simeq\Qt[w'](M)$.}\label{eq:sur}
\eneq

Now, since $\Qt[w'](M)$ is simple,
$\Qr(M)$ is simple by \eqref{it:p1}. Hence, 
\eqref{it:p2} implies that there exists $\La\in\pwtl$ such that
$N\seteq M\hconv \dM(w\La,s_i\La)\in\catC_{w,s_i}$ and
$\Qr(N)\isoto\Qr(\dM(w\La,s_i\La)\conv M)$.
Hence, we have
\eqn
\Qt[w'](N)&&\simeq\Qt[w']\bl\dM(w\La,s_i\La)\conv M\br
\simeq\tdC^{w'}_{s_i\La}\conv\Qt[w'](M)
\eneqn
by \eqref{it:p2}. Hence, we obtain
$$\Qt[w'](M)\simeq \Qt[w'](N)\conv\tdC^{w'}_{-s_i\La}
\simeq \Est\bl\Qt(N)\conv\tdC_{-\La}\br.$$
Thus we obtain \eqref{eq:sur}.
\QED
Therefore, by Main Theorem~\ref{th:Main} and Proposition~\ref{prop:bij} we obtain 
\MTh\label{thm:main2}
$\CP\cl\Irr(\tCw)\to\Bw$ is an isomorphism of crystals.
\enmth

\subsection{Connectedness of $\Bw$}
As an application of Main Theorem~\ref{thm:main2} in previous subsections, let us show the connectedness of 
the crystal $\Bw$ as a crystal graph (Definition~\ref{def:graph}). 
\Th
The crystal $\Bw$ is connected.
\enth
Note that in \cite{Kana-N}, it is shown that $\Bw$ is connected for a semi-simple $\g$.
\Proof
Since $\Bw\simeq\Irr(\tCw)$, we shall show that
$\Irr(\tCw)$ is connected.

Since by Theorem~\ref{th:main1}, we have the crystal isomorphism 
\[
{\psi_*}\cl\bl\Irr(\tCw),\st{\tE_i,\tF_i}_{i\in I}\br\simeq
\bl\Irr(\tCw[w^{-1}],\st{\tEs_i,\tFs_i}_{i\in I}\br),
\]
the assertion is a consequence of the following lemma.
\QED

\Lemma
The crystal $\Irr(\tCw)$ with respect to the crystal operators
$\{\tFs_i,\tEs_i\}_{i\in I}$ is connected.
\enlemma
\Proof
Since as we have seen in Lemma~\ref{lm:stable} that
$\Irr(\Cw)\sqcup\{0\}$ is stable by $\tEs_i$,  
it follows from Lemma~\ref{lem:tEcomp} and  $\tEsm_{i_1}\cdots\tEsm_{i_\ell}M\simeq\one$
that
$\Qt(M)$ is connected to $\one$ for any $M\in\Irr(\Cw)$.
Since $\tEs_i(\tdC_\la\conv X)\simeq\tdC_\la\conv \tEs_i(X)$ for any $X\in\Irr(\tCw)$ and $\la\in\wtl$, an object
$\tdC_{-\La}\conv\Qt(M)\in \Irr(\tCw)$ ($\La\in\pwtl$, $M\in\Irr(\Cw)$)
is connected to
$\tdC_{-\La}\conv\one$,
which is connected with
$\tdC_{-\La}\conv\tdC_{\La}\simeq\one$ by applying the above argument to
$\tdC_{\La}=\Qt(\dC_\La)\in \Qt(\Irr(\Cw))$ (see \eqref{det-in-C}).
\QED

Here, we give an example of $\tCw$, as a monoidal categorification of
a cluster algebra:
\Ex[{cf.\ }\cite{KKKO18}]
Let $\g=A_3$, $w=s_2w_0=s_1s_2s_3s_2s_1$,
$\Cw=\st{M\in R\gmod\mid \E_2M\simeq0}$.
$\dC_{\La_1}=\ang{321}$, $\dC_{\La_2}=\ang{132}$, $\dC_{\La_3}=\ang{123}$ 
are frozen variables.
In this case, $\tCw$ is a monoidal categorification of the cluster algebra
with four clusters: 
$$\xymatrix{\ang{1}\ar@{-}[d]_{C_1}\ar@{-}[r]^{C_4}&\ang{12}\ar@{-}[d]^{C_3}\\
  \ang{3}\ar@{-}[r]^{C_2}&\ang{32},}$$ 
\scalebox{.95}{\parbox{\textwidth}{
\eqn
\CP(\ang{1}^x\conv\ang{3}^y\conv\ang{321}^a\conv\ang{132}^b\conv\ang{123}^c)
&&=(c+b,c,y+a+b+c,a+b,x+a),\\
\CP(\ang{3}^y\conv\ang{32}^x\conv\ang{321}^a\conv\ang{132}^b\conv\ang{123}^c)
&&=(c+b,c,y+x+a+b+c,x+a+b,a),\\
\CP(\ang{32}^x\conv\ang{12}^y\conv\ang{321}^a\conv\ang{132}^b\conv\ang{123}^c)
&&=(y+c+b,c,x+a+b+c,x+y+a+b,a),\\
\CP(\ang{12}^y\conv\ang{1}^x\conv\ang{321}^a\conv\ang{132}^b\conv\ang{123}^c)
&&=(y+c+b,c,a+b+c,y+a+b,a+x).
\eneqn}}

Note that the frozen variables correspond to the following elements in $\Bw$:
\eqn
&&\dC_{\La_1}=\ang{321}=\tF_3\tF_2\tF_1\one\longleftrightarrow \til f_3\til f_2\til f_1(0,0,0,0,0)=(0,0,1,1,1),\\
&&\dC_{\La_2}=\ang{132}=\tF_1\tF_3\tF_2\one\longleftrightarrow \til f_1\til f_3\til f_2(0,0,0,0,0)=(1,0,1,1,0),\\
&&\dC_{\La_3}=\ang{123}=\tF_1\tF_2\tF_3\one\longleftrightarrow \til f_1\til f_2\til f_3(0,0,0,0,0)=(1,1,1,0,0).
\eneqn

Here $x,y\in\Z_{\ge0}$ and $a,b,c\in\Z$.
Each cluster $C_k$ ($1\le k\le 4$) consists of a commuting family
of simple modules. 

Hence, the images of the four clusters  are given as follows:
\eqn
\ba{lll}
\CP(C_1)&=\{(c_1,c_2,c_3,c_4,c_5)\in \Bw
\mid c_3\ge c_2+c_4,\hs{1ex} &c_1+c_5\ge c_2+c_4\},\\
\CP(C_2)&=\{(c_1,c_2,c_3,c_4,c_5)\in \Bw
\mid c_3\ge c_2+c_4,\hs{1ex} &c_1+c_5\le c_2+c_4\},\\
\CP(C_3)&=\{(c_1,c_2,c_3,c_4,c_5)\in \Bw
\mid c_3\le c_2+c_4,\hs{1ex} &c_1+c_5\le c_2+c_4\},\\
\CP(C_4)&=\{(c_1,c_2,c_3,c_4,c_5)\in \Bw
\mid c_3\le c_2+c_4,\hs{1ex} &c_1+c_5\ge c_2+c_4\}.\\
\ea
\eneqn
$$\Bw\qquad\xymatrix@C=-2ex@R=1ex{
&&&&\raisebox{1.5ex}{$c_3-c_2-c_4$}\\
&&&C_2&&C_1\\
\ar[rrrrrrr]&\akew[5ex]&&&&&\akew[5ex]&\akew[2ex]c_1+c_5-c_2-c_4\\
&&&C_3&&C_4\\
&&&&\ar[uuuu]
}$$
\enex

\bibliographystyle{amsplain}
\providecommand{\bysame}{\leavevmode\hbox to3em{\hrulefill}\thinspace}
\providecommand{\MR}{\relax\ifhmode\unskip\space\fi MR }
\providecommand{\MRhref}[2]{%
  \href{http://www.ams.org/mathscinet-getitem?mr=#1}{#2}
}
\providecommand{\href}[2]{#2}

\end{document}

\bibitem{VV09}
M. Varagnolo and E. Vasserot,
 \emph{Canonical bases and KLR algebras},
J. Reine Angew. Math. \textbf{659} (2011), 67--100.

\end{thebibliography}

\end{document}

\end{thebibliography}
\end{document}